\documentclass[10pt]{amsart}
\usepackage{amssymb,amsthm,amsmath}
\usepackage{graphicx}
\usepackage{color}
\usepackage{array}
\usepackage{appendix}
\usepackage{url}

 \numberwithin{equation}{section}
\newtheorem{thm}{Theorem}[section]
\newtheorem{cor}[thm]{Corollary}
\newtheorem{lem}[thm]{Lemma}
\newtheorem{prop}[thm]{Proposition}
\newtheorem{defn}[thm]{Definition}

\newtheorem{prob}[thm]{Problem}
\theoremstyle{definition}
\newtheorem{rmk}[thm]{Remark}

\newenvironment{rem}{%
\bigskip
\noindent \textsl{{\sl Remark. }}}{\bigskip}

\newcommand{\ignore}[1]{}

\newcommand{\ip}[2]{{\langle#1,#2\rangle}}

\newcommand{\norm}[1]{{\|#1\|}}

\newcommand{\DDD}{{\bf D}}

\def \R{\mathbb{R}}
\def \C{\mathbb{C}}

\newcommand{\fc}{{\mathcal{F}}}

\newcommand{\Acr}{{\mathcal{A}}}

\newcommand{\Pp}{{\mathcal P}}

\newcommand{\Ssc}{\mathfrak{S}}

\renewcommand{\SS}{\mathcal{S}}

\newcommand{\Rnd}{{\R}^{n\times d}}
\newcommand{\Sn}{{\mathcal S}_n}
\newcommand{\hRnd}{{\widehat{\Rnd}}}

\newcommand{\halpha}{{\hat{\alpha}}}
\newcommand{\hbeta}{{\hat{\beta}}}
\newcommand{\Ddelta}{\mathbb{B}}

\newcommand{\XX}{X}
\newcommand{\YY}{Y}

\newcommand{\XXb}{{\bf x}}
\newcommand{\TTb}{{\bf t}}

\renewcommand{\Mc}{{\mathcal{M}_b}}
\newcommand{\Xinf}{{X^{\infty}}}
\newcommand{\Yinf}{{Y^{\infty}}}

\title[Permutation Invariant Representations]{Permutation Invariant Representations with Applications to Graph Deep Learning}
\author{Radu Balan}
\author{Naveed Haghani}
\author{Maneesh Singh}

\address{Department of Mathematics, University of Maryland, College Park, MD 20742, USA}
\email{rvbalan@umd.edu}
\address{Applied Mathematics and Statistics and Scientific Computation Program, University of Maryland, College Park, MD 20742, USA}
\email{nhaghan1@umd.edu}
\address{Verisk Analytics,  Jersey City, NJ 07310, USA}
\email{dr.maneesh.singh@ieee.org}

\begin{document}

\maketitle

\begin{abstract}
This paper presents primarily two Euclidean embeddings of the quotient space 
generated by matrices that are identified modulo arbitrary row permutations. 
The original application is in deep learning on graphs where the learning task 
is invariant to node relabeling. Two embedding schemes are introduced, one based on 
sorting and the other based on algebras of multivariate polynomials. While both 
embeddings exhibit a computational complexity exponential in problem size, the sorting based embedding is globally bi-Lipschitz and admits a low dimensional target space. Additionally, an almost everywhere injective scheme can be implemented with minimal redundancy and low computational cost. In turn, this proves that almost any classifier can be implemented with an arbitrary small loss of performance. 
Numerical experiments are carried out on two data sets,   a chemical compound  data set (QM9) and a proteins 
data set (PROTEINS\_FULL). 
\end{abstract}
    
\section{Introduction\label{sec1}}

This paper is motivated by a class of problems in graph deep learning, where the
primary task is either graph classification or graph regression. 
In either case, the result should be invariant to arbitrary permutations of graph nodes.

As we explain below, the mathematical problem analyzed in this paper is a special case 
of the permutation invariance issue described above. To set the notations consider the
vector space $\Rnd$ of $n\times d$ matrices endowed with the Frobenius norm 
 $\norm{X}=\left(trace(XX^T)\right)^{1/2}$
and its associated Hilbert-Schmidt scalar product, $\ip{X}{Y}=trace(XY^T)$.
 Let $\Sn$ denote the symmetric group of $n\times n$ permutation matrices. 
 $\Sn$ is a finite group of size $|\Sn|=n!$.

On $\Rnd$ we consider the equivalence relation $\sim$ 
induced by the symmetric group of permutation matrices $\Sn$ as follows. Let $X,Y\in\Rnd$. 
Then we say $X\sim Y$ if there is $P\in\Sn$ so that $Y=PX$. In other words, two matrices are equivalent if one is a row permutation of the other. 
The equivalence relation induces a natural distance on the quotient space 
$\hRnd:=\Rnd/\sim$,
\begin{equation}
    \label{eq:1.1}
d: \hRnd\times \hRnd \rightarrow\R  ~~,~~d(\hat{X},\hat{Y})=\min_{\Pi\in\Sn}\norm{X-\Pi Y}   
\end{equation}
This makes $(\hRnd,d)$ a complete metric space.

Our main problem can now be stated as follows:
\begin{prob}\label{prob1}
Given $n,d\geq 1$ positive integers, find $m$ and a bi-Lipschitz map
 $\hat{\alpha}:(\hRnd,d)\rightarrow(\R^m,\norm{\cdot}_2)$.
\end{prob}
Explicitly the problem can be restated as follows. One is asked to construct a 
map $\alpha:\Rnd\rightarrow\R^m$ that satisfies the following conditions:
\begin{enumerate}
    \item If $X,Y\in\Rnd$ so that $X\sim Y$ then $\alpha(X)=\alpha(Y)$
    \item If $X,Y\in\Rnd$ so that $\alpha(X)=\alpha(Y)$ then $X\sim Y$
    \item There are constants $0<a_0\leq b_0$ so that for any $X,Y\in\Rnd$,
\begin{equation}
    \label{eq:1.2}
    a_0 \min_{\Pi\in\Sn}\norm{X-\Pi Y}\leq \norm{\alpha(X)-\alpha(Y)}_2 \leq
    b_0 \min_{\Pi\in\Sn}\norm{X-\Pi Y}
\end{equation}
\end{enumerate}
Condition (1) allows us to lift $\alpha$ to the quotient space $\hRnd$. 
Thus $\hat{\alpha}(\hat{X}) = \alpha(X)$ is well-defined. Condition (2) says that $\hat{\alpha}$ is injective (or, that $\alpha$ is faithful with respect to the equivalence relation $\sim$).
Condition (3) says that $\hat{\alpha}$ is bi-Lipschitz with constants $a_0,b_0$.
By a slight abuse of notation, when $\alpha$ satisfies (1) 
we shall use the same letter to denote
the map $\alpha:\Rnd\rightarrow\R^m$ as well as the induced map on the quotient space
 $\alpha:\hRnd\rightarrow\R^m$. 

For $X,Y\in\Rnd$, $d(X,Y)$ denotes the same quantity in (\ref{eq:1.1}) .
In this case $d$ is only a semi-distance on $\Rnd$, i.e., 
it is symmetric, non-negative and satisfies the triangle inequality but fails the 
positivity condition. 

One approach to embedding $\hRnd$ is to consider the convex set of probability measures on $\R^d$, $\Pp(\R^d)$, and the map
\begin{equation}\label{eq:measure} \alpha_\infty:\Rnd\rightarrow \Pp(\R^d)~~,~~\alpha_\infty(X) = \frac{1}{n}\sum_{k=1}^n \delta(\cdot-x_k) 
\end{equation}
where $[x_1,\ldots,x_n]=X^T$, i.e., $x_k$ is the $k^{th}$ row of $X$ reshaped as a vector, and $\delta$ denotes the Dirac measure. 
When $\Pp(\R^d)$ is endowed with the Wasserstein-1 distance (the Earth Moving Distance), known also as the Kantorovich-Rubinstein metric, 
\[ d_{KR}(p,q) = \inf_{\begin{array}{c}
\mbox{$\pi\in\Pp(\R^d\times\R^d):$} \\ \mbox{$\pi(\cdot,\R^d)=p$} \\
\mbox{$\pi(\R^d,\cdot)=q$}
\end{array}} \int_{\R^d \times\R^d} \norm{x-y}d\pi(x,y) \] 
the distance between $a_{\infty}(X)$
and $a_{\infty}(Y)$ becomes
\[ d_{KR}(a_{\infty}(X),a_{\infty}(Y)) =
\min_{\Pi\in \Sn} \sum_{k=1}^n \norm{x_k-(\Pi Y)_k}. \]
By the Kantorovich-Rubinstein theorem (\cite{Villani}Theorem 1.14), $d_{KR}$ extends to a norm on the linear space of bounded signed Borel measures on $\R^d$, $\Mc(\R^d)$. It is easy to verify that
\[ d(\hat{X},\hat{Y}) \leq
d_{KR}(a_{\infty}(X),a_{\infty}(Y)) \leq \sqrt{n} d(\hat{X},\hat{Y}) \]
which proves that $a_{\infty}$ provides an embedding into a normed linear space. Yet this embedding does not solve the problem since the linear space $\Mc(\R^d)$ is infinite dimensional. 

Instead of the previous infinite dimensional embedding, we consider two different classes of  embeddings. To illustrate these two constructions, consider the simplest case $d=1$. 
\begin{enumerate}
    \item {\em Algebraic Embedding}. For $x\in\R^n$, $x=(x_1,\ldots,x_n)^T$, 
    construct the polynomial $P_x(z)=(z-x_1)\cdots(z-x_n)$ and then expand the product: $P_x(z)=z^n + c_1(x)z^{n-1}+\cdots+c_n(x)$. Using Vieta's formulas and Newton-Girard identities, an algebraically equivalent description of $P_x$ is given by the symmetric polynomials:
\begin{equation}\label{eq:poly} 
\alpha:\R^n\rightarrow\R^n~~,~~
    \alpha(x)=\left( \sum_{k=1}^n x_k , 
    \sum_{k=1}^n x_k^2,\ldots,\sum_{k=1}^n x_k^n
    \right). 
\end{equation}
    It is not hard to see that this map satisfies Conditions (1) and (2) and therefore lifts to an injective continuous map $\hat{\alpha}$ on $\hat{\R^n}$. Yet it is not Lipschitz, let alone bi-Lipschitz. The approach in \cite{Cahill19} can be used to modify $\alpha$ to a Lipschitz continuous map, but, for the same reason as described in that paper, it cannot be ``fixed" to a bi-Lipschitz embedding. In Section \ref{sec2} we show how to construct an algebraic Lipschitz embedding in the case $d>1$.  
    \item {\em Sorting Embedding}.
    For $x\in\R^n$, consider the sorting map
\begin{equation}\label{eq:ord}
    \downarrow:\R^n\rightarrow\R^n ~~,~~
    \downarrow(x)=(x_{\pi(1)},x_{\pi(2)},\ldots,
    x_{\pi(n)})^T 
\end{equation}
    where the permutation $\pi$ is so that
     $x_{\pi(1)}\geq x_{\pi(2)}\geq\cdots\geq x_{\pi(n)}$. It is obvious that $\downarrow$ satisfies Conditions (1) and (2) and therefore lifts to an injective map on $\hRnd$. As we see in Section \ref{sec3}, the map $\downarrow$ is bi-Lipschitz. In fact it is isometric, and hence produces an ideal embedding. Our work in Section \ref{sec3} is to extend such construction to the more general 
     case $d>1$.
\end{enumerate}
The algebraic embedding is a special case of the more general {\em kernel method} that can be thought of as a projection of the measure 
$a_{\infty}(X)$ onto a finite dimensional space, e.g., the space of polynomials spanned by $\{X,X^2,\cdots,X^n\}$. In applications such kernel method is known as a ``Readout Map" \cite{deepsets}, based on ``Sum Pooling".

The sorting embedding has been used in applications under the name of ``Pooling Map" \cite{deepsets}, based on ``Max Pooling". A na\"{\i}ve extension of the unidimensional map (\ref{eq:ord}) to the case $d>1$ might employ the lexicographic order: order monotone decreasing the rows according to the first column, and break the tie by going to the next column. While this gives rise to an injective map, it is easy to see it is not even continuous, let alone Lipschitz. The main work in this paper is to extend the sorting embedding to the case $d>1$ using a three-step procedure, first embed $\Rnd$ into a larger vector space $\R^{n\times D}$,  then apply $\downarrow$ in each column independently, and then perform a dimension reduction by a linear map into $\R^{2nd}$. Similar to the phase retrieval problem (\cite{bcmn,bod,balan16}), the redundancy introduced in the first step counterbalances the loss of information (here, relative order of one column with respect to another) in the second step.   

A summary of main results presented in this paper is contained in the following result.
\begin{thm}\label{t1}
Consider the metric space $(\hRnd,d)$.
\begin{enumerate}
\item (Polynomial Embedding) There exists a Lipschitz injective map
\[ \halpha:\hRnd\rightarrow\R^m \]
with $m=\left( \begin{array}{c}
\mbox{$d+n$} \\
\mbox{$d$}
\end{array} \right)$. Two explicit constructions of this map are given in   (\ref{eq:alpha1}) and (\ref{eq:alpha2}).
\item (Sorting based Embedding) There exists a class of bi-Lipschitz maps 
\[ \hbeta_{A,B}:(\hRnd,d)\rightarrow(\R^m,\norm{\cdot}) ~,~ \hbeta_{A,B}(\hat{X})=B\left(\hbeta_A(\hat{X})\right) \]
with $m=2nd$, where each map $\hbeta_{A,B}$ is the composition of two bi-Lipschitz maps: a full-rank linear operator $B:\R^{n\times D}\rightarrow \R^m$, with the nonlinear bi-Lipschitz map $\hbeta_A:\hRnd\rightarrow\R^{n\times D}$
 parametrized by a matrix $A\in\R^{d\times D}$ called "key". 
 Explicitly,
$\hbeta(\hat{X})=\downarrow(XA)$, where $\downarrow$ acts column-wise.
 These maps are characterized by the following properties:
 \begin{enumerate}
\item For $D=1+(d-1)n!$, any
$A\in\R^{d\times (1+(d-1)n!)}$ whose columns form a full spark frame defines a bi-Lipschitz map $\hbeta_A$ on $\hRnd$. 
Furthermore, a lower Lipschitz constant is given by the smallest $d^{th}$ singular value among all $d\times d$ sub-matrices of $A$,
 $\min_{J\subset[D],|J|=d}s_d(A[J])$.
\item For any matrix (``key") $A\in\R^{d\times D}$ such that the map 
$\hbeta_A$ is injective, then $\hbeta_A:(\hRnd,d)\rightarrow(\R^{n\times D},\norm{\cdot})$ is bi-Lipschitz. Furthermore, an upper Lipschitz constant is given by $s_1(A)$, the largest singular value of $A$.
\item Assume  $A\in\R^{d\times D}$ is such that the map 
$\hbeta_A$ is injective (i.e., a "universal key"). Then for almost any linear map $B:\R^{n\times D}\rightarrow\R^{2nd}$ the map $\hbeta_{A,B}=B\circ\hbeta_A$ is
bi-Lipschitz.
\end{enumerate}
\end{enumerate}
\end{thm}

An immediate consequence of this result is the following corollary whose proof is included in subsection \ref{subsec4.4}:
\begin{cor}
\label{c0}
Let $\beta:\R^{n\times d}\rightarrow\R^m$ induce a bi-Lipschitz embedding $\hbeta:\hRnd\rightarrow\R^m$ of the
metric space $(\hRnd,d)$ into $(\R^m,\norm{\cdot}_2)$. 
\begin{enumerate}
\item For any continuous function  $f:\R^{n\times d}\rightarrow\R$ 
invariant to row-permutation (i.e., $f(PX)=f(X)$ for every 
$X\in\R^{n\times d}$ and $P\in\Sn$) there exists a continuous
function $g:\R^m\rightarrow\R$ such that $f=g\circ\beta$.
Conversely, for any $g:\R^m\rightarrow\R$ continuous function, the 
function $f=g\circ\beta:\R^{n\times d}\rightarrow\R$ is continuous
and row-permutation invariant.
\item For any Lipschitz continuous function  $f:\R^{n\times d}\rightarrow\R$ 
invariant to row-permutation (i.e., $f(PX)=f(X)$ for every 
$X\in\R^{n\times d}$ and $P\in\Sn$) there exists a Lipschitz continuous
function $g:\R^m\rightarrow\R$ such that $f=g\circ\beta$.
Conversely, for any $g:\R^m\rightarrow\R$ Lipschitz continuous function, the 
function $f=g\circ\beta:\R^{n\times d}\rightarrow\R$ is Lipschitz continuous
and row-permutation invariant.
\end{enumerate}
\end{cor}
\vspace{5mm}

The structure of the paper is as follows. Section \ref{sec2} contains the algebraic embedding method and encoders $\alpha$ described at part (1) of Theorem \ref{t1}. Corollary \ref{cor2} contains part (1) of the main result stated above. Section \ref{sec3} introduces the sorting based embedding procedure and describes the key-based encoder $\beta$. Necessary and sufficient conditions for key universality are presented in Proposition
\ref{prop3.8}; the injectivity of the encoder described at part (2.a) of Theorem \ref{t1} is proved in Theorem \ref{t4}; the bi-Lipschitz property of any universal key described at part (2.b) of Theorem \ref{t1} is shown in Theorem \ref{t5}; the dimension reduction statement (2.c) of Theorem \ref{t1}
is included in Theorem \ref{t6}. Proof of Corollary \ref{c0} is presented in subsection \ref{subsec4.4}. Section \ref{sec4} contains applications to graph deep learning. These application use Graph Convolution Networks and the numerical experiments are carried out on two graph data sets: a chemical compound data set (QM9) and a protein data set (PROTEINS\_FULL). 

While the motivation of this analysis is provided by graph deep learning applications,
this is primarily a mathematical paper. Accordingly the formal theory is presented first, and then is followed by the machine learning application. Those interested in the application (or motivation) can skip directly to Section \ref{sec4}. 

{\bf Notations}. For an integer $d\geq 1$, $[d]=\{1,2,\ldots,d\}$. For a matrix $X\in\Rnd$,
 $x_1,\ldots x_d\in\R^n$ denote its columns, $X=[x_1\vert\cdots\vert x_d]$. All norms are Euclidean; for a matrix $X$, $\norm{X}=\sqrt{trace(X^TX)}=\sqrt{\sum_{k,j}|X_{k,j}|^2}$ denotes the Frobenius norm; for vectors $x$, $\norm{x}=\norm{x}_2=\sqrt{\sum_{j} |x_j|^2}$.

\subsection{Prior Works}

Several methods for representing orbits of vector spaces under the action of permutation (sub)groups have been studied in literature. Here we describe some of these results, without claiming an exhaustive literature survey.

A rich body of literature emanated from the early works on 
symmetric polynomials and group invariant representations of 
Hilbert, Noether, Klein and Frobenius. They are part of standard
commutative algebra and finite group representation theory. 

Prior works on permutation invariant mappings have predominantly employed some form of summing procedure, though some have alternatively employed some form of sorting procedure.

The idea of summing over the output nodes of an equivariant network has been well studied. 
The algebraic invariant theory goes back to Hilbert and Noether (for finite groups) and then continuing with the continuous invariant function theory of 
Weyl and Wigner (for compact groups), 
who posited that a generator function $\psi:X\rightarrow\R$ gives rise to a function $E:X\rightarrow\R$ invariant to the action of a finite group $G$ on $X$, $(g,x)\mapsto g.x$, via the averaging formula $E(x)=\frac{1}{|G|}\sum_{g\in G} \psi(g.x)$.

More recently, this approach provided the framework for universal approximation results of $G$-invariant functions.  \cite{maron2018invariant} showed that invariant or equivariant networks must satisfy a fixed point condition. The equivariant condition is naturally realized by GNNs. The invariance condition is realized by GNNs when followed by summation on the output layer, as was further shown in \cite{keriven2019universal}, \cite{pmlr-v97-maron19a} and \cite{lipman2022}. Subsequently, \cite{yarotsky2021universal} proved universal approximation results over compact sets for continuous functions invariant to the action of finite or continuous groups. In \cite{geerts2022}, the authors
obtained bounds on the separation power of GNNs in terms of the Weisfeiler-Leman (WL) tests by tensorizing the input-output mapping. 
\cite{sannai2020universal} studied approximations of equivariant maps, while \cite{NEURIPS2019_71ee911d} showed that if a GNN with sufficient expressivity is well trained, it can solve the graph isomorphism problem.

The authors of \cite{OrderMatters_2015arXiv151106391V} designed an algorithm for processing sets with no natural orderings. The algorithm applies an attention mechanism to achieve permutation invariance with the attention keys being generated by a Long-Short Term Memory (LSTM) network. Attention mechanisms amount to a weighted summing and therefore can be considered to fall within the domain of summing based procedures.

In \cite{GGsNN_2015arXiv151105493L}, the authors designed a permutation invariant mapping for graph embeddings. The mapping employs two separate neural networks, both applied over the feature set for each node. One neural network produces a set of new embeddings, the other serves as an attention mechanism to produce a weighed sum of those new embeddings.

Sorting based procedures for producing permutation invariant mappings over single dimensional inputs have been addressed and used by \cite{deepsets}, notably in their {\it max pooling} procedure.

The authors of \cite{qi2017pointnet} developed a permutation
invariant mapping 
$pointnet$ for point sets that is based on a $max$ function. The mapping takes in a set of vectors, processes each vector through a neural network followed by an scalar output function, and takes the maximum of the resultant set of scalars.

The paper \cite{zhang2018end} introduced {\it SortPooling}. {\it SortPooling} orders the latent embeddings of a graph according to the values in a specific, predetermined column. All rows of the latent embeddings are sorted according to the values in that column. While this gives rise to an injective map, it is easy to see it is not even continuous, let alone Lipschitz. The same issue
arises with any lexicographic ordering, including the well-known Weisfeiler-Leman embedding \cite{wl}.
Our paper introduces a novel method that bypasses this issue.

As shown in \cite{pmlr-v97-maron19a}, the sum pooling-based GNNs provides universal approximations for of any permutation invariant continuous function but only on \emph{compacts}. Our sorting based embedding removes the compactness restriction as well as it extends to all Lipschitz maps.

While this paper is primarily mathematical in nature, methods developed here are applied to two graph data sets, QM9 and PROTEINS\_FULL. Researchers have applied various graph deep learning techniques to both data sets. In particular, \cite{Gilmer_2017arXiv170401212G} studied extensively the QM9 data set, and compared their method with many other algorithms
proposed by that time.

\section{Algebraic Embeddings\label{sec2}}

The algebraic embedding presented in this section can be thought of a special kernel to project equation (\ref{eq:measure}) onto.

\subsection{Kernel Methods}
The kernel method employs a family of continuous kernels (test) functions, $\{K(x;y)~;~x\in\R^d~,~y\in Y\}$ parametrized/indexed by a set $Y$. 
The measure representation $\mu=a_{\infty}(X)$ in (\ref{eq:measure}) yields a nonlinear map
\[ \alpha:\R^{n\times d} \rightarrow C(Y)
~~,~~X \mapsto F(y)=\int_{R^d} K(x;y)d\mu \]
given by
\[ \alpha(X)(y)= \frac{1}{n}\sum_{k=1}^n K(x_k;y) \]
The embedding problem \ref{prob1}) can be restated as follows. One is asked
to find a finite family of kernels $\{K(x;y)~;~x\in\R^d~,~y\in Y\}$, 
 $m=|Y|$ so that
\begin{equation}
\label{eq:kernel}
\halpha:(\hRnd,d) \rightarrow l^2(Y)\sim (\R^m,\norm{\cdot}_2) ~~,~~ (\halpha(\hat{X}))_y = \frac{1}{n} \sum_{k=1}^n K(x_k;y)
\end{equation}
is injective, Lipschitz or bi-Lipschitz. 

Two natural choices for the kernel $K$ are the Gaussian kernel and the complex exponential (or, the Fourier) kernel:
\[ K_{G}(x,y) = e^{-\norm{x-y}^2/\sigma^2} ~~,
K_{F}(x,y) = e^{2\pi i \ip{x}{y}}
\]
where in both cases $Y\subset\R^d$. 
In this paper we analyze a different kernel, namely the polynomial kernel $K_P(x,y)=x_1^{y_1}x_2^{y_2}\cdots x_d^{y_d}$, $Y\subset\{0,1,2,\ldots,n\}^d$. 

\subsection{The Polynomial Embedding}

Since the polynomial representation is intimately related to the Hilbert-Noether algebraic invariants theory \cite{compuinvar} and the Hilbert-Weyl theorem, it is advantageous to start our construction from a different perspective.  

The linear space $\Rnd$ is isomorphic to $\R^{nd}$ by stacking the columns one on top of each other. In this case, the action of the permutation group $S_n$ can be recast as the action of the subgroup $I_d\otimes S_n$ of the bigger group $S_{nd}$ on $\R^{nd}$. Specifically, let us denote by $\sim_G$ the equivalence relation
\[ x,y\in\R^{nd}~~,~~x\sim_G y \Longleftrightarrow  y=\Pi x~,~{\rm for ~ some}~\Pi\in G \]
induced by a subgroup $G$ of $S_{nd}$. In the case
 $G=I_d\otimes S_n=\{diag_d(P)~,~P\in S_n\}$ of block diagonal permutation obtained by repeating $d$ times the same $P\in S_n$ permutation along the main diagonal, two vectors $x,y\in\R^{nd}$ are $\sim_G$ equivalent iff there is a permutation matrix $P\in S_n$ so that $y(1+(k-1)n:kn) = Px(1+(k-1)n:kn)$ for each $1\leq k\leq d$. In other words, each disjoint $n$-subvectors in $y$ and $x$ are related by the same permutation. In this framework, the Hilbert-Weyl theorem (Theorem 4.2, Chapter XII, in \cite{BifTheory2}) states that the ring of invariant polynomials is finitely generated.  The G\"{o}bel's algorithm (Section 3.10.2 in \cite{compuinvar}) provides a recipe to find a complete set of invariant polynomials. In the following we provide a direct approach to construct a complete set of  polynomial invariants. 
 
 Let $\R[\XXb_1,\XXb_2,...,\XXb_d]$ denote the algebra of polynomials in $d$-variables with real coefficients. 
Let us denote $\XX\in\Rnd$ a generic data matrix.
Each row of this matrix defines a 
linear form over $\XXb_1,...\XXb_d$,
 $\lambda_k = \XX_{k,1}\XXb_1+\cdots + \XX_{k,d}\XXb_d$.
 Let us denote by $\R[\XXb_1,\ldots,\XXb_d][\TTb]$ the algebra of polynomials in variable $\TTb$ with coefficients in the ring $\R[\XXb_1,\ldots,\XXb_d]$. Notice $\R[\XXb_1,\XXb_2,\ldots,\XXb_d][\TTb]=\R[\TTb,\XXb_1,\XXb_2,\ldots,\XXb_d]$ 
 by rearranging the terms according to degree in $\TTb$. 
 Thus $\lambda_k\in\R[\XXb_1,\ldots,\XXb_d]\subset\R[\XXb_1,\ldots,\XXb_d][\TTb]$ can be encoded as zeros of a polynomial $P_\XX$ of degree $n$ in variable $\TTb$ with coefficients in $\R[\XXb_1,\ldots,\XXb_d]$:
 \begin{equation}
     \label{eq:polyencoding}
     P_\XX(\TTb,\XXb_1,\ldots,\XXb_d) = \prod_{k=1}^n (\TTb-\lambda_k(\XXb_1,\ldots,\XXb_d))
     =\prod_{k=1}^n (\TTb-\XX_{k,1}\XXb_1-\ldots -\XX_{k,d}\XXb_d)
 \end{equation}
 Due to identification $\R[\XXb_1,\XXb_2,\ldots,\XXb_d][\TTb]=\R[\TTb,\XXb_1,\XXb_2,\ldots,\XXb_d]$,
 we obtain that \\
 $P_\XX\in \R[\TTb,\XXb_1,\XXb_2,\ldots,\XXb_d]$ is a homogeneous polynomial of degree $n$ in $d+1$ variables. Let $\R_n[\TTb,\XXb_1,\ldots,\XXb_d]$ denote the vector space of homogeneous polynomials in $d+1$ variables of degree $n$ with real coefficients. Notice the real dimension of this vector space is 
 \begin{equation}
     \label{eq:dimRn}
     \dim_\R \R_n[\TTb,\XXb_1,\ldots,\XXb_d] = \left( 
     \begin{array}{c}
     n+d \\
     d
     \end{array}
     \right) = \left(
     \begin{array}{c}
     n+d \\
     n
     \end{array}
     \right).
 \end{equation}
By noting that $P_\XX$ is monic in $\TTb$ (the coefficient of $\TTb^n$ is always 1)  we obtain an injective embedding of $\hRnd$ into $\R^m$ with 
$m=\dim_\R \R_n[\TTb,\XXb_1,\ldots,\XXb_d]-1$ via the coefficients of $P_\XX$ similar to (\ref{eq:poly}). This is summarized in the following theorem:
\begin{thm}
\label{t2}
The map $\alpha_0:\Rnd\rightarrow\R^{m-1}$ with $m=\left(\begin{array}{c}n+d \\ d \end{array} \right)$ given by the (non-trivial) coefficients of polynomial $P_\XX\in\R_n[\TTb,\XXb_1,\ldots,\XXb_d]$ lifts to an analytic embedding $\halpha_0$ of  $(\hRnd,d)$ into $\R^m$. Specifically, for $\XX\in\Rnd$ expand the polynomial
\begin{equation} \label{eq:PA}
P_\XX(\TTb,\XXb_1,\ldots,\XXb_d) = \prod_{k=1}^n (\TTb-\XX_{k,1}\XXb_1-\ldots -\XX_{k,d}\XXb_d)
 = \TTb^n + \hspace{-10mm}\sum_{\begin{array}{c}
 \mbox{$p_0,p_1,...,p_d\geq 0$} \\
 \mbox{$p_0+\cdots+p_d=n$} \\ 
 \mbox{$p_0<n$}
 \end{array}
}
\hspace{-5mm} c_{p_0,p_1,\ldots,p_d}\TTb^{p_0}\XXb_1^{p_1}\cdots \XXb_d^{p_d}
 \end{equation}
Then 
\begin{equation}\label{eq:alpha0}
\XX\in\Rnd\mapsto \alpha_0(\XX) = (c_{p_0,p_1,\ldots,p_d})_{(p_0,p_1,\ldots,p_d)\in I_{n,d}}
\end{equation}
where the index set is given by
\begin{equation} \label{eq:Ind}
I_{n,d}=\{(p_0,p_1,\ldots,p_d)~,~0\leq p_0,p_1,\ldots,p_d~,~p_0<n~,~p_0+p_1+\cdots+p_d=n\}
\end{equation}
and is of cardinal $|I_{n,d}|=m-1$. The map $\halpha_0:\hRnd\rightarrow\R^{m-1}$ is the lifting of $\alpha_0$ to the quotient space.
\end{thm}
{\bf Proof}

Since for any permutation $\pi$ with associated permutation matrix $\Pi\in\Sn$, 
\[ P_{\Pi \XX}(\TTb,\XXb_1,\cdots,\XXb_d) = \prod_{k=1}^n 
(\TTb-\XX_{\pi(k),1}\XXb_1-\ldots -\XX_{\pi(k),d}\XXb_d) =
P_\XX(\TTb,\XXb_1,\ldots,\XXb_d), \] 
it follows that $\alpha_0$ is invariant to the action of $\Sn$, $\alpha_0(\XX)=\alpha_0(\Pi\XX)$.
Thus $\alpha_0$ lifts to a map $\halpha_0$ on $\hRnd$. 

The coefficients of polynomial $P_\XX$ depend analytically on its roots (Vieta's formulas), hence on entries of matrix $\XX$. 

The only remaining claim is that if $\XX,\YY\in\Rnd$ so that $\alpha_0(\XX)=\alpha_0(\YY)$ then there is $\Pi\in\Sn$ so that $\YY=\Pi\XX$. Assume $P_\XX=P_\YY$. For each choice $(\XXb_1,\XXb_2,\ldots,\XXb_d)=(f(1),\ldots,f(d))$ in $\R^d$, 
the $n$ real zeros of the two polynomials in $\TTb$, $P_\XX(\TTb,f(1),\ldots,f(d))$ and $P_\YY(\TTb,f(1),\ldots,f(d))$ coincide. Therefore $\XX f\sim \YY f$ for each $f\in\R^d$, 
Let $D=1+(d-1)n!$ and choose $F\in\R^{d\times D}$ so that
each subset of $d$ columns are linearly independent, in other words, the set $\fc=\{f_1,f_2,\ldots,f_D\}$ formed by the $D$ columns of $F$ is a full spark frame in $\R^d$, see \cite{spark}. As proved in \cite{spark}, almost every such set is a full spark frame. Then for each $1\leq k\leq D$ there is a permutation $\Pi_k\in\Sn$ so that $\XX f_k=\Pi_k \YY f_k$. 
By the pigeonhole principle, since $|\Sn|=n!$, there are $1\leq k_1< k_2 < \cdots < k_d\leq D$ so that $\Pi_{k_1}=\Pi_{k_2}=\cdots=\Pi_{k_d}$.
Then $(\XX-\Pi_{k_1}\YY)f_{k_j}=0$ for every $1\leq j\leq d$. Since $\{f_{k_1},\ldots,f_{k_d}\}$ is linearly independent it follows that $\XX-\Pi_{k_1}\YY=0$. Thus $\XX\sim \YY$ which
 ends the proof of this result. $\qed$.

\begin{rmk}
The set of invariants produced by map $\alpha_0$ are proportional to those produced by the G\"{o}bel's algorithm in \cite{compuinvar}, \S 3.10.2. Indeed, the $nd$ 
primary invariants are given by 
\[ \{c_{p,n-p,0,\ldots,0}~,~0\leq p\leq n-1\}\cup\cdots\cup \{c_{p,0,\ldots,0,n-p}~,~0\leq p\leq n-1\}
\]
corresponding to the elementary symmetric polynomials in entries of each column. The secondary invariants correspond to the remaining coefficients that have at least 2 nonzero indices among $p_1,\ldots,p_d$.
\end{rmk}

The embedding provided by $\alpha_0$ is analytic and injective but is not globally Lipschitz because of the polynomial growth rate. Next we show how a simple modification of this map 
will make it Lipschitz.
First, let us denote by $L_0$ the Lipschitz constant of $\alpha_0$ when restricted to the closed unit ball $B_1(\Rnd):\{X\in \Rnd~,~\norm{X}\leq 1\}$ of $\Rnd$, i.e. $\norm{\alpha_0(\XX)-\alpha_0(\YY)}\leq L_0\norm{\XX-\YY}$ for any $\XX,\YY\in\Rnd$ with $\norm{\XX},\norm{\YY}\leq 1$. 
 Second, let
\begin{equation}
    \label{eq:phi}
    \varphi_0:\R\rightarrow [0,1]~,~
    \varphi_0(x)=min(1,\frac{1}{x})=\left\{ \begin{array}{rcl}
    1 & if & \mbox{$x\leq 1$} \\
    \mbox{$\frac{1}{x}$} & if & \mbox{$x > 1$}
    \end{array}\right.
\end{equation}
be a Lipschitz monotone decreasing function with Lipschitz constant 1.
\begin{cor}\label{cor2}
Consider the map:
\begin{equation}\label{eq:alpha1}
\alpha_1:\Rnd\rightarrow\R^m~~,~~
\alpha_1(\XX) = \left( \begin{array}{c}
\mbox{$\alpha_0\bigg (\varphi_0(\norm{\XX})\XX \bigg )$} \\
\mbox{$\norm{\XX}$}
\end{array}\right),
\end{equation}
with $m=\left(\begin{array}{c}n+d \\ d \end{array} \right)$.
The map $\alpha_1$ lifts to an injective and globally Lipschitz map $\halpha_1:\hRnd\rightarrow\R^m$ with Lipschitz constant $Lip(\halpha_1) \leq \sqrt{1+L_0^2}$.
\end{cor}
{\bf Proof}

Clearly $\alpha_1(\Pi \XX)=\alpha_1(\XX)$ for any $\Pi\in\Sn$ and $\XX\in\Rnd$. Assume now that $\alpha_1(\XX)=\alpha_1(\YY)$. Then $\norm{\XX}=\norm{\YY}$ and since $\halpha_0$ is injective on $\hRnd$ it follows $\varphi(\norm{\XX})\XX = \Pi \varphi(\norm{\YY})\YY$ for some $\Pi\in\Sn$. Thus $\XX\sim \YY$ which proves $\alpha_1$ lifts to an injective map on $\hRnd$. 

Now we show $\halpha_1$ is Lipschitz on $(\hRnd,d)$ of appropriate Lipschitz constant. Let $\XX,\YY'\in\Rnd$ and $\Pi_0\in\Sn$ so that $d(\hat{\XX},\hat{\YY'})=\norm{\XX-\Pi_0 \YY'}$. Let $\YY=\Pi_0 \YY'$ so that $d(\hat{\XX},\hat{\YY})=\norm{\XX-\YY}$. 

Choose two matrices $\XX,\YY\in\Rnd$. We claim $\norm{\alpha_1(\XX)-\alpha_1(\YY)}\leq \sqrt{1+L_0^2}
\norm{\XX-\YY}$.
This follows from two observations: 

(i) The map
\[ \XX \mapsto \rho(\XX):=\varphi_0(\norm{\XX})\XX \]
is the nearest-point map to (or, the metric projection map onto) the convex closed set $B_1(\Rnd)$. This means $\norm{\varphi_0(\norm{\XX})\XX - Z}\leq \norm{\XX-Z}$ for any $Z\in B_1(\Rnd)$. 

(ii) The nearest-point map to a convex closed subset of a Hilbert space is Lipschitz with constant 1, i.e. it shrinks distances, see \cite{phelps56}.

These two observations yield:
\begin{multline*} 
\norm{\alpha_1(\XX)-\alpha_1(\YY)}^2  =  \norm{\alpha_0(\rho(\YY))
- \alpha_0(\rho(\YY) )}^2 + |\norm{\XX}-\norm{\YY}|^2  \\
 \leq  
L_0^2 \norm{\rho(\XX)-\rho(\YY) }^2 + \norm{\XX-\YY}^2 \leq (1+L_0^2)\norm{\XX-\YY}^2 .
\end{multline*}
This concludes the proof of this result. $\qed$
\vspace{5mm}

A simple modification of $\phi_0$ can produce a $C^\infty$ map by smoothing it out around $x=1$.

On the other hand the lower Lipschitz constant of $\halpha_1$ is 0 due to terms of the form $\XX_{i,j}^k$ with $k\geq 2$. 
In \cite{Cahill19}, the authors built a Lipschitz map by  a retraction to the unit sphere instead of unit ball. 
Inspired by their construction, a modification of $\alpha_0$ in their spirit reads:
\begin{equation}
    \label{eq:alpha2}
\alpha_2:\Rnd\rightarrow\R^m~~,~~
\alpha_2(\XX)=\left( \begin{array}{c}
\mbox{$\norm{\XX}\alpha_0\bigg ( \frac{\XX}{\norm{\XX}} \bigg )$} \\
\mbox{$\norm{\XX}$}
\end{array}\right)~,~if~\XX\neq 0~~,~and~~\alpha_2(0)=0.
\end{equation}
It is easy to see that $\alpha_2$ satisfies the non-parallel property in \cite{Cahill19} and is Lipschitz with a slightly better constant than $\alpha_1$ (the constant is determined by the tangential derivatives of $\alpha_0$). 
But, for the same reasons as in \cite{Cahill19} this map is not bi-Lipschitz. 

\subsection{Dimension reduction in the case $d=2$ and consequences}

In this subsection we analyze the case $d=2$. 
The embedding dimension for $\alpha_0$ is $\left( \begin{array}{c} n \\ 2 \end{array}\right)-1=\frac{n(n-1)}{2}-1$. 
On the other hand, consider the following approach. 
Each row of $\XX$ defines a complex number $z_1=\XX_{1,1}+i\,\XX_{1,2}$, ... , $z_n=\XX_{n,1}+i\,\XX_{n,2}$ that
can be encoded by one polynomial of degree $n$ with complex coefficients $Q\in\C_n[t]$,
\[ Q(\TTb) = \prod_{k=1}^n (\TTb-z_k) = \TTb^n + \sum_{k=0}^{n-1}
\TTb^k q_k \]
The coefficients of $Q$ provide a $2n$-dimensional real embedding $\zeta_0$,
\[ \zeta_0:\R^{n\times 2}\rightarrow\R^{2n}~~,~~\zeta_0(\XX)=(Re(q_{n-1}),Im(q_{n-1}),\ldots,Re(q_{0}),Im(q_0)) \]
with properties similar to those of $\alpha_0$. 
One can similarly modify this embedding to obtain a globally Lipschitz embedding $\hat{\zeta}_1$ of $\hat{R_{n,2}}$ 
into $\R^{2n+1}$. 

It is instructive to recast this embedding in the framework of commutative algebras. Indeed, let $\langle \XXb_1-1,\XXb_2^2+1 \rangle$ denote
the ideal generated by polynomials $\XXb_1-1$ and $\XXb_2^2+1$
in the algebra $\R[\TTb,\XXb_1,\XXb_2]$. Consider the quotient space
 $\R[\TTb,\XXb_1,\XXb_2]/\langle \XXb_1-1,\XXb_2^2+1 \rangle$ and the quotient map
 $\sigma:\R[\TTb,\XXb_1,\XXb_2]\mapsto \R[\TTb,\XXb_1,\XXb_2]/\langle\XXb_1-1,\XXb_2^2+1\rangle$.
 In particular, let $S=\sigma(\R_n[\TTb,\XXb_1,\XXb_2])$ denote the vector space projected through this quotient map.
Then a basis for $S$ is given by $\{1,\TTb,\ldots,\TTb^n,\XXb_2,\XXb_2 \TTb,\ldots,\XXb_2 \TTb^{n-1},\XXb_2 \TTb^n\}$. Thus $\dim S=2n+2$. 
Let 
$\Ssc=\{P_\XX~,~\XX\in\R^{n\times 2} \}\subset\R_2[\TTb,\XXb_1,\XXb_2]$ 
denote the set of polynomials realizable as in (\ref{eq:PA}).
Then the fact that $\hat{\zeta}_0:\R^{n\times 2}\rightarrow\R^{2n}$ 
is injective is equivalent to the fact that $\sigma{\vert}_{\Ssc}:\Ssc\rightarrow S$ is injective.
On the other hand note 
\[
\sigma(\Ssc)\subset \TTb^n+
span_\R\{1,\TTb,\ldots,\TTb^{n-1},\XXb_2,\XXb_2 \TTb, \ldots,\XXb_2 \TTb^{n-1} \} \]
where the last linear subspace is of dimension $2n$. 

In the case $d=2$ we obtain the identification
$\R[\TTb,\XXb_1,\XXb_2]/\langle \XXb_1-1,\XXb_2^2+1 \rangle = \C[\TTb]$ due to uniqueness of polynomial factorization.

This observation raises the following {\em open problem}:

For $d>2$, is there a non-trivial ideal 
$I=\langle Q_1,\ldots,Q_r \rangle \subset\R[\TTb,\XXb_1,\ldots,\XXb_d]$
so that the restriction $\sigma{\vert}_{\Ssc}$
of the quotient map $\sigma:\R[\TTb,\XXb_1,\ldots,\XXb_d]\rightarrow
\R[\TTb,\XXb_1,\ldots,\XXb_d]/I$ is injective? Here $\Ssc$ denote the set of polynomials in $\R_n[\TTb,\XXb_1,\ldots,\XXb_d]$ realizable via (\ref{eq:PA}).
\begin{rmk}
One may ask the question whether the quaternions can be
utilized in the case $d=4$. While the quaternions form an associative division algebra, unfortunately polynomials have in general an infinite number of factorization. This prevents an immediate extension of the previous construction to the case $d=4$.   
\end{rmk}

\begin{rmk}
Similar to the construction in \cite{Cahill19}, a linear dimension reduction technique may be applicable here (which, in fact, may answer the open problem above) which would reduce the embedding dimension to $m=2nd+1$ (twice the intrinsec dimension plus one for the homogenization variable). 
However we did not explore this approach since, even if possible, it would not produce a bi-Lipschitz embedding. 
Instead we analyze the linear dimension reduction technique in the next section in the context of sorting based embeddings. 
\end{rmk}

\section{Sorting based Embedding\label{sec3}}

In this section we present the extension of the sorting embedding (\ref{eq:ord}) to the case $d>1$.

The embedding is performed by a linear-nonlinear transformation that resembles the phase retrieval problem. 
Consider a matrix $A\in\R^{d\times D}$ and the induced nonlinear 
transformation:

\begin{equation}
\label{eq:qA}
\beta_A:\Rnd\rightarrow\R^{n\times D}~~,~~\beta_A(X)=\downarrow (XA)
\end{equation}
where $\downarrow$ is the monotone decreasing sorting operator acting in each column independently. Specifically, let 
$Y=XA\in\R^{n\times D}$ and note its column vectors
$Y=[y_1,y_2,\ldots,y_D]$. Then 
\[ \beta_A(X)=\left[ \begin{array}{cccc}
\mbox{$\Pi_1 y_1$} & \mbox{$\Pi_2 y_1$} & \cdots & \mbox{$\Pi_D y_D$}
\end{array} \right] \]
for some $\Pi_1,\Pi_2,\ldots,\Pi_D\in\Sn$ so that each column is sorted monotonically decreasing:
\[ (\Pi_k y_k)_1\geq (\Pi_k y_k)_2\geq \cdots\geq (\Pi_k y_k)_n. \]
Note the obvious invariance $\beta_A(\Pi X)=\beta_A(X)$ for any $\Pi\in\Sn$ and $X\in\Rnd$. Hence $\beta_A$ 
lifts to a map $\hat{\beta_A}$ on $\hRnd$. 
\begin{rmk}
Notice the similarity to the phase retrieval problem, e.g.,  \cite{balan16}, where the data is obtained via a linear transformation of the 
input signal followed by the nonlinear operation of taking the absolute value of the frame coefficients. Here the nonlinear transformation is implemented by sorting the coefficients. 
In both cases it represents the action of a particular
subgroup of the unitary group.   
\end{rmk}

In this section we analyze necessary and sufficient conditions so that maps of type (\ref{eq:qA}) are injective, or injective almost everywhere. 
First a few definitions.

\begin{defn}
A matrix $A\in\R^{d\times D}$ is called a \emph{universal key} (for $\Rnd$) if $\hat{\beta_A}$ is injective
on $\hRnd$.
\end{defn}
In general we refer to $A$ as a {\em key} for encoder $\beta_A$. 
\begin{defn}
Fix a matrix $X\in\Rnd$. A matrix $A\in\R^{d\times D}$ is said \emph{admissible} (or an {\em admissible key}) for $X$ if for any $Y\in\Rnd$ so that $\beta_A(X)=\beta_A(Y)$ then $Y=\Pi X$ for some $\Pi\in\Sn$.  
\end{defn}
In other words, $\hat{\beta_A}^{-1}(\hat{\beta_A}(\hat{X}))=\{\hat{X}\}$.
We let $\Acr_{D}(X)$, or simply $\Acr(X)$, denote the set of admissible keys for $X$. 
\begin{defn}
Fix $A\in\R^{d\times D}$. A matrix $X\in\Rnd$ is said to be {\em separated} by $A$ if $A\in\Acr(X)$.
\end{defn}
For a key $A$, we let $\Ssc_{n}(A)$, or simply $\Ssc(A)$, denote the set of {\em matrices separated by $A$}. Thus a matrix $X\in\Ssc_n(A)$ if and only if, for any matrix $Y\in\R^{n\times d}$, if $\beta_A(X)=\beta_A(Y)$ then $X\sim Y$.

Thus a key $A$ is universal if and only if $\Ssc_n(A)=\Rnd$.

Our goal is to produce keys that are admissible for all matrices in $\Rnd$, or at least for almost every data matrix.
As we show in Proposition \ref{prop3.6} below this requires that $D\geq d$ and $A$ is full rank. In particular this means that the columns of $A$ form a frame for $\R^d$.  

\subsection{Characterizations of $\Acr(X)$ and $\Ssc(A)$}

We start off with simple linear manipulations of sets of admissible keys and separated data matrices.

\begin{prop}\label{prop3.5}
Fix $A\in\R^{d\times D}$ and $X\in\Rnd$.
\begin{enumerate}
    \item For an invertible $d\times d$ matrix $T\in\R^{d\times d}$,
    \begin{equation}\label{eq:TA}
        \Ssc_n(TA) = \Ssc_n(A)T^{-1}.
    \end{equation}
    In other words, if $X$ is separated by $A$ then $XT^{-1}$ is separated by $TA$.
    
    \item For any permutation matrix $L\in\SS_D$ and diagonal invertible matrix $\Lambda\in\R^{D\times D}$,
    \begin{equation}\label{eq:AL}
        \Ssc_n(AL\Lambda)=\Ssc_n(A\Lambda L) = \Ssc_n(A).
    \end{equation}
    In other words, if $X$ is separated by $A$ then $X$ is separated also by $AL\Lambda$ as well as by $A\Lambda L$.
    
    \item Assume $T\in\R^{d\times d}$ is a $d\times d$ invertible matrix. Then
    \begin{equation}\label{eq:XT}
        \Acr_D(XT)=T^{-1}\Acr_D(X).
    \end{equation}
    In other words, if $A$ is an admissible key for $X$ then $T^{-1}A$ is an admissible key for $XT$.
\end{enumerate}
\end{prop}
{\bf Proof}

The proof is immediate, but we include it here for convenience of the reader. 

(1) Denote $B=TA$. Let $Y\in\R^{n\times d}$. Then
\[ \beta_B(Y)=\beta_B(X) \Longleftrightarrow \downarrow(XB)=\downarrow(YB) \Longleftrightarrow \downarrow(XTA)=\downarrow(YTA)
\Longleftrightarrow\beta_A(XT)=\beta_A(YT). \]
Thus, if $X\in\Ssc_n(A)$ and $Y'\in\Rnd$ so that $\beta_B(Y')=\beta_B(X')$ with $X'=XT^{-1}$, then $\beta_A(Y'T)=\beta_A(X)$. Therefore there exists $\Pi\in\Sn$ so that $Y'T=\Pi X$. Thus $Y'\sim X'$. Hence  $X'\in\Ssc_n(B)$.
This shows $\Ssc_n(A)T^{-1}\subset \Ssc_n(TA)$. The reverse include follows by replacing $A$ with $TA$ and $T$ with $T^{-1}$. Together they prove (\ref{eq:TA}).

(2) Let $Y\in\Rnd$ such that $\beta_{AL\Lambda}(X)=\beta_{AL\Lambda}(Y)$. 
For every $1\leq j\leq D$ let $k\in [D]$ be so that $L_{jk}=1$. 

If $\Lambda_{kk}>0$ then $\downarrow((XA)_j)=\downarrow((YA)_j)$. 

If $\Lambda_{kk}<0$ then  $\downarrow(-(XA)_j)=\downarrow(-(YA)_j)$.
But this implies also $\downarrow((XA)_j)=\downarrow((YA)_j)$ since
$\downarrow(-z)=L_0\downarrow(z)$ where $L_0$ is the permutation matrix that has 1 on its main antidiagonal.

Either way, $\downarrow((XA)_j)=\downarrow((YA)_j)$. Hence 
$\downarrow(XA)=\downarrow(YA)$. Therefore $X\sim Y$ and thus $X\in \Ssc_n(AL\Lambda)$. This shows $\Ssc_n(A)\subset\Ssc_n(AL\Lambda)$.
the reverse inclusion follows by a similar argument.
Finally, notice $\{L\Lambda\}$ forms a group since $L^{-1}\Lambda L$ is also a diagonal matrix. This shows $\Ssc_n(A\Lambda L)=\Ssc(AL\Lambda')$ for some diagonal matrix $\Lambda'$, and the conclusion (\ref{eq:AL}) then follows.

(3) The relation (\ref{eq:XT}) follows from noticing $\beta_{T^{-1}A}(Y)=\beta_A(YT)$.   $\qed$

Relation (\ref{eq:AL}) shows that, since $A$ is assumed full rank, without loss of generality we can assume the first $d$ columns are linearly independent. Let $V$ denote the first $d$ columns of $A$ so that
\begin{equation}\label{eq:AA}
A = V\left[ \begin{array} {ccc}
\mbox{$I$} & \mbox{$\vert$} & \mbox{$\tilde{A}$} 
\end{array} \right] 
\end{equation}
where $\tilde{A}\in\R^{d\times (D-d)}$. 
The following result shows that, unsurprisingly, when $D=d>1$, almost every matrix $X$ is not separated by $A$. By Proposition \ref{prop3.5} we can reduce the analysis to the case $A=I$ by a change of coordinates.
\begin{prop}\label{prop3.6}
Assume $D=d>1$, $n>1$. Then
\begin{enumerate}
    \item The set of data matrices not separated by $I_d$ includes:
\begin{equation}\label{eq:SA1}
    \Ddelta:=\{ X\in\R^{n\times d}~,~\exists i,j,k,l~,~ 1\leq i<j\leq n,1\leq k<l\leq d~\Rightarrow ~ X_{i,k}\neq X_{j,k}~\& ~X_{i,l}\neq X_{j,l} \} \subset \Ssc_n(I_d)^c.
\end{equation}
\item The set $\Ddelta$ is generic with respect to Zariski topology, i.e., open and dense. Specifically, its complement is the zero set of the polynomial
\[ P(X) = \sum_{1\leq i<j\leq n}~\sum_{1\leq k<l\leq d}(X_{i,k}-X_{j,k})^2(X_{i,l}-X_{j,l})^2. \]
\item For an invertible matrix $A\in\R^{d\times d}$, 
\[ \Ddelta \cdot A^{-1} \subset \Ssc_n(A)^c. \] 
Hence almost every matrix (w.r.t. Lebesgue measure) 
$X\in\Rnd$ is \emph{not} separated by $A$. 
\end{enumerate}
\end{prop}

{\bf Proof}

(1) 
We need to show that any matrix $X$ that on some columns $k$ and $l$ has distinct elements on same row positions is not separated by $I_d$. Indeed if $X$ is such a matrix, let $Y$ denote a copy of $X$ except on those 4 entries where we set
\[ Y_{i,k}=X_{j,k}~~,~~Y_{j,k}=X_{i,k}~~,~~Y_{i,l}=X_{i,l}~~,~~Y_{j,l}=X_{j,l}. \]
Note $X\not\sim Y$ yet $\downarrow(X)=\downarrow(Y)$. Hence such matrices are not separated by $I_d$.

(2) By negation, the complement of $\Ddelta$ is given by
\[ \Ddelta^c = \{ X\in\R^{n\times d} ~,~\forall i,j,k,l~,~
1\leq i<j\leq n,1\leq k<l\leq d \& (X_{i,k}=X_{j,k}~or~X_{i,l}=X_{j,l}) \} \]
This shows $\Ddelta^c$ is the zero set
of polynomial $P$ as claimed. Thus $\Ddelta^c$ is a closed Zariski set. Its complement is generic with respect to the Zariski topology 
since $\Ddelta^c\neq \R^{n\times d}$. 

(3) The inclusion is immediate. Density claim follows from this inclusion.
$\Box$
\vspace{3mm}

On the other hand, extending the identity matrix by
only one column produces an almost universal key:
\begin{prop}\label{prop3.7} Assume $d\geq 2$ and $n\geq 3$.

Let $a\in\R^d$ be a vector with non-zero entries, i.e., $\prod_{i=1}^d a_i \neq 0$. Let $A=\left[ \begin{array}{ccc} I_d & \vert & a \end{array} \right]\in\R^{d\times (d+1)}$ be a key. 
Then $\Ssc_n(A)$ is generic with respect to the Zariski topology (i.e., open and dense), however $\Ssc_n(A)\neq \R^{n\times d}$. In particular, its complement $\Ssc_n(A)^c:=\Rnd\setminus \Ssc_n(A)$ is non-empty but has 
Lebesgue measure zero. Thus almost every matrix $X\in\R^{n\times d}$
 is separated by $A$.
\end{prop}

{\bf Proof}

First we show that $\Ssc_n(A)\neq \R^{n\times d}$. 
Consider the matrices $X,Y\in\R^{n\times d}$ full of zeros except for the 3x2 top left corner where:
\[ X_{1,1}=Y_{1,1}=\frac{1}{a_1} ~~,~~
X_{2,1}=Y_{2,1} = -\frac{1}{a_1} ~~,~~X_{3,1}=Y_{3,1}=0\]
\[ X_{1,2} = Y_{3,2} = -\frac{1}{a_2} ~~,~~
 X_{2,2} = Y_{1,2} = 0 ~~,~~
X_{3,2}=Y_{2,2}=\frac{1}{a_2} \]
Clearly $\beta_A(X)=\beta_A(Y)$ (the two left columns and the last column contain $1$, $0$ repeated $n-2$ times and $-1$) and yet $X\not\sim Y$.

Next we show that $\Ssc_n(A)^c$ is included in a finite union 
of linear spaces each of positive codimension. This proves
the clam.

To simplify notation we introduce the following two operators. Let $\Pi,\Pi_0,\Pi_1,\cdots,\Pi_d\in \Sn$
denote permutation matrices of size n. For $X\in\R^{n\times d}$ denote by $x_1,\ldots,x_d$ its columns. Thus $X=\left[ x_1 \vert x_2 \vert \cdots \vert x_d \right]$.
\[ L_{\Pi_0,\Pi_1,\ldots,\Pi_d} :\R^{n\times d}\rightarrow\R^d ~~,~~ L_{\Pi_0,\Pi_1,\ldots,\Pi_d} X = \Pi_0 Xa - (a_1 \Pi_1 x_1 + \cdots a_d \Pi_d x_d) \]
and
\[ M_{\Pi,\Pi_1,\ldots,\Pi_d}:\R^{n\times d} \rightarrow \R^{n\times d} ~~,~~
M_{\Pi,\Pi_1,\ldots,\Pi_d} X = \Pi X - \left[ 
\begin{array}{ccccc} \Pi_1 x_1 & \vert & \cdots & \vert & \Pi_d x_d \end{array} \right]
\]
A matrix $X\in\Rnd$, $X=[x_1\vert\cdots\vert x_d]$ 
is not separated by $A=[I_d\vert a]\in\R^{d\times (d+1)}$, i.e., 
$X\in\Ssc_n(A)^c$ if there are permutation matrices $\Pi_1,\ldots,\Pi_d\in\Sn$ such that the matrix
 $Y=[\Pi_1 x_1\vert\cdots\vert\Pi_d x_d]$ satisfies:
 \[ Y\not\sim X~~and~~\downarrow(X\cdot a) = \downarrow(Y\cdot a) \]
 This is equivalent to say:
 \[ \exists \Pi_0\in\Sn ~~,~~\Pi_0 Xa - (a_1\Pi_1 x_1
 +\cdots + a_d \Pi_d x_d) = 0 \]
 \[ \forall \Pi\in\Sn \exists k\in [d]~~,~~(\Pi-\Pi_k)x_k \neq 0 \]
 Hence
\[ \Ssc_n(A)^c=\bigcup_{(\Pi_0,\Pi_1,\ldots,\Pi_d)\in \Sn^{d+1}} \left(
\ker\, L_{\Pi_0,\Pi_1,\ldots,\Pi_d} \setminus \left( \bigcup_{\Pi\in\Sn}
\ker\, M_{\Pi,\Pi_1,\ldots,\Pi_d} \right) \right) \]

Let $\Delta$ denote the diagonal in $\Sn^{d+1}$,
\[ \Delta = \{ (\Pi_0,\Pi_1,\Pi_2,\ldots,\Pi_d)\in\Sn^{d+1}
~~,~~\Pi_1=\Pi_2=\cdots=\Pi_d \} \]
parametrized by the first two permutation matrices. 
For any $(\Pi_0,\Pi_1,\ldots,\Pi_d)\in\Delta$, we have
$\ker\, M_{\Pi_1,\Pi_1,\ldots,\Pi_d}=\Rnd$. Thus
\[ \ker\, L_{\Pi_0,\Pi_1,\ldots,\Pi_d} \setminus \left( \bigcup_{\Pi\in\Sn}
\ker\, M_{\Pi,\Pi_1,\ldots,\Pi_d} \right) = \emptyset \]
It follows:
\[ \Ssc_n(A)^c=\bigcup_{(\Pi_0,\Pi_1,\ldots,\Pi_d)\in \Sn^{d+1}
\setminus\Delta} \left(
\ker\, L_{\Pi_0,\Pi_1,\ldots,\Pi_d} \setminus \left( \bigcup_{\Pi\in\Sn}
\ker\, M_{\Pi,\Pi_1,\ldots,\Pi_d} \right) \right) \]
Consider now $(\Pi_0,\Pi_1,\ldots,\Pi_d)\in\Sn^{d+1}\setminus\Delta$. Then
\[ \ker\, L_{\Pi_0,\Pi_1,\ldots,\Pi_d} = \ker\, L_{I,\Pi_0^{-1}\Pi_1,\ldots,\Pi_0^{-1}\Pi_d} \]
Hence there is $k\in [d]$ so that $\Pi_0^{-1}\Pi_k\neq I$.
Choose $x_k\in\R^n$ so that $(\Pi_0^{-1}\Pi_k)x_k\neq x_k$.
 Set $x_j=0$ for $j\in[d]$, $j\neq k$ and consider the matrix 
 $X=[x_1\vert\cdots\vert x_d]$.
 Then $L_{\Pi_0,\Pi_1,\ldots,\Pi_d}X=a_k(\Pi_0-\Pi_k)x_k\neq 0$.
 This show that $\ker\, L_{\Pi_0,\Pi_1,\ldots,\Pi_d}\neq \Rnd$ 
 and hence it is a subspace of positive codimension. We obtain:
 \[ \Ssc_n(A)^c\subset \bigcup_{(\Pi_0,\Pi_1,\ldots,\Pi_d)\in \Sn^{d+1}
\setminus\Delta} \ker\, L_{\Pi_0,\Pi_1,\ldots,\Pi_d} \]
This shows that $\Ssc_n(A)^c$ is included in a finite union 
of proper subspaces of $\Rnd$ which in turn is a closed set with respect to the Zariski topology of empty interior. This ends the proof of this result. $\Box$

\vspace{5mm}

The next result provides a characterization of the set $\Ssc_n(A)$.
To do so we need to introduce additional notation that extends the operators $L_{\Pi_0,\ldots,\Pi_d}$ and $M_{\Pi,\ldots,\Pi_d}$ defined in the proof of 
Proposition \ref{prop3.7}.
For $E_1,E_2,\ldots,E_d\in\R^{n\times n}$ and $b\in\R^d$, with
$b=\left(b_1,b_2,\cdots,b_d \right)^T$ define
\begin{equation}
    \label{eq:LAb}
    L_{E_1,E_2,\ldots,E_d;b}:\R^{n\times d}\rightarrow\R^n ~~,~~
X=\left[x_1\vert x_2\vert\cdots\vert x_d\right] \mapsto    L_{E_1,E_2,\ldots,E_d;b}(X)=b_1 E_1 x_1+\cdots+b_d E_d x_d.
\end{equation}
\begin{prop}
\label{prop3.8} Fix $a_1,\ldots,a_{D-d}\in\R^{d}$ and consider the
key $A=\left[ I_d\vert a_1\vert\cdots\vert a_{D-d} \right] \in\R^{d\times D}$. Let $X=\left[x_1\vert x_2\vert\cdots\vert x_d\right]$.
\begin{enumerate}
    \item  $X\in \Ssc_n(A)^c:=\Rnd\setminus \Ssc_n(A)$ if and only if
   there are $\Pi_1,\Pi_2,\ldots,\Pi_d,\Xi_1,\ldots,\Xi_{D-d}\in\Sn$ such that:
   \begin{enumerate}
       \item $\forall j\in[D-d]$, $\left[(\Xi_j-\Pi_1)x_1,\ldots,(\Xi_j)-\Pi_d)x_d\right]a_j = 0$
       \item $\forall \Pi\in\Sn~\exists k\in[d]$ so that $(\Pi_k-\Pi)x_k\neq 0$.
   \end{enumerate}
   \item The following hold true:
   \begin{equation}
       \label{eq:SA}
       \Ssc_n(A) = \bigcap_{\begin{array}{c} \Pi_1,\ldots,\Pi_d\in\Sn \\
       \Xi_1,\ldots,\Xi_{D-d}\in\Sn \end{array}} \left[ 
\bigcup_{j=1}^{D-d}\left( \ker\,\ L_{\Xi_j-\Pi_1,\ldots,\Xi_j-\Pi_d;a_j}\right)^c \bigcup
\bigcup_{\Pi\in\Sn} 
\bigcap_{k=1}^d \ker\, L_{\Pi_1-\Pi,\ldots,\Pi_d-\Pi;\delta_k} \right]
   \end{equation}
   and
   \begin{equation}
       \label{eq:SAc}
       \Ssc_n(A)^c = \bigcup_{\begin{array}{c} \Pi_1,\ldots,\Pi_d\in\Sn \\
       \Xi_1,\ldots,\Xi_{D-d}\in\Sn \end{array}} \left( 
\bigcap_{j=1}^{D-d} \ker\,\ L_{\Xi_j-\Pi_1,\ldots,\Xi_j-\Pi_d;a_j}\right) \bigcap \left(
\bigcup_{\Pi\in\Sn} 
\bigcap_{k=1}^d \ker\, L_{\Pi_1-\Pi,\ldots,\Pi_d-\Pi;\delta_k} \right)^c
   \end{equation}
   where $\delta_k=(0,\ldots,0,1,0,\ldots,0)^T$ has only one 1 on the $k^{th}$ position.
\end{enumerate}
\end{prop}

{\bf Proof}

The proof is a consequence of linear algebra analysis of 
map $\beta_A$.

(1) Assume $X$ is  not separated by $A$. Then there is $Y\in\Rnd$
so that $\beta_A(X)=\beta_A(Y)$ yet $X\not\sim Y$.

Let $Y=[y_1\vert\cdots\vert y_d]$. Then $\beta_A(X)=\beta_A(Y)$ implies that there are permutation
matrices $\Pi_1,\ldots,\Pi_d,\Xi_1,\ldots,\Xi_{D-d}\in\Sn$
 so that:
 \[ y_1=\Pi_1 x_1,\ldots, y_d=\Pi_d x_d, Ya_1=\Xi_1 X a_1,\ldots,
 Ya_{D-d}=\Xi_{D-d}Xa_{D-d} \]
Substituting the expressions for $y_1,\ldots,y_d$ provided by the first $d$ equations into the latter  $D-d$ equations, we obtain
 part 1.(a). 
 
For same $Y$, the condition $X\not\sim Y$ implies that for
every $\Pi\in\Sn$, $Y-\Pi X\neq 0$. Thus part 1(b) is proved.

(3) Equation (\ref{eq:SAc}) is a transcription of part 1.

(2) Equation (\ref{eq:SA}) follows from (\ref{eq:SAc}) 
by taking the complement. $\Box$

\subsection{Construction of universal keys}
In this subsection we construct universal keys. Proposition \ref{prop3.8} provides us with an algorithm to check whether a key $A$ is universal. Unfortunately the
algorithm has an exponential complexity in data size. 

If the key $A\in\R^{d\times D}$ is universal
then $A$ must have full rank. Therefore there are permutation matrix $L\in\SS_D$
 and invertible $T\in GL(d,\R)$ so that $A=T\left[ I_d~ \tilde{A}\right]L$,
  with $\tilde{A}\in\R^{d\times (D-d)}$.
Proposition \ref{prop3.5} shows that $A$ is a universal key if and only if
 $\left[ I_d ~\tilde{A}\right]$ is a universal key. 
This observation allows us to prove the main result of this subsection stated earlier as part b of  Theorem \ref{t2}. Recall a set of vectors $\{f_1,\ldots,f_m\}$ in a 
linear space $V$ of finite dimension $n\leq m$ is called a 
{\em full spark frame} if
any subset of $n$ vectors is linearly independent. See \cite{spark,oussa} for more
information and explicit constructions of full spark frames.
\begin{thm}
 \label{t4}
Consider the metric space $(\hRnd,d)$. Set $D=1+(d-1)n!$ and let 
$A\in\R^{d\times D}$ be a matrix whose columns form a full spark frame, i.e., 
any subset of $d$ columns is linearly independent. Then the key $A$ is
universal and the induced map $\hbeta_A:\hRnd\rightarrow\R^{n\times D}$,
 $X\mapsto \beta_A(X)=\downarrow(XA)$ is injective. Furthermore, $\hbeta_A$
 is bi-Lipschitz, with estimates of the bi-Lipschitz constants 
 $a_0= \min_{J\subset[D],|J|=d}s_d(A[J])$ 
 and $b_0= s_1(A)$, where $s_1(A)$ denotes the largest singular value of $A$, 
 $A[J]$ denotes the submatrix of $A$ formed by columns indexed by $J$, and
  $s_d(A[J])$ denotes the $d^{th}$ singular value (in this case, the smallest)
   of $A[J]$.
 Specifically, for any $X,Y\in\Rnd$,
 \begin{equation}
     \label{eq:Lipbeta}
      a_0 d(\hat{X},\hat{Y})\leq \norm{\beta_A(X)-\beta_A(Y)}\leq b_0
      d(\hat{X},\hat{Y})
 \end{equation}
 where all norms are Frobenius norms.
\end{thm}
 
 {\bf Proof}
 
 Let $a_1,\ldots,a_D$ denote the columns of $A$,  $A=[a_1\vert\cdots\vert a_D]$. 

Fix $X,Y\in\Rnd$ two matrices. Then there are permutation matrices $P_0,\Pi_1,\ldots,\Pi_D,\Xi_1,\ldots,\Xi_D\in\Sn$ so that $d(\hat{X},\hat{Y})=\norm{P_0 X-Y}$ and
 \[
 \beta_A(X)=\left[
 \begin{array}{ccccc}
 \Pi_1 X a_1 & \vert & \cdots & \vert & \Pi_D X a_D
 \end{array} \right] 
 ~,~
 \beta_A(Y)=\left[
 \begin{array}{ccccc}
 \Xi_1 Y a_1 & \vert & \cdots & \vert & \Xi_D Y a_D
 \end{array} \right] 
 \]
 Thus
 \begin{equation}\label{eq:betaA} \norm{\beta_A(X)-\beta_A(Y)}^2 = \sum_{k=1}^D \norm{(\Pi_k X - \Xi_k Y) a_k}_2^2 =\sum_{k=1}^D \norm{(\Xi_k^T\Pi_k X-Y)a_k}_2^2
 \end{equation}
Permutations $\Pi_k$ and $\Xi_k$ satisfy the optimality condition:
$\norm{\Pi_k Xa_k - \Xi_k Ya_k}_2 = \min_{P\in\Sn} \norm{P Xa_k - Ya_k}_2$. 
Hence $\norm{\Pi_k Xa_k - \Xi_k Ya_k}_2 \leq \norm{P_0 Xa_k - Ya_k}_2$.
 Therefore:
 \begin{equation}
 \label{eq:betaAA}
     \norm{\beta_A(X)-\beta_A(Y)}^2 \leq
 \sum_{k=1}^D \norm{P_0 Xa_k-Ya_k}_2^2 = 
 \norm{(PX-Y)A}^2 \leq s_1(A)^2 \norm{P_0 X-Y}^2
 \end{equation} 
 from where we obtain the upper bound in (\ref{eq:Lipbeta}).
 
The lower bound in (\ref{eq:Lipbeta}) follows from the pigeonhole principle similar to the one employed in the proof of Theorem \ref{t2}. In equation (\ref{eq:betaA}) there are $D=1+(d-1)n!$ terms. Since only $n!$ permutations are distinct, there is a permutation $Q$ that repeats at least $d$ times. Say $J=\{j_1,j_2,\ldots,j_d\}\subset [D]$ is a set of indices so that $\Xi_{j_1}^T\Pi_{j_1}=\cdots=\Xi_{j_d}^T\Pi_{j_d}=Q$. Then
\begin{align*}
\norm{\beta_A(X)-\beta_A(Y)}^2 &\geq\sum_{k=1}^d \norm{(\Xi_{j_k}^T\Pi_{j_k}X-Y)a_{j_k}}_2^2 =
\norm{(QX-Y)A[J]}^2 \\
 &\geq s_d(A[J])^2\norm{QX-Y}^2
 \geq s_d(A[J])^2 \norm{P_0X-Y}^2
 \geq a_0^2 d(\hat{X},\hat{Y})^2.
\end{align*} 
The lower bound in (\ref{eq:Lipbeta}) implies that $\hbeta_A:\hRnd\rightarrow\R^{n\times D}$ is injective and hence $A$ is a universal key.
This ends the proof of Theorem \ref{t4}. $\Box$

\subsection{Bi-Lipschitz properties of universal keys}

In this subsection we prove that any universal key defines a bi-Lipschitz encoding map, regardless of $D$.

\begin{thm}
 \label{t5}
 Assume the key $A\in\R^{d\times D}$ is
universal, i.e., the induced map $\hbeta_A:\hRnd\rightarrow\R^{n\times D}$,
 $X\mapsto \beta_A(X)=\downarrow(XA)$ is injective. Then $\hbeta_A$ is bi-Lipschitz, that is, there are constants $a_0>0$ and $b_0>0$ so that
  for all $X,Y\in\Rnd$,
  \begin{equation}
     \label{eq:Lipbeta2}
      a_0\, d(\hat{X},\hat{Y})\leq \norm{\beta_A(X)-\beta_A(Y)}\leq b_0\,
      d(\hat{X},\hat{Y})
 \end{equation}
 where all are Frobenius norms.
  Furthermore, an estimate for $b_0$ is provided by the largest singular value of $A$, $b_0= s_1(A)$.
\end{thm}

{\bf Proof}

The upper bound in (\ref{eq:Lipbeta2}) follows as in the proof of Theorem \ref{t4}, from equations (\ref{eq:betaA}) and (\ref{eq:betaAA}). Notice that
no property is assumed  in order to obtain the upper Lipschitz bound.

The lower bound in (\ref{eq:Lipbeta2}) is more difficult. 
It is shown by contradiction following the strategy 
utilized in the 
Complex Phase Retrieval problem \cite{balazou}.

Assume $\inf_{X\not\sim  Y}\frac{\norm{\beta_A(X)-\beta_A(Y)}_2^2}{d(\hat{X},\hat{Y})^2}=0$. 

{\em Step 1: Reduction to local analysis.} 
Since $d(\hat{tX},\hat{tY})=t\,d(\hat{X},\hat{Y})$ for all $t>0$, the 
quotient $\frac{\norm{\beta_A(X)-\beta_A(Y)}_2}{d(\hat{X},\hat{Y})}$ 
is scale invariant. Therefore, there are sequences $(X^t)_t,(Y^t)_t$
with $\norm{Y^t}\leq\norm{X^t}=1$ and $d(\hat{X^t},\hat{Y^t})>0$ so that
$\lim_{t\rightarrow\infty} \frac{\norm{\beta_A(X^t)-\beta_A(Y^t)}_2}{d(\hat{X^t},\hat{Y^t})} = 0$.
By compactness of the closed unit ball, one can extract convergence subsequences. For easiness of notation, assume $(X^t)_t,(Y^t)_t$ are
these subsequences. Let $\Xinf=\lim_t X^t$ and $\Yinf = \lim_t Y^t$ denote their limits. Notice $\lim_t \norm{\beta_A(X^t)-\beta_A(Y^t)}_2=0$.
This implies $\norm{\beta_A(\Xinf)-\beta_A(\Yinf)}=0$ and thus $\beta_A(\Xinf)=\beta_A(\Yinf)$. Since $\widehat{\beta_A}$ is assumed injective, it follows that $\widehat{\Xinf}=\widehat{\Yinf}$. 

This means that, if the lower Lipschitz bound vanishes, then this 
is achieved by vanishing of a local lower Lipschitz bound. To follow the terminology in \cite{balazou}, the type I local lower Lipschitz bound vanishes at some 
$Z_0\in\Rnd$, with $\norm{Z_0}=1$:
\begin{equation}
    \label{eq:lb}
{A}(Z_0):= \lim_{r\rightarrow 0} \inf_{
\begin{array}{c}
\hat{X}\neq\hat{Y} \\
d(\hat{X},\hat{Z_0})<r \\
d(\hat{Y},\hat{Z_0})<r
\end{array}
}
\frac{\norm{\beta_A(X)-\beta_A(Y)}_2^2}{d(\hat{X},\hat{Y})^2} = 0.
\end{equation}
Note that, in general, the infimum of the type I local lower Lipschitz
bound over the unit sphere may be strictly larger than the global lower Lipschitz bound 
(see Theorems 2.1 and Theorem 2.2 in \cite{balazou} and Theorem 4.3 in \cite{balwan}). The compactness argument forces the local lower Lipschitz bound to vanish when the global lower bound vanishes.

\vspace{5mm}

{\em Step 2. Local Linearization.}  
The following stability subgroups of $\Sn$ play an important role:
\[ G = \{ P\in\Sn ~:~ PZ_0=Z_0\}~~,~~H_j=\{ P\in\Sn~:~PZ_0 a_j = Z_0 a_j\}~,~1\leq j\leq D. \]
Obviously $I_n\in G\subset H_j\subset \Sn$, for every $j\in[D]$. The group $G$ is the stabilizer of $Z_0$, whereas $H_j$ is the stabilizer of $Z_0a_j$. 
Let $\delta_0=\min_{P\in\Sn\setminus G} \norm{(I_n-P)Z_0}$ denote the smallest variation of $Z_0$ under row permutations. Note $\delta_0>0$ by the definition of $G$.

Consider $X=Z_0+U$ and $Y=Z_0+V$ where $U,V\in\Rnd$ are ``aligned" in the sense that $d(\hat{X},\hat{Y})=\norm{U-V}$. This property requires that $\norm{U-V}\leq\norm{PX-Y}$, for every $P\in\Sn$. 
Next result replaces equivalently this condition  
by requirements involving $(U,V)$ and the group $G$ only.
\begin{lem}
\label{l3.1}
Assume $\norm{U},\norm{V}<\frac{1}{4}\delta_0$, where $\delta_0=\min_{P\in\Sn\setminus G} \norm{(I_n-P)Z_0}$. Let $X=Z_0+U$, $Y=Z_0+V$. 
Then:
\begin{enumerate}
    \item $d(\hat{X},\hat{Z_0})=\norm{U}$ and $d(\hat{Y},\hat{Z_0})=\norm{V}$.
    \item $d(\hat{X},\hat{Y})=\min_{P\in G}\norm{U-PV}=\min_{P\in G}\norm{PU-V}$
    \item The following
are equivalent:
\begin{enumerate}
    \item $d(\hat{X},\hat{Y})=\norm{U-V}$.
    \item For every $P\in G$, $\norm{U-V}\leq \norm{PU-V}$.
    \item For every $P\in G$, $\ip{U}{V}\geq \ip{PU}{V}$.
\end{enumerate}
\end{enumerate}
\end{lem}
{\bf Proof of Lemma \ref{l3.1}}.
(1)

Note that is $U=0$ then the claim follows. Assume $U\neq 0$. Then
\[ d(\hat{X},\hat{Z_0})=\min_{P\in\Sn} \norm{X-PZ_0}
= \min_{P\in \Sn}\norm{(I_n-P)Z_0 + U}\leq \norm{U} \]
On the other hand, assume the minimum is achieved for a permutation $P_0\in\Sn$. If $P_0\in G$ then
$d(\hat{X},\hat{Z_0})=\norm{(I_n-P_0)Z_0+U}=\norm{U}$. If $P_0\not\in G$ then 
\[ d(\hat{X},\hat{Z_0})\geq \norm{(I_n-P_0)Z_0}-\norm{U}>\frac{3\delta_0}{4}>\norm{U}\geq d(\hat{X},\hat{Z_0}) \]
which yields a contradiction.
Hence $d(\hat{X},\hat{Z_0})=\norm{U}$. Similarly, one shows $d(\hat{X},\hat{Z_0})=\norm{V}$. 

(2) Clearly
\[ d(\hat{X},\hat{Y})=\min_{P\in\Sn}\norm{PX-Y}\leq \min_{P\in G}\norm{PX-Y}=\min_{P\in G}\norm{PU-V} \]
On the other hand, for $P\in\Sn\setminus G$ and $Q\in G$,
\[ \norm{PX-Y}=\norm{(P-I_n)Z_0 + PU-V}\geq \norm{(I_n-P)Z_0} - \norm{U}-\norm{V}\geq \]
\[ \geq \delta_0
-2\norm{U}-2\norm{V}+\norm{QU-V}\geq \min_{Q\in G}\norm{QU-V}\geq d(\hat{X},\hat{Y}). \]

(3)

(a)$\Rightarrow$(b).

If $d(\hat{X},\hat{Y})=\norm{U-V}$ then
\[ \norm{U-V}\leq \norm{PX-Y}=\norm{(P-I_n)Z_0 + PU-V}
~~,~~\forall P\in \Sn. \]
In particular, for $P\in G$, $(P-I_n)Z_0=0$ and
the above inequality reduces to (b).

(b)$\Rightarrow$(a).

Assume (b). For $P\in G$,
\[ \norm{U-V}=\norm{X-Y}\leq\norm{PU-V}=\norm{PX-Y}. \]
For $P\in\Sn\setminus G$,
\[ \norm{PX-Y}=\norm{(P-I_n)Z_0 + PU-V}\geq \norm{(I_n-P)Z_0} - \norm{U}-\norm{V}\geq \]
\[ \geq \delta_0
-2\norm{U}-2\norm{V}+\norm{U-V}\geq \norm{U-V}=\norm{X-Y}. \]
This shows $d(\hat{X},\hat{Y})=\norm{X-Y}=\norm{U-V}$.

(b)$\Longleftrightarrow$(c). This is immediate after squaring (b) and simplifying the terms.

$\Box$
\vspace{3mm}

Consider now sequences $(\hat{X^t})_t,(\hat{Y^t})_t$ that converge to $\hat{Z_0}$ 
and achieve lower bound 0 as in (\ref{eq:lb}).
Choose
representatives  $X_t$ and $Y_t$ in their equivalence classes that satisfy the hypothesis of Lemma \ref{l3.1} so that $X_t=Z_0+U_t$,  $Y_t=Z_0+V_t$, $\norm{U_t},\norm{V_y}<\frac{1}{4}\delta_0$,
 $d(\hat{X_t},\hat{Z_0})=\norm{U_t}$, $d(\hat{Y_t},\hat{Z_0})=\norm{V_t}$ and $d(\hat{X_t},\hat{Y_t})=\norm{U_t-V_t}>0$.
 With $A=[a_1 |\cdots|a_D]$ we obtain:
\[ \norm{\beta_A(X_t)-\beta_A(Y_t)}_2^2 = \sum_{j=1}^D \norm{\downarrow(X_t a_j)-\downarrow(Y_t a_j)}_2^2 =
\sum_{j=1}^D \norm{(Z_0+U_t)a_j-\Pi_{j,t}(Z_0+V_t)a_j}_2^2 \]
for some $\Pi_{j,t}\in\Sn$. In fact $\Pi_{j,t}\in argmin_{\Pi\in H_j}\norm{U_t-\Pi V_t)a_j}_2$. 
Pass to sub-sequences (that will be indexed by $t$ for an easier notation) so that $\Pi_{j,t}=\Pi_j$ for some $\Pi_j\in\Sn$. Thus
\[ \norm{\beta_A(X_t)-\beta_A(Y_t)}_2^2 =
\sum_{j=1}^D \norm{(I_n-\Pi_j)Z_0a_j + (U_t-\Pi_j V_t)a_j}_2^2 \]
Since the above sequence must converge to $0$ as $t\rightarrow\infty$, while $U_t,V_t\rightarrow 0$, it follows that necessarily $\Pi_j\in H_j$ and the
expressions simplify to
\[ \norm{\beta_A(X_t)-\beta_A(Y_t)}_2^2 =
\sum_{j=1}^D \norm{(U_t-\Pi_j V_t)a_j}_2^2  \]
Thus equation (\ref{eq:lb}) implies that for
every $j\in[D]$,
\begin{equation}
\label{eq:lb2}
\lim_{t\rightarrow\infty} \frac{\norm{(U_t-\Pi_j V_t)a_j}_2^2}{\norm{U_t-V_t}^2} = 0
\end{equation} 
where $\Pi_j\in H_j$, $\norm{U_t},\norm{V_t}\rightarrow 0$, and $U_t,V_t$  are aligned so that $\ip{U_t}{V_t}\geq \ip{PU_t}{V_t}$ for every $P\in G$.
Equivalently, relation (\ref{eq:lb}) can be restated as:
\begin{equation}
    \label{eq:opt2}
    \inf_{\begin{array}{c} U,V\in\Rnd \\ s.t. \\
    U\neq V \\
    \ip{U}{V}\geq \ip{PU}{V} , \forall P\in G
    \end{array} } \frac{\sum_{j=1}^D \norm{(U-\Pi_j V)a_j}_2^2}{\norm{U-V}^2} = 0
\end{equation}
for some permutations $\Pi_j\in H_j$, $j\in[D]$.
By Lemma \ref{l3.1} the constraint in the optimization problem above implies $\norm{U-V}=\min_{P\in G}\norm{U-PV}$. Hence (\ref{eq:opt2}) implies:
\begin{equation}
    \label{eq:opt3}
     \inf_{\begin{array}{c} U,V\in\Rnd \\ s.t. \\
    U\neq P V , \forall P\in G
    \end{array} }
    \max_{P\in G} \frac{\sum_{j=1}^D \norm{(U-\Pi_j V)a_j}_2^2}{\norm{U-P V}^2} = 0
\end{equation}
for same permutation matrices $\Pi_j$'s.
While the above optimization problem seems a relaxation of (\ref{eq:opt2}), in fact (\ref{eq:opt3}) implies (\ref{eq:opt2})
with a possibly change of permutation matrices $\Pi_j$, but 
remaining still in $H_j$.
\vspace{5mm}

{\em Step 3. Existence of a Minimizer.}

The optimization problem (\ref{eq:opt2}) is a Quadratically Constrained Ratio of Quadratics (QCRQ) optimization problem. A significant number of papers 
have been published on this topic \cite{teb06,teb10}. 
In particular, \cite{QCRQbook} presents
a formal setup for analysis of QCRQ problems. 
Our interest is to utilize some of these techniques in order to establish the existence of a minimizer for (\ref{eq:opt2}) or (\ref{eq:opt3}). Specifically we show:
\begin{lem}\label{l3.2}
Assume the key $A$ has linearly independent rows (equivalently, the columns of $A$ form a frame for $\R^d$) and the lower Lipschitz bound of $\hbeta_A$ is $0$. Then there are $\tilde{U},\tilde{V}\in\Rnd$ so that:
\begin{enumerate}
    \item $\tilde{U}\neq P \tilde{V}$, for every $P\in G$;
    \item For every $j\in[D]$, $(\tilde{U}-\Pi_j \tilde{V})a_j=0$.
\end{enumerate}
\end{lem}
{\bf Proof of Lemma \ref{l3.2}}

We start with the formulation (\ref{eq:opt3}). Therefore there are sequences 
$(U_t,V_t)_{t\geq 1}$ so that $U_t\neq PV_t$ for any $P\in G, t\geq 1$, and yet for any $P\in G$,
\[ \lim_{t\rightarrow\infty} \frac{\sum_{j=1}^D\norm{(U_t-\Pi_j V_t)a_j}_2^2}{\norm{U_t-P V_t}^2} = 0. \]
Let $E=\{(U,V)\in\Rnd\times\Rnd~,~(U-\Pi_j)V)a_j=0~,~\forall j\in[D]\}$ denote the null space of the linear operator 
\[ T:\Rnd\times\Rnd\rightarrow \R^D~,~(U,V)\mapsto \left[\begin{array}{ccccc}
(U-\Pi_1 V)a_1 & \vert & \cdots & \vert & (U-\Pi_D V)a_D
\end{array}\right],
\]
associated to the numerator of the above quotient. Let $F_P=\{(U,V)\in\Rnd\times\Rnd~,~U-PV=0\}$ be the null space of the linear operator 
\[ R_P:\Rnd\times\Rnd\rightarrow \Rnd~,~(U,V)\mapsto U-P V. \]
A consequence of (\ref{eq:opt3}) is that for every $P\in G$, $E\setminus F_P\neq\emptyset$. 
In particular, $F_p\cap E$ is a subspace of $E$ of positive codimension. Using the 
Baire category theorem (or more elementary linear algebra arguments), we conclude that
\[ E\setminus \left(\cup_{P\in G} F_P\right) \neq \emptyset. \]
Let $(\tilde{U},\tilde{V})\in E\setminus\left(\cup_{P\in G}F_P\right)$. This pair satisfies the
conclusions of Lemma \ref{l3.2}.

 $\Box$

\vspace{5mm}

{\em Step 4. Contradiction with the universality property of the key.}

So far we obtained that if the lower Lipschitz bound of $\hbeta_A$ vanishes than there are $Z_0,\tilde{U},\tilde{V}\in\Rnd$ with $Z_0\neq 0$ and $\tilde{U}\neq P \tilde{V}$, for all $P\in G$ that satisfy the conclusions of Lemma \ref{l3.2}. Notice $\ip{Z_0}{Z_0}=\ip{PZ_0}{Z_0}$
 for all $P\in G$ and $(Z_0-\Pi_j Z_0)a_j=0$ for all $j\in[D]$. Choose $s>0$ but small enough so that $s\norm{\tilde{U}},s\norm{\tilde{V}}<\frac{1}{4}\delta_0$ with $\delta_0=\min_{P\in\Sn\setminus G} \norm{(I_n-P)Z_0}$.
 Let $X=Z_0+s \tilde{U}$ and $Y=Z_0+s\tilde{V}$.
 Then Lemma \ref{l3.1} implies $d(\hat{X},\hat{Y})=\min_{P\in G}\norm{\tilde{U}-P\tilde{V}}>0$. 
 Hence $\hat{X}\neq\hat{Y}$. On the other hand,
 for every $j\in[D]$, $Xa_j = \Pi_j Ya_j$. Thus
  $\hbeta_A(\hat{X})=\hbeta_A(\hat{Y})$. 
  Contradiction with the assumption that $\hbeta_A$ is injective.
  
  This ends the proof of Theorem \ref{t5}.

$\Box$
\ignore{
\begin{rem}
The proof of the previous theorem provides estimates for 
both type I and type II local lower and upper Lipschitz bounds.
\end{rem}
}

\subsection{Dimension Reduction}

Theorem \ref{t4} provides an Euclidean bi-Lipschitz embedding of very high dimension, $D=1+(d-1)n!$. On the other hand, Theorem \ref{t5} shows that any universal key $A\in\R^{d\times D}$ for $\hRnd$, 
and hence any injective map $\hat{\beta}_A$ is bi-Lipschitz. In this subsection we show that 
any bi-Lipschitz Euclidean embedding $\hat{\beta}_A:\hRnd\rightarrow\R^{n\times D}$ with $D>2d$ 
can be further compressed to a smaller dimension space $\R^m$ with $m=2nd$ thus yielding
bi-Lipschitz Euclidean embeddings of redundancy 2. This is shown in the next result.

\begin{thm}
 \label{t6} Assume $A\in\R^{d\times D}$ is a universal key for $\hRnd$ with $D\geq 2d$. 
 Then, for $m\geq 2nd$, a generic linear operator $B:\R^{n\times D}\rightarrow\R^{m}$ with respect to Zariski topology on
 $\R^{n\times D\times m}$, the map
 \begin{equation}
 \label{eq:AB1}
     \hat{\beta}_{A,B}:\hRnd\rightarrow\R^{2nd}~,~ \hat{\beta}_{A,B}(\hat{X})=B\left(\hat{\beta}_A(\hat{X})\right)
 \end{equation}
 is bi-Lipschitz. In particular, almost every full-rank linear operator $B:\R^{n\times D}\rightarrow\R^{2nd}$ produces such a 
 bi-Lipschitz map.
\end{thm}

\begin{rmk}
The proof shows that, in fact, the complement set of linear operators $B$ that produce bi-Lipschitz embeddings is included 
in the zero-set of a polynomial. 
\end{rmk}

\begin{rmk}
Putting together Theorems \ref{t4}, \ref{t5}, \ref{t6} we obtain that the metric space $\hRnd$ admits
a global bi-Lipschitz embedding in the Euclidean space $\R^{2nd}$. This result is compatible
with a Whitney embedding theorem (see \S 1.3 in \cite{hirsh}) with the important caveat that the Whitney embedding result
applies to smooth manifolds, whereas here $\hRnd$ is merely a non-smooth algebraic variety.
\end{rmk}

\begin{rmk}
These three theorems are summarized in part two of the 
Theorem \ref{t2} presented in 
the first section.
\end{rmk}

\begin{rmk}
While the embedding dimension grows linearly in $nd$, in fact $m=2nd$, the computational complexity of constructing $\hbeta_{A,B}$ is NP due to the $1+(d-1)n!$ intermediary dimension.
\end{rmk}

\begin{rmk}
As the proofs show, for $D\geq 1+(d-1)n!$, a generic $(A,B)$ with
respect to Zariski topology, $A\in\R^{d\times D}$ and linear map $B:\R^{n\times D}\rightarrow\R^{2nd}$, produces a bi-Lipschitz embedding $(\hbeta_{A,B},d)$ of $\hRnd$ into $(\R^{2nd},\norm{\cdot}_2)$.  
\end{rmk}
{\bf Proof of Theorem \ref{t6} }

The proof follows a similar approach as in Theorem 3 of \cite{Cahill19}.
See also \cite{DUFRESNE20091979}.

Without loss of generality, assume $m<nD$.

Notice $\beta_A:\R^{n\times d}\rightarrow\R^{n\times D}$ is already homogeneous of degree 1 (with respect to positive scalars). 
Let $\Delta:\R^{n\times d}\times\R^{n\times d}\rightarrow\R^{n\times D}$ be defined by $\Delta(X,Y)=\beta_A(X)-\beta_A(Y)$. 
Denote $E=Ran(\Delta)=\{\beta_A(X)-\beta_A(Y)~,~X,Y\in\R^{n\times d} \}$. 

Recall $A=[a_1\vert\cdots\vert a_D]$ is a notation for the columns of key $A$. Notice that 
\[ \Delta(X,Y)=\left[P_1Xa_1-Q_1Ya_1\vert\cdots\vert P_DXa_D-Q_DYa_D\right]  \]
for some $P_1,\ldots,P_D,Q_1,\ldots,Q_D\in\Sn$, 
so that for each $k\in[D]$, $P_k,Q_k$ are permutations that 
sort monotone decreasingly vectors $Xa_k$ and $Ya_k$, respectively. 
In particular,
\[ E\subset F:=\bigcup_{\gamma \in(\Sn)^{2D}} F_\gamma~~,~~ F_\gamma:= Ran(L_{\gamma})
 \]
where the $(n!)^{2D}$ linear operators $L_{\gamma}:\R^{n\times d}\times\R^{n\times d}\rightarrow\R^{n\times D}$, are defined by
\[ L_{\gamma}(X,Y)=
\left[P_1Xa_1-Q_1Ya_1\vert\cdots\vert P_DXa_D-Q_DYa_D\right] \]
when $\gamma=(P_1,\ldots,P_D,Q_1,\ldots,Q_D)\in (\Sn)^{2D}$.

Claim:  We claim that, for $m\geq 2nd$ and
a generic linear operator $B:\R^{n\times D}\rightarrow\R^m$ we have
 $\ker(B)\cap F=\{0\}$. Such a generic linear operator has the kernel
 of dimension $\dim(\ker(B))=nD-m\leq n(D-2d)$. It is therefore
 sufficient to show that, for a generic subspace $V\subset \R^{n\times D}$
 of dimension $r\leq n(D-2d)$, for every $\gamma\in(\Sn)^{2D}$, 
 $V\cap F_\gamma =\{0\}$. This last claim follows from the
 observation $\dim(F_\gamma)\leq 2nd$.
 
We now show how this claim proves the Theorem. 
Let $B$ be such a linear map, and let $\beta_{A,B}:\R^{n\times d}\rightarrow\R^{m}$ be the map $\beta_{A,B}(X)=B(\downarrow(XA))$. Then $\beta_{A,B}(X)=\beta_{A,B}(Y)$ implies $\Delta(X,Y)=\beta_A(X)-\beta_A(Y)\in \ker(B)$. Thus $\Delta(X,Y)=0$
which implies $\beta_A(X)=\beta_A(Y)$. Since $\hbeta_A$ is injective on $\hRnd$ it follows $\hat{X}=\hat{Y}$. Thus $\hbeta_{A,B}$ is
injective. On the other hand, for each $\gamma=(P_1,\ldots,P_D,Q_1,\ldots,Q_D)\in\Sn^{2D}$, the restriction of $B$ to the linear space $Ran(L_{\gamma})$ is injective, and thus bounded below as a linear map: there is $a_\gamma>0$ so that for every $X,Y\in\R^{n\times d}$, 
$\norm{B(L_{\gamma}(X,Y))}\geq a_\gamma \norm{L_{\gamma}(X,Y)}$.
Let $a_\infty = \min_{\gamma}a_\gamma >0$. Thus
\[ \norm{\beta_{A,B}(X) - \beta_{A,B}(Y)} =\norm{B(L_{\gamma_0}(X,Y))}
\geq a_\infty \norm{L_{\gamma_0}(X,Y)}=a_\infty \norm{\beta_A(X)-\beta_A(Y)} \]
where $\gamma_0\in(S_n)^{2D}$ is a particular $2D$-tuple of permutations. This shows that 
$B{\vert}_{\beta_A(\R^{n\times d})}:\beta_A(\R^{n\times d})\rightarrow\R^m$ is bi-Lipschitz.
By Theorem \ref{t5}, the map $\hbeta_A$ is bi-Lipschitz. Therefore
we get $\hbeta_{A,B}$ is bi-Lipschitz as well.
$\Box$

\subsection{Proof of Corollary \ref{c0}\label{subsec4.4}}

(1) It is clear that any continuous $f$ induces a continuous $\varphi:\beta(\R^{n\times d})\rightarrow\R$ via $\varphi(\beta(X))=f(X)$. Furthermore, 
$F:=\beta(\R^{n\times d})=\hbeta(\hRnd)$
is a closed subset of $\R^m$ since $\hbeta$ is bi-Lipschitz. 
Then a consequence of Tietze extension theorem 
(see problem 8 in \S 12.1 of \cite{roydenfitzpatrick})
implies that $\varphi$ admits a continuous extension $g:\R^m\rightarrow\R$. Thus $g(\beta(X))=f(X)$ 
for all $X\in\R^{n\times d}$. The converse is trivial.

(2) As at part (1), the Lipschitz continuous function $f$ induces a Lipschitz continuous function $\varphi:F\rightarrow\R$. Since $F\subset\R^m$ is a subset of a Hilbert space, by Kirszbraun
extension theorem (see \cite{WelWil75}), $\varphi$ 
admits a Lipschitz continuous extension 
(even with the same Lipschitz constant!)
 $g:\R^m\rightarrow\R$ so that $g(\beta(X))=f(X)$ for every $X\in\R^{n\times d}$. The converse is trivial. $\Box$

\section{Applications to Graph Deep Learning\label{sec4}}

In this section we take an empirical look at the permutation invariant mappings presented in this paper. We focus on the problems of graph classification, for which we employ the PROTEINS\_FULL dataset \cite{DobsonDoing_proteins}, and graph regression, for which we employ the quantum chemistry QM9 dataset
\cite{ramakrishnan2014quantum}. In both problems we want to estimate a function $F: (A,Z) \rightarrow p$, where $(A,Z)$ characterizes a graph where $A \in \R^{n\times n}$ is an adjacency matrix and $Z \in \R^{n\times r}$ is an associated feature matrix where the $i^{th}$
row encodes an array of $r$ features associated with the $i^{th}$ node.  $p$ is a scalar output where we have $p \in \{0,1\}$ for binary classification and $p \in \R_+$ for regression.

\begin{figure}[!htbp]
    \centering
	\includegraphics[width=.8\linewidth]{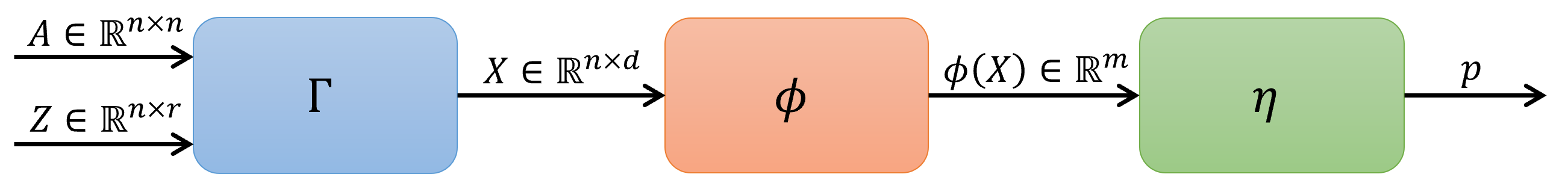}
	\caption[.]{.}
	\label{fig:gcn_end2end}
\end{figure}

We estimate $F$ using a deep network that is trained in a supervised manor. The network is comprised of three successive components applied in series: $\Gamma$, $\phi$, and $\zeta$. $\Gamma$ represents a graph deep network \cite{GCN},
which produces a set of embeddings $X \in \R^{N\times d}$ across the nodes in the graph. Here $N\geq n$ is chosen to accommodate the graph with the largest number of nodes. In this case, the last $N-n$ rows of $Y$ are filled with 0's. $\phi: \R^{N\times d} \rightarrow \R^{m}$ represents a permutation invariant mapping such as those proposed in this paper. $\zeta: \R^{m} \rightarrow \R$ is a fully connected neural network. The entire end-to-end network is shown in Figure \ref{fig:gcn_end2end}.

In this paper, we model $\Gamma$ using a Graph Convolutional Network (GCN) outlined in \cite{GCN}. 
Let $\DDD \in \mathbb{R}^{n \times n}$ be the associated degree matrix for our graph $\mathcal{G}$.  Also let $\tilde{A}$ be the associated adjacency matrix of $\mathcal{G}$ with added self connection: $\tilde{A}=I+A$, where $I$ is the $n \times n$ identity matrix, and $\tilde{\DDD}=\DDD+I$.  Finally, we define the modified adjacency matrix $\hat{A}=\tilde{\DDD}^{-1/2} \tilde{A} \tilde{\DDD}^{-1/2}$.  A GCN layer is defined as $H^{(l+1)}=\sigma (\hat{A}H^{(l-1)}W^{(l)})$.
Here $H^{(l-1)}$ represents the GCN state coming into the $l^{th}$ layer, $\sigma$ represents a chosen nonlinear element-by-element operation such as ReLU, and $W^{(l)}$ represents a matrix of trainable weights assigned to the $l^{th}$ layer whose number of rows match the number of columns in $H^{l}$ and number of columns is set to the size of the embeddings at the (l+1)’th layer.  The initial state $H^{(0)}$ of the network is set to the feature set of the nodes of the graph $H^{(0)}=Z$.

For $\phi$ we employ seven (7) different methods that are described next.
\begin{enumerate}
    \item ordering: For the ordering method, we set $D=d+1$, $\phi_{ordering}(X)=\beta_A(X)=\downarrow(XA)$ with $A=[I~1]$ the identity matrix followed by a column of ones. The ordering and identity-based mappings have the notable disadvantage of not producing the same output embedding size for different sized graphs. To accommodate this and have consistently sized inputs for $\eta$, we choose to zero-pad $\phi(X)$ for these methods to produce a vector in $\R^{m}$, where $m=ND=N(d+1)$ and $N$ is the size of the largest graph in the dataset.
    \item kernels: For the kernels method, 
    $$(\phi_{kernel}(X))_j=\sum_{k=1}^n K_G(x_k,a_j)=\sum_{k=1}^n exp(-\norm{x_k-a_j}^2),
    ~~j\in[m],$$ 
     for $X=[x_1|\cdots|x_n]^T$, where kernel vectors $a_1,\ldots,a_m\in\R^d$ 
     are generated randomly, each element of each vector is drawn from a standard normal distribution. Each resultant vector is then normalized to produce a kernel vector of magnitude one. When inputting the embedding $X$ to the kernels mapping, we first normalized the embedding for each respective node.
    \item identity: In this case $\phi_{id}(X)=X$, which is obviously not a permutation invariant map.
    \item data augmentation: In this case $\phi_{data\;augment}(X)=X$ but data augmentation is used. Our data augmentation scheme works as follows. We take the training set and create multiple permutations of the adjacency and associated feature matrix for each graph in the training set. We add each permuted graph to the training set to be included with the original graphs. In our experiments we use four added permutations for each graph when employing data augmentation.
    \item sum pooling: The sum pooling method sums the feature values across the set of nodes: $\phi_{sum\;pooling}(X)=\mathbf{1}_{n\times 1}^T X$.
    \item sort pooling: The sort pooling method flips entire rows of $X$ so that the last column is ordered descendingly, $\phi_{sort\;pool}(X)=\Pi X$ where $\Pi\in\Sn$ so that  $\Pi\,X(:,d)=\downarrow(X(:,d))$. 
    \item set-2-set: This method employs a recurrent neural network
    that achieves permutation invariance through attention-based weighted summations. It has been introduced in \cite{OrderMatters_2015arXiv151106391V}.
\end{enumerate}

For our deep neural network $\eta$ we use a simple multilayer perceptron of size described below.

Size parameters related to $\Gamma$ and $\zeta$ components are largely held constant across the different implementations.
However the network parameters are trained independently for each method.

\subsection{Graph Classification}

\subsubsection{Methodology}
For our experiments in graph classification we consider the PROTEINS\_FULL dataset obtained from \cite{KKMMN2016} and originally introduced in \cite{DobsonDoing_proteins}. 
The dataset consists of 1113 proteins falling into one of two classes: those that function as enzymes and those that do not. Across the dataset there are 450 enzymes in total. 
The graph for each protein is constructed such that the nodes represent amino acids and the edges represent the bonds between them. The number of amino acids (nodes) vary from around 20 to a maximum of 620 per protein with an average of 39.06. 
Each protein comes with a set of features for each node. 
The features represent characteristics of the associated amino acid represented by the node. The number of features is $r=29$.
We run the end-to-end model with three GCN layers in $\Gamma$, each with 50 hidden units. 
$\zeta$ consists of three dense multi-layer perceptron layers, each with 150 hidden units. 
We set d equal to 1, 10, 50 and 100.

For each method and embedding size we train for 300 epochs. Note though that the data augmentation method will have experienced five times as many training steps due to the increased size of its training set. We use a batch size of 128 graphs. The loss function minimized during training is the binary cross entropy loss (BCE) defined as
\begin{equation}
\label{eq:BCE}
BCE = -\frac{1}{B}\sum_{t=1}^B p_t log(\sigma(\eta(\phi(X^{(t)}))))+(1-p_t)log(1-\sigma(\eta(\phi(X^{(t)}))))    
\end{equation} 
where $B=128$ is the batch size, $p_t=1$ when the $t^{th}$ graph
(protein) is an enzyme and $p_t=0$ otherwise, $\sigma(x)=\frac{1}{1+e^{-x}}$ is the sigmoid function that maps the output $\eta(\phi(X^{(t)})$ of the 3-layer fully connected network $\eta$ to $[0,1]$. Three performance metrics were computed: accuracy (ACC), area under the receiver operating characteristic curve (AUC), and average precision (AP) as area under the precision-recall curve from precision scores. These measure are defined as follows (see   sklearn.metrics module documentation in pytorch, or \cite{roc}).

For a threshold $\tau\in[0,1]$, the classification decision $\hat{p}_t(\tau)$ is given by:
\begin{equation}
    \hat{p}_t(\tau) = \left\{
    \begin{array}{rcl}
    1 & \mbox{if} & \mbox{$\sigma(\eta(\phi(X^{(t)}))\geq \tau$} \\
    0 & \mbox{if} & \mbox{otherwise}
    \end{array}\right. .
\end{equation}
By default $\tau=\frac{1}{2}$. For a given threshold, one computes the four scores, true positive (TP), false positive (FP), true negative (TN) and false negative (FN):
\begin{equation}
    TP(\tau) = \frac{1}{B_1}\sum_{t=1}^B 1_{\hat{p}_t(\tau) = 1}1_{p_t = 1} ~~, ~~
    TN(\tau) = \frac{1}{B_0} 1_{\hat{p}_t(\tau) = 0}1_{p_t = 0}
\end{equation}
\begin{equation}
    FP(\tau) = \frac{1}{B_0}\sum_{t=1}^B 1_{\hat{p}_t(\tau) = 1}1_{p_t = 0} = 1-TN(\tau) ~~,~~
    FN(\tau) = \frac{1}{B_1}\sum_{t=1}^B 1_{\hat{p}_t(\tau) = 0}1_{p_t = 1} = 1-TP(\tau)
\end{equation}
where $B_0=\sum_{t=1}^B 1_{p_t = 0}$ and $B_1=\sum_{t=1}^B 1_{p_t = 1}=B-B_0$.

These four statistics predict Precision $P(\tau)$, Recall $R(\tau)$ (also known as sensitivity or true positive rate), and Specificity $S(\tau)$ (also known as true negative rate)
\begin{equation} P(\tau) = \frac{TP(\tau)}{TP(\tau)+FP(\tau)}
~~,~~R(\tau) =  \frac{TP(\tau)}{TP(\tau)+FN(\tau)}
~~,~~S(\tau) = \frac{TN(\tau)}{TN(\tau)+FP(\tau)}
\end{equation}

Accuracy (ACC) is defined as the fraction of correct classification for default threshold $\tau=\frac{1}{2}$ over the set of batch samples:
\begin{equation}
\label{eq:ACC}
    ACC = \frac{1}{B}\sum_{t=1}^B 1_{p_t = \hat{p}_t(\frac{1}{2})} =\frac{B_0}{B} TN(\frac{1}{2}) + \frac{B_1}{B} TP(\frac{1}{2}) 
\end{equation}
Area under the receiver operating characteristic curve (AUC) is computed from prediction scores as the area under true positive rate (TPR) vs. false positive rate (FPR) curve, i.e. the recall vs. 1-specificity curve
\begin{equation}
\label{eq:AUC}
    AUC = \frac{1}{2}\sum_{k=1}^K (S(\tau_{k-1})-S(\tau_k))(R(\tau_{k-1})+R(\tau_k))
\end{equation}
where $K$ is the number of thresholds.
Average precision (AP) summarizes a precision-recall curve as the weighted mean of precision achieved at each threshold, with the increase in recall from the previous thresholds used as the weight:
\begin{equation}
\label{eq:AP}
    AP = \sum_{k=1}^K (R(\tau_k) - R(\tau_{k-1}))P(\tau_k).
\end{equation}

We track the binary cross entropy (BCE) through training and we compute it on the holdout set and a random node permutation of the holdout set (see Figures \ref{fig:prot1} and \ref{fig:prot2}). The lower the value the better.

We look at the three performance metrics on the training set, the holdout set, and a random node permutation of the holdout set: see Figures \ref{fig:prot3}, and \ref{fig:prot4} for accuracy (ACC); see Figures \ref{fig:prot5}, and \ref{fig:prot6} for area under the receiver operating characteristic curve (AUC); and see Figures \ref{fig:prot7}, and \ref{fig:prot8} for average precision (AP). For all these performance metrics, the higher the score the better.

\subsubsection{Discussion}

Tables \ref{table:t1}-\ref{table:t12} list values of the three performance metrics (ACC, AUC, AP) at the end of training (after 300 epochs). 
Performances over the course of training are plotted in Figures \ref{fig:prot1} through \ref{fig:prot8}.

The authors of \cite{KKMMN2016} utilized a Support Vector Machine (1-layer perceptron) for classification and 
obtained an accuracy (ACC) of 77\% on the entire data set 
using 52 features, and an accuracy of 80\% on a smaller set of 36 features. By comparison, our data augmentation method for $d=100$ achieved an accuracy of 97.5\% on training data set,
but dropped dramatically to 73\% on holdout data, and 72\% on 
holdout data set with randomly permuted nodes. 
On the other hand, both the kernels method and the sum-pooling
method with $d=50$ achieved an accuracy of around 79\% on
training data set, while dropping accuracy performance by 
only 2\% to around
77\% on holdout data (as well as holdout data with nodes permuted).

For $d=1$, data augmentation performed the best on the training set with an area under the receiver operating characteristic  (AUC) of 0.896, followed closely by the identity method with an AUC of 0.886. On the permuted holdout set however, sort-pooling performed the best with an AUC of 0.803.

For $d=10$, sum-pooling, ordering, and kernels performed well on the permuted holdout set with AUC's of 0.821, 0.820, and 0.818 respectively. The high performance of the identity method, data augmentation, and sort-pooling on the training set did not translate to the permuted holdout set at $d=10$. By $d=100$, sum-pooling still performed the best on the permuted holdout set with an AUC of 0.817. This was followed by the kernels method which achieved an AUC of 0.801 on the permuted holdout set.

For experiments where $d>1$, the identity method and data augmentation show a notable drop in performance from the training set to the holdout set. This trend is also, to a lesser extent, visible in the sort pooling and ordering methods. In the holdout permuted set we see significant oscillations in the performance of both the identity and data augmentation methods.

\subsection{Graph Regression}

\subsubsection{Methodology}
For our experiments in graph regression we consider the qm9 dataset \cite{ramakrishnan2014quantum}. This dataset consists of 134 thousand molecules represented as graphs, where the nodes represent atoms and edges represent the bonds between them. 

Each graph has between 3 and 29 nodes, $3\leq n\leq 29$. Each node has 11 features,  $r=11$. We hold out 20 thousand of these molecules for evaluation purposes. The dataset includes 19 quantitative features for each molecule.

For the purposes of our study, we focus on electron energy gap (units $eV$), which is $\Delta\varepsilon$ in \cite{DFTpaper} whose chemical accuracy is $0.043 eV$ and whose prediction performance of any machine learning technique
is worse than any other feature.
The best existing estimator for this feature is enn-s2s-ens5 from \cite{Gilmer_2017arXiv170401212G}
 and has a mean absolute error (MAE) of $0.0529eV$ which is $1.23$ larger than the chemical accuracy. 
 We run the end to end model with three GCN layers in $\Gamma$, each with 50 hidden units. $\eta$ consists of three multi-layer perceptron layers, each with 150 hidden units. We use rectified linear units as our nonlinear activation function. Finally, we vary $d$, the size of the node embeddings that are outputted by $\Gamma$. We set $d$ equal to 1, 10, 50 and 100.

For each method and embedding size we train for 300 epochs. Note though that the data augmentation method will have experienced five times as many training steps due to the increased size of its training set. We use a batch size of 128 graphs. The loss function minimized during training is the mean square error (MSE) between the ground truth and the network output 
 (see Figures \ref{fig:a1}, \ref{fig:a2}) 
\begin{equation}
\label{eq:MSE}
MSE = \frac{1}{B}\sum_{t=1}^B |\Delta\varepsilon_t -\eta(\phi(X^{(t)}))))|^2
\end{equation}
where $B=128$ is the batch size of 128 graphs and $\Delta\varepsilon_t$ is the electron energy gap of the $t^{th}$ graph (molecule). The performance
metric is Mean Absolute Error (MAE)
\begin{equation}
\label{eq:MAE}
MAE = \frac{1}{B}\sum_{t=1}^B |\Delta\varepsilon_t -\eta(\phi(X^{(t)}))))|.
\end{equation}
We track the mean absolute error through the course of training. We look at this performance metric on the training set, the holdout set, and a random node permutation of the holdout set (see Figures \ref{fig:b1}, and \ref{fig:b2}).

\subsubsection{Discussion}

Numerical results at the end of training (after 300 epochs) are included in Tables \ref{table:a}, \ref{table:b}, \ref{table:c} and \ref{table:d}.
From the results we see that the ordering method performed best for $d=100$
followed closely by the data augmentation method, while both the ordering method and the kernels method performed well for $d=10$, though both fell slightly short of data augmentation which performed marginally better on both the training data and the holdout data, though with significantly more training iterations. For $d=1$, the kernels method failed to train adequately. The identity mapping performed relatively well on training data (for $d=100$ it achieved the smallest MAE among all methods and all parameters) and even the holdout data, however it lost its performance on the permuted holdout data. The identity mapping's failure to generalize across permutations of the holdout set is likely exacerbated by the fact that the QM9 data as presented to the network comes ordered in its node positions from heaviest atom to lightest. Data augmentation notably kept its performance despite this due to training on many permutations of the data. 

For $d=100$, our ordering method achieved a MAE of $0.155eV$ on training data set and $0.187eV$ on holdout data set, which are $3.6$ and $4.35$ times larger than the chemical accuracy ($0.043eV$\ignore{ cf. Supplementary material of \cite{Gilmer_2017arXiv170401212G}}), respectively. This is worse than the enn-s2s-ens5 method in \cite{Gilmer_2017arXiv170401212G} (current best method) that achieved a MAE $0.0529$ (eV), $1.23$ larger than the chemical accuracy, 
but better than the Coulomb Matrix (CM) representation in \cite{PhysRevLett.108.058301} that achieved a  MAE $5.32$ larger than the chemical accuracy  whose features were optimized for this task.

\bibliographystyle{amsplain}
\bibliography{main.bib}
\newpage
\appendix

\section{Results for the PROTEINS\_FULL dataset}


\begin{table}[hbtp]
	\scalebox{1.0}{
	\begin{tabular}{|w{c}{2cm}|m{1.3cm} m{1.3cm} m{1.3cm} m{1.4cm} m{1.3cm} m{1.3cm} m{1.3cm}|}
		\hline
		d = 1 & ordering & kernels & identity & data augment & sum-pooling & sort-pooling & set-2-set  \\
		\hline
		Training & 76 & 72 & 80 & 81.6 & 76.2 & 78 & 72.4 \\
		\hline
		Holdout & 74 & 74 & 72.5 & 76.5 & 70.5 & 74.5 & 72 \\
		\hline
        Holdout Perm & 74 & 74 & 67.5 & 75 & 70.5 & 74.5 & 72 \\
		\hline
	\end{tabular}} 
	\caption{Accuracy ACC(\%) for enzyme/non-enzyme classification of the seven algorithms on PROTEINS\_FULL dataset after 300 epochs for embedding dimension $d=1$}
	\label{table:t1}
\end{table}

\begin{table}[hbtp]
	\scalebox{1.0}{
	\begin{tabular}{|w{c}{2cm}|m{1.3cm} m{1.3cm} m{1.3cm} m{1.4cm} m{1.3cm} m{1.3cm} m{1.3cm}|}
		\hline
		d = 10 & ordering & kernels & identity & data augment & sum-pooling & sort-pooling & set-2-set  \\
		\hline
		Training & 84.5 & 78.2 & 87 & 90.6 & 77.8 & 85.2 & 72.5 \\
		\hline
		Holdout & 74 & 75.5 & 73 & 76 & 75 & 71 & 74.5 \\
		\hline
        Holdout Perm & 74 & 75.5 & 62.5 & 73.5 & 75 & 71 & 74.5 \\
		\hline
	\end{tabular}} 
	\caption{Accuracy ACC(\%) for enzyme/non-enzyme classification of the seven algorithms on PROTEINS\_FULL dataset after 300 epochs for embedding dimension $d=10$}
	\label{table:t2}
\end{table}

\begin{table}[hbtp]
	\scalebox{1.0}{
	\begin{tabular}{|w{c}{2cm}|m{1.3cm} m{1.3cm} m{1.3cm} m{1.4cm} m{1.3cm} m{1.3cm} m{1.3cm}|}
		\hline
		d = 50 & ordering & kernels & identity & data augment & sum-pooling & sort-pooling & set-2-set  \\
		\hline
		Training & 83.1 & 78.8 & 91 & 96 & 79.2 & 83.7 & 76.7\\
		\hline
		Holdout & 71.5 & 76.5 & 72.5 & 71 & 77 & 71 & 76 \\
		\hline
        Holdout Perm & 71.5 & 76.5 & 69.5 & 72 & 77 & 71 & 76\\
		\hline
	\end{tabular}} 
	\caption{Accuracy ACC(\%) for enzyme/non-enzyme classification of the seven algorithms on PROTEINS\_FULL dataset after 300 epochs for embedding dimension $d=50$}
	\label{table:t3}
\end{table}

\begin{table}[hbtp]
	\scalebox{1.0}{
	\begin{tabular}{|w{c}{2cm}|m{1.3cm} m{1.3cm} m{1.3cm} m{1.4cm} m{1.3cm} m{1.3cm} m{1.3cm}|}
		\hline
		d = 100 & ordering & kernels & identity & data augment & sum-pooling & sort-pooling & set-2-set  \\
		\hline
		Training & 88 & 77 & 97.5 & 97.5 & 78.1 & 87.3 & 76.6 \\
		\hline
		Holdout & 71 & 74.5 & 72.5 & 73 & 75.5 & 69.5 & 74.5 \\
		\hline
        Holdout Perm & 71 & 74.5 & 68.5 & 72 & 75.5 & 69.5 & 74.5 \\
		\hline
	\end{tabular}} 
	\caption{Accuracy ACC(\%) for enzyme/non-enzyme classification of the seven algorithms on PROTEINS\_FULL dataset after 300 epochs for embedding dimension $d=100$}
	\label{table:t4}
\end{table}

\begin{table}[hbtp]
	\scalebox{1.0}{
	\begin{tabular}{|w{c}{2cm}|m{1.3cm} m{1.3cm} m{1.3cm} m{1.4cm} m{1.3cm} m{1.3cm} m{1.3cm}|}
		\hline
		d = 1 & ordering & kernels & identity & data augment & sum-pooling & sort-pooling & set-2-set  \\
		\hline
		Training & 0.846 & 0.758 & 0.886 & 0.896 & 0.818 & 0.858 & 0.778\\
		\hline
		Holdout & 0.794 &  0.775 & 0.766 & 0.796 & 0.777 & 0.803 & 0.788\\
		\hline
        Holdout Perm & 0.794 & 0.775 & 0.747 & 0.785 & 0.777 & 0.803 & 0.788\\
		\hline
	\end{tabular}} 
	\caption{Area under the receiver operating characteristic curve (AUC)  for enzyme/non-enzyme classification of the seven algorithms on PROTEINS\_FULL dataset after 300 epochs for embedding dimension $d=1$}
	\label{table:t5}
\end{table}

\begin{table}[hbtp]
	\scalebox{1.0}{
	\begin{tabular}{|w{c}{2cm}|m{1.3cm} m{1.3cm} m{1.3cm} m{1.4cm} m{1.3cm} m{1.3cm} m{1.3cm}|}
		\hline
		d = 10 & ordering & kernels & identity & data augment & sum-pooling & sort-pooling & set-2-set  \\
		\hline
		Training & 0.913 & 0.849 & 0.941 & 0.970 & 0.842 & 0.930 & 0.787 \\
		\hline
		Holdout & 0.820 & 0.817 & 0.782 & 0.796 & 0.821 & 0.798 & 0.779 \\
		\hline
        Holdout Perm & 0.820 & 0.817 & 0.668 & 0.784 & 0.821 & 0.798 & 0.779 \\
		\hline
	\end{tabular}} 
	\caption{Area under the receiver operating characteristic curve (AUC) for enzyme/non-enzyme classification of the seven algorithms on PROTEINS\_FULL dataset after 300 epochs for embedding dimension $d=10$}
	\label{table:t6}
\end{table}

\begin{table}[hbtp]
	\scalebox{1.0}{
	\begin{tabular}{|w{c}{2cm}|m{1.3cm} m{1.3cm} m{1.3cm} m{1.4cm} m{1.3cm} m{1.3cm} m{1.3cm}|}
		\hline
		d = 50 & ordering & kernels & identity & data augment & sum-pooling & sort-pooling & set-2-set  \\
		\hline
		Training & 0.922 & 0.847 & 0.965 & 0.994 & 0.856 & 0.920 & 0.820\\
		\hline
		Holdout & 0.791 & 0.818 & 0.775 & 0.768 & 0.821 & 0.791 & 0.777\\
		\hline
        Holdout Perm & 0.791 & 0.818 & 0.716 & 0.768 & 0.821 & 0.791 & 0.777\\
		\hline
	\end{tabular}} 
	\caption{Area under the receiver operating characteristic curve (AUC) for enzyme/non-enzyme classification of the seven algorithms on PROTEINS\_FULL dataset after 300 epochs for embedding dimension $d=50$}
	\label{table:t7}
\end{table}

\begin{table}[hbtp]
	\scalebox{1.0}{
	\begin{tabular}{|w{c}{2cm}|m{1.3cm} m{1.3cm} m{1.3cm} m{1.4cm} m{1.3cm} m{1.3cm} m{1.3cm}|}
		\hline
		d = 100 & ordering & kernels & identity & data augment & sum-pooling & sort-pooling & set-2-set  \\
		\hline
		Training & 0.949 & 0.832 & 0.997 & 0.997 & 0.849 & 0.948 & 0.842 \\
		\hline
		Holdout & 0.754 & 0.801 & 0.766 & 0.775 & 0.817 & 0.784 & 0.776 \\
		\hline
        Holdout Perm & 0.754 & 0.801 & 0.708 & 0.775 & 0.817 & 0.784 & 0.776 \\
		\hline
	\end{tabular}} 
	\caption{Area under the receiver operating characteristic curve (AUC) for enzyme/non-enzyme classification of the seven algorithms on PROTEINS\_FULL dataset after 300 epochs for embedding dimension $d=100$}
	\label{table:t8}
\end{table}


\begin{table}[hbtp]
	\scalebox{1.0}{
	\begin{tabular}{|w{c}{2cm}|m{1.3cm} m{1.3cm} m{1.3cm} m{1.4cm} m{1.3cm} m{1.3cm} m{1.3cm}|}
		\hline
		d = 1 & ordering & kernels & identity & data augment & sum-pooling & sort-pooling & set-2-set  \\
		\hline
		Training & 0.788 & 0.709 & 0.844 & 0.857 & 0.754 & 0.811 & 0.707\\
		\hline
		Holdout & 0.720 & 0.698 & 0.692 & 0.725 & 0.636 & 0.680 & 0.708\\
		\hline
        Holdout Perm & 0.720 & 0.698 & 0.622 & 0.710 & 0.636 & 0.680 & 0.708\\
		\hline
	\end{tabular}} 
	\caption{Average precision (AP) for enzyme/non-enzyme classification of the seven algorithms on PROTEINS\_FULL dataset after 300 epochs for embedding dimension $d=1$}
	\label{table:t9}
\end{table}

\begin{table}[hbtp]
	\scalebox{1.0}{
	\begin{tabular}{|w{c}{2cm}|m{1.3cm} m{1.3cm} m{1.3cm} m{1.4cm} m{1.3cm} m{1.3cm} m{1.3cm}|}
		\hline
		d = 10 & ordering & kernels & identity & data augment & sum-pooling & sort-pooling & set-2-set  \\
		\hline
		Training & 0.890 & 0.804 & 0.922 & 0.961 & 0.797 & 0.904 & 0.722 \\
		\hline
		Holdout & 0.738 & 0.749 & 0.631 & 0.646 & 0.753 & 0.693 & 0.693 \\
		\hline
        Holdout Perm & 0.738 & 0.749 & 0.497 & 0.664 & 0.753 & 0.693 & 0.693 \\
		\hline
	\end{tabular}} 
	\caption{Average precision (AP) for enzyme/non-enzyme classification of the seven algorithms on PROTEINS\_FULL dataset after 300 epochs for embedding dimension $d=10$}
	\label{table:t10}
\end{table}

\begin{table}[hbtp]
	\scalebox{1.0}{
	\begin{tabular}{|w{c}{2cm}|m{1.3cm} m{1.3cm} m{1.3cm} m{1.4cm} m{1.3cm} m{1.3cm} m{1.3cm}|}
		\hline
		d = 50 & ordering & kernels & identity & data augment & sum-pooling & sort-pooling & set-2-set  \\
		\hline
		Training & 0.899 & 0.797 & 0.950 & 0.991 & 0.814 & 0.891 & 0.757\\
		\hline
		Holdout & 0.700 & 0.738 & 0.627 & 0.589 & 0.750 & 0.676 & 0.666 \\
		\hline
        Holdout Perm & 0.700 & 0.738 & 0.520 & 0.600 & 0.750 & 0.676 & 0.666 \\
		\hline
	\end{tabular}} 
	\caption{Average precision (AP) for enzyme/non-enzyme classification of the seven algorithms on PROTEINS\_FULL dataset after 300 epochs for embedding dimension $d=50$}
	\label{table:t11}
\end{table}

\begin{table}[hbtp]
	\scalebox{1.0}{
	\begin{tabular}{|w{c}{2cm}|m{1.3cm} m{1.3cm} m{1.3cm} m{1.4cm} m{1.3cm} m{1.3cm} m{1.3cm}|}
		\hline
		d = 100 & ordering & kernels & identity & data augment & sum-pooling & sort-pooling & set-2-set  \\
		\hline
		Training & 0.933 & 0.777 & 0.995 & 0.995 & 0.806 & 0.927 & 0.782 \\
		\hline
		Holdout & 0.627 & 0.729 & 0.601 & 0.622 & 0.747 & 0.637 & 0.704 \\
		\hline
        Holdout Perm & 0.627 & 0.729 & 0.529 & 0.656 & 0.747 & 0.637 & 0.704 \\
		\hline
	\end{tabular}} 
	\caption{Average precision (AP) for enzyme/non-enzyme classification of the seven algorithms on PROTEINS\_FULL dataset after 300 epochs for embedding dimension $d=100$}
	\label{table:t12}
\end{table}


\begin{figure}[p]
	\includegraphics[width=0.49\linewidth]{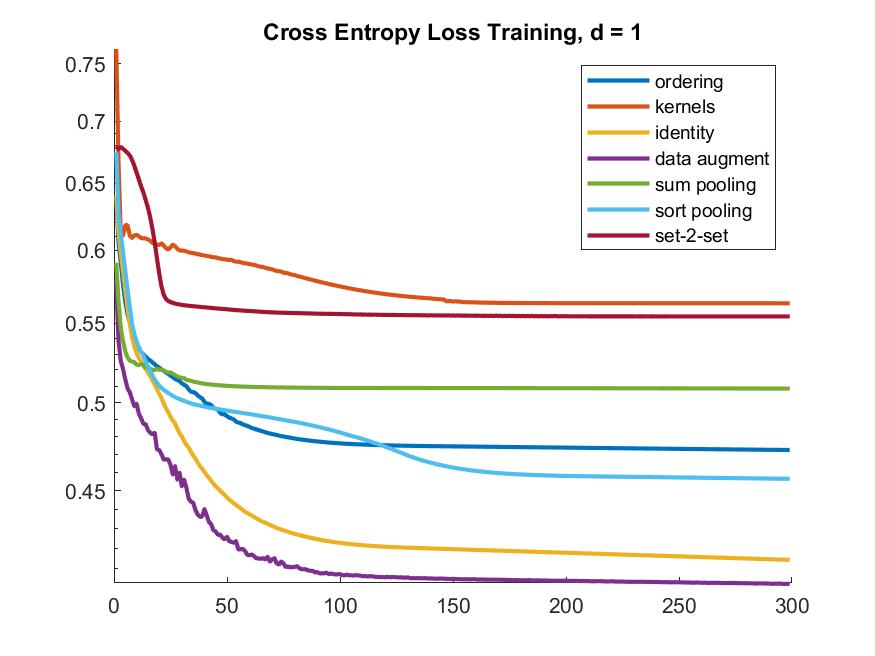}
	\includegraphics[width=0.49\linewidth]{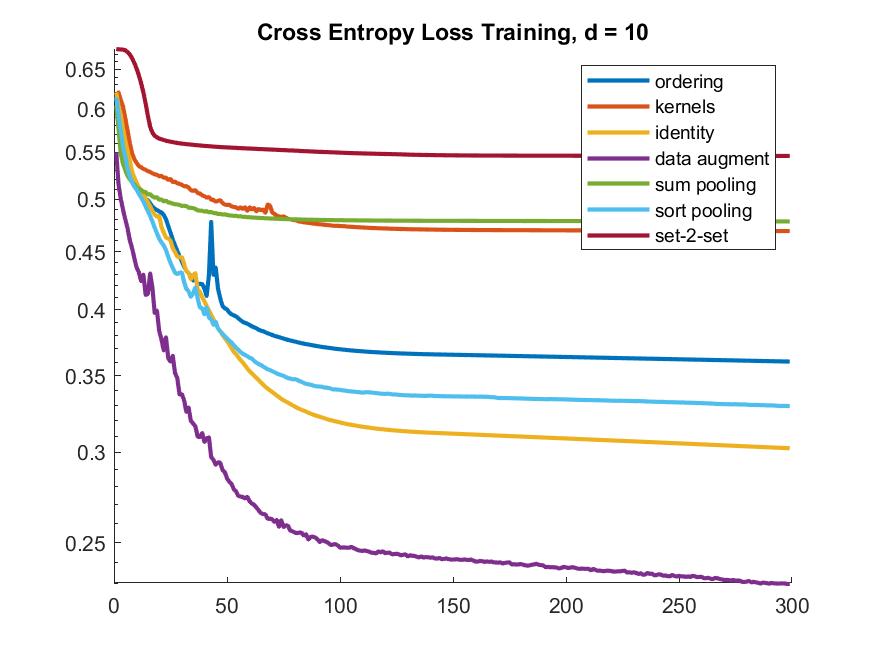}
	\\
	\includegraphics[width=0.49\linewidth]{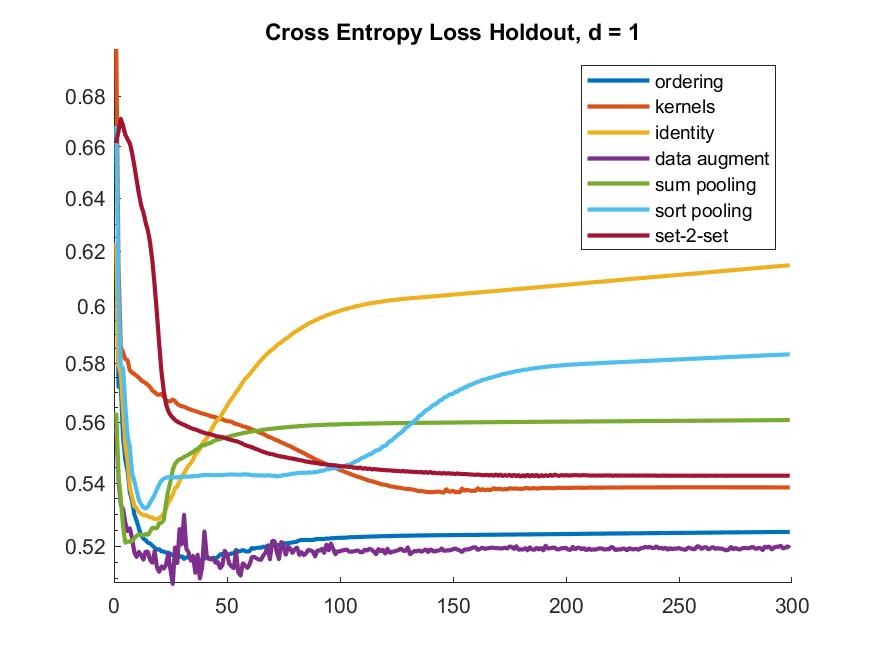}
	\includegraphics[width=0.49\linewidth]{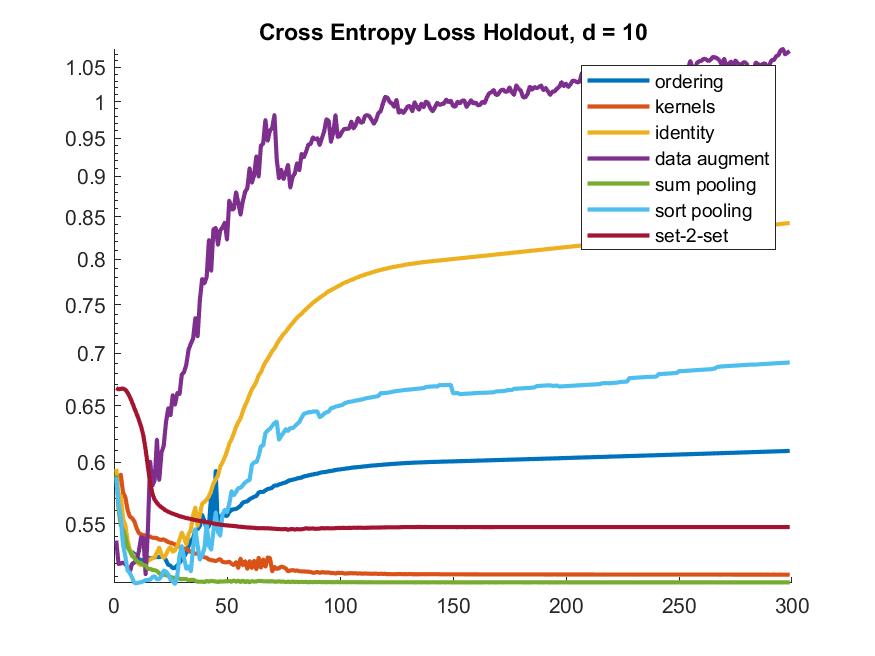}
    \\
	\includegraphics[width=0.49\linewidth]{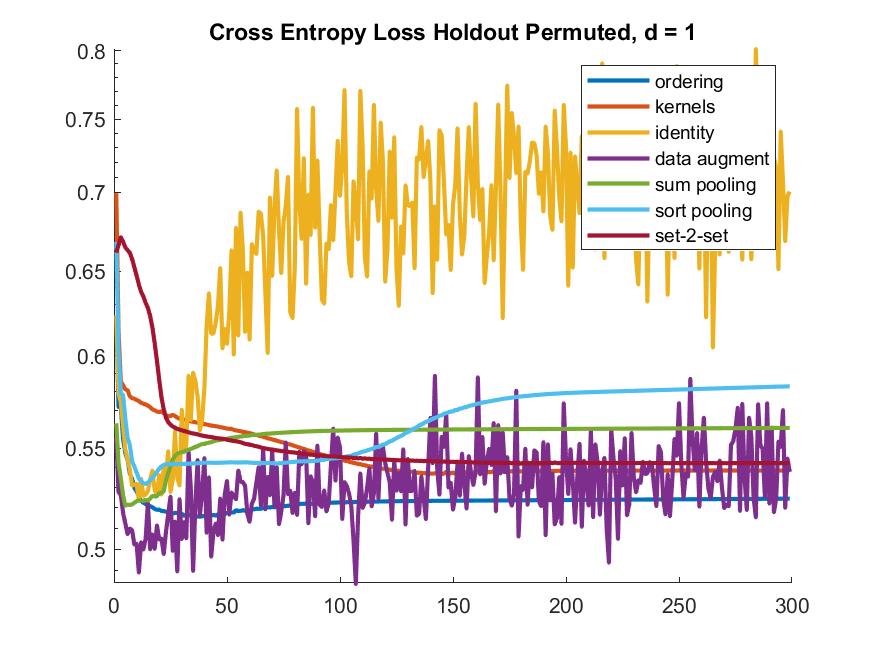}
	\includegraphics[width=0.49\linewidth]{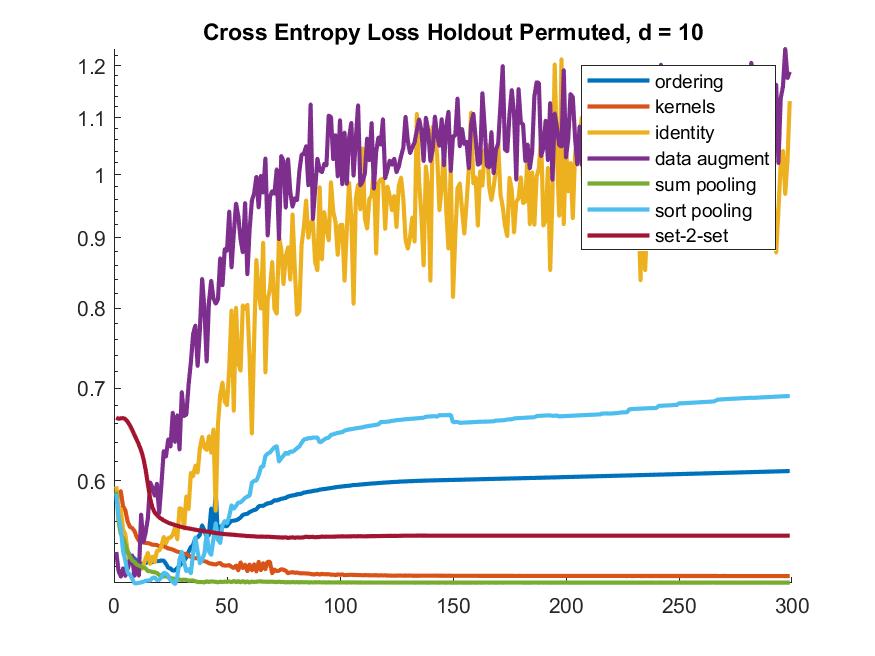}
	\caption[caption me]{Binary Cross Entropy  (BCE) for enzyme/non-enzyme classification on PROTEINS\_FULL dataset using; $d=1$ left column; $d=10$ right column.}
	\label{fig:prot1}
\end{figure}

\begin{figure}[p]
	\includegraphics[width=0.49\linewidth]{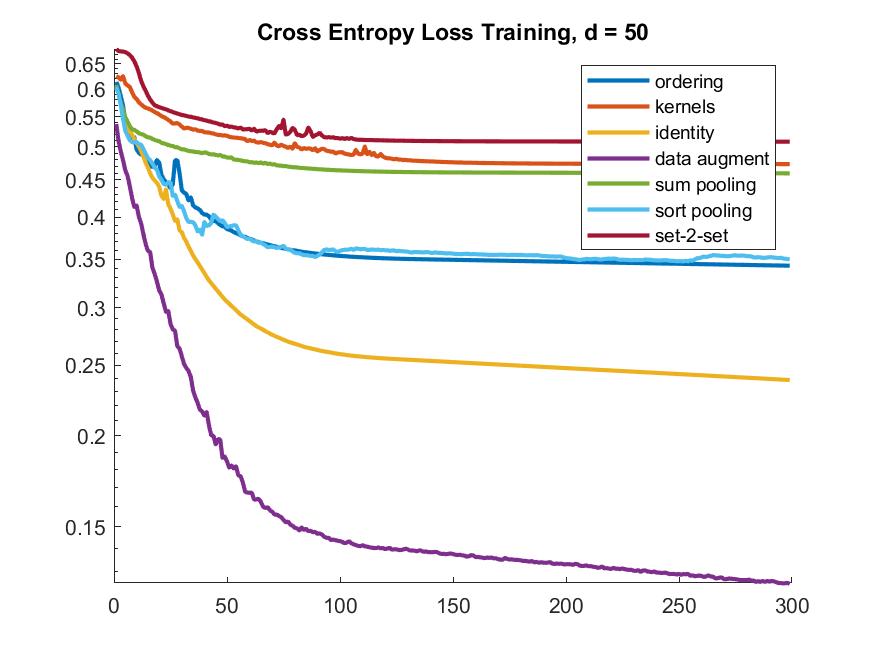}
	\includegraphics[width=0.49\linewidth]{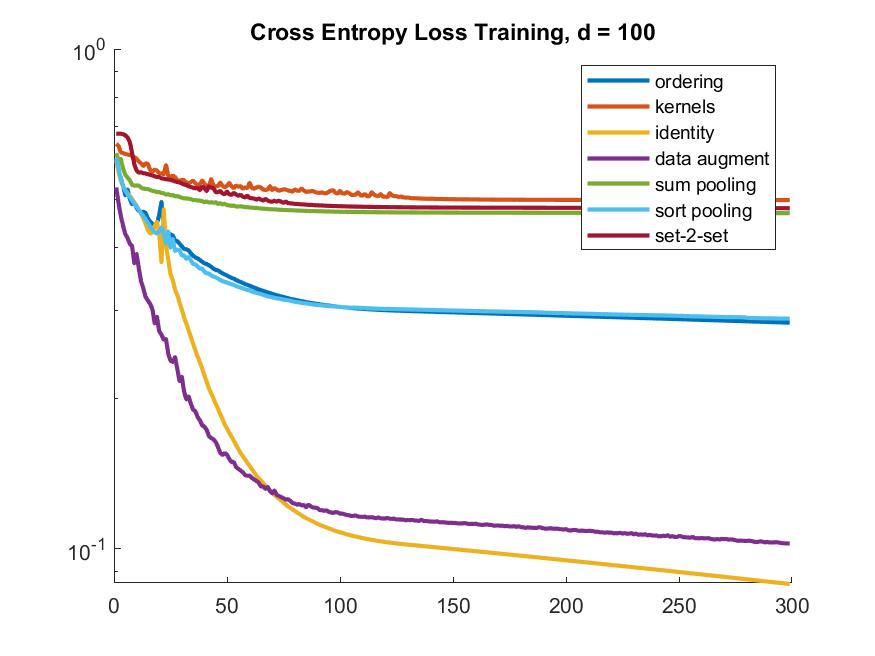}
	\\
	\includegraphics[width=0.49\linewidth]{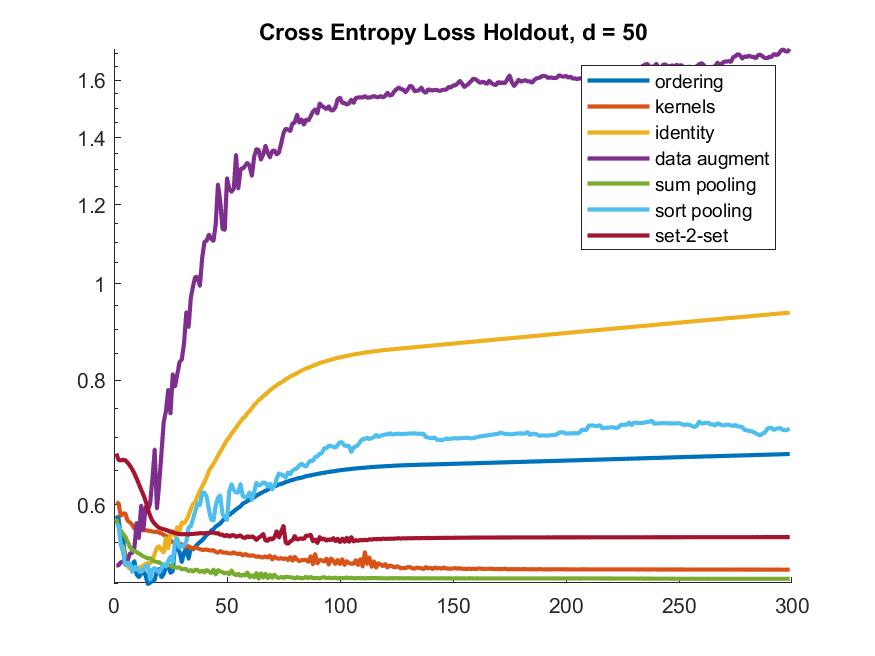}
	\includegraphics[width=0.49\linewidth]{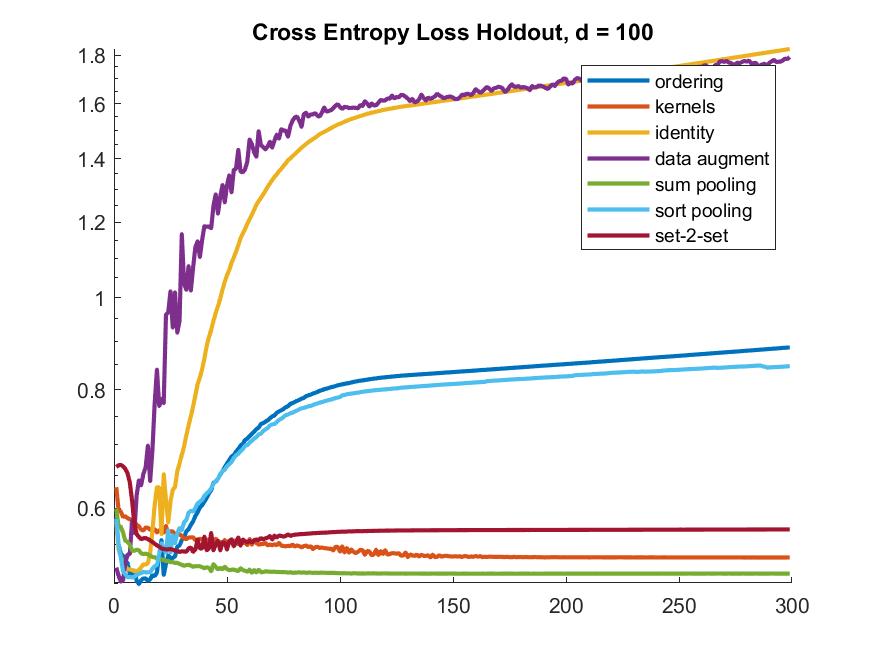}
    \\
	\includegraphics[width=0.49\linewidth]{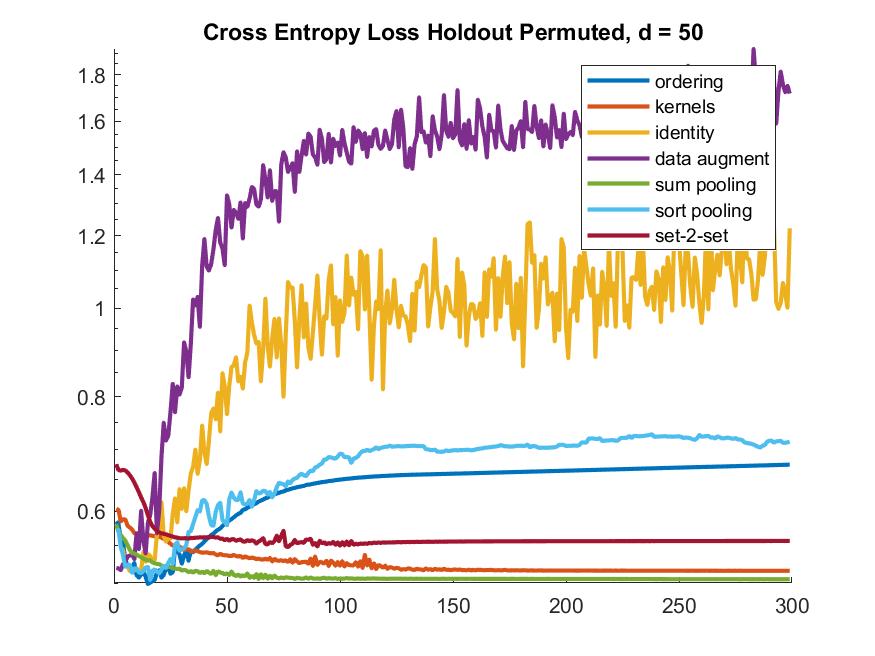}
	\includegraphics[width=0.49\linewidth]{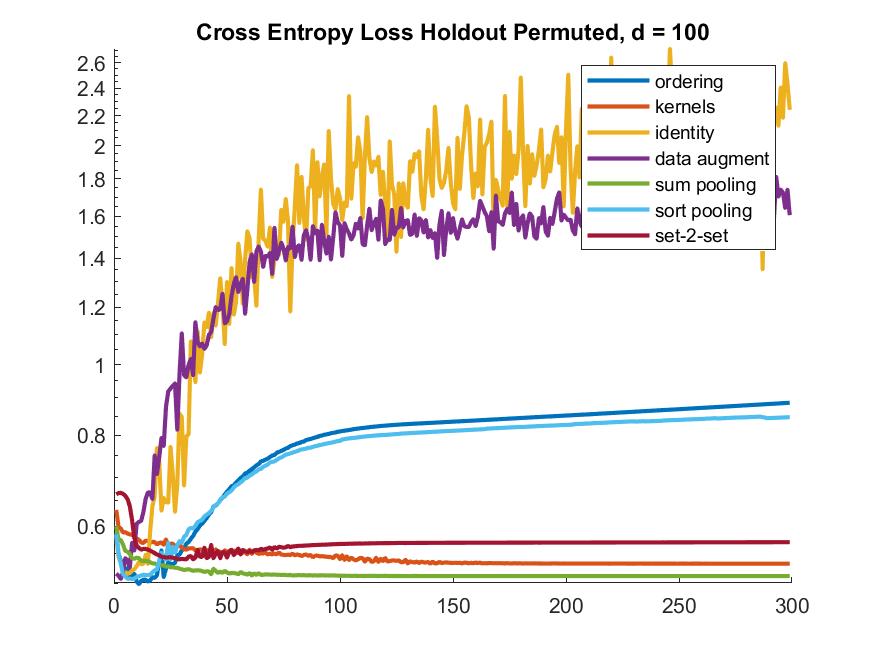}
	\caption[caption me]{Binary Cross Entropy (BCE) for enzyme/non-enzyme classification on PROTEINS\_FULL dataset using; $d=50$ left column; $d=100$ right column.}
	\label{fig:prot2}
\end{figure}

\begin{figure}[p]
	\includegraphics[width=0.49\linewidth]{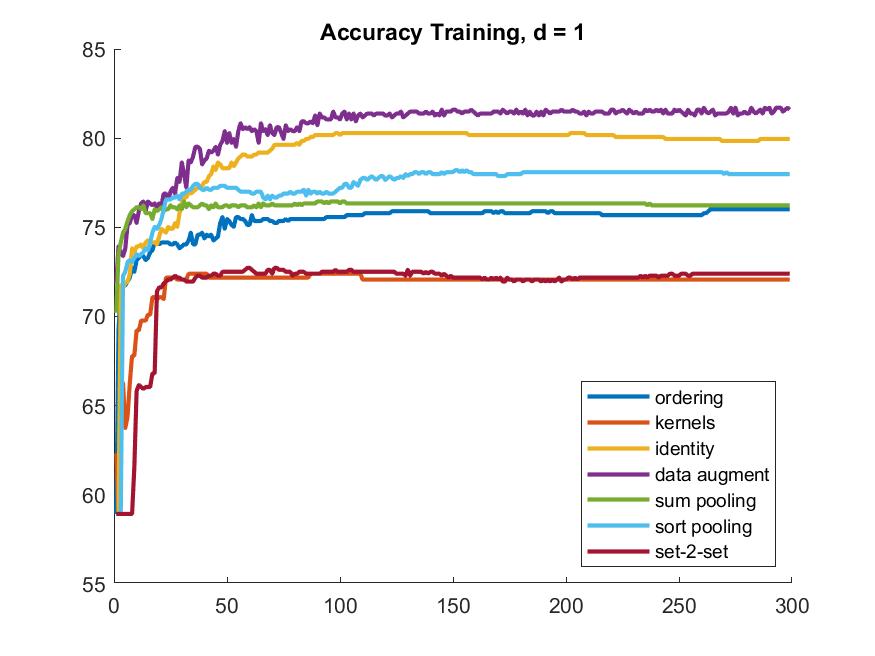}
	\includegraphics[width=0.49\linewidth]{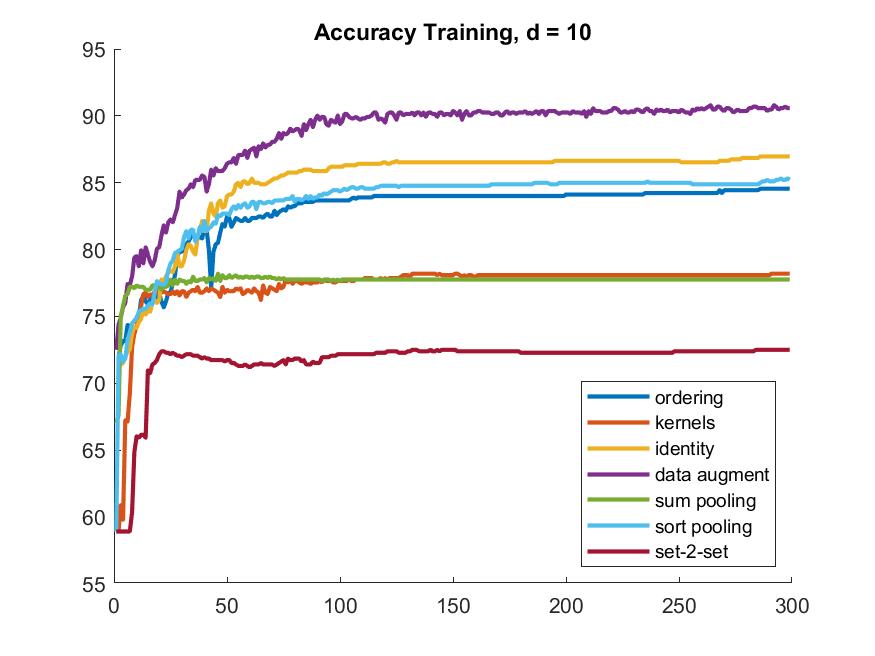}
	\\
	\includegraphics[width=0.49\linewidth]{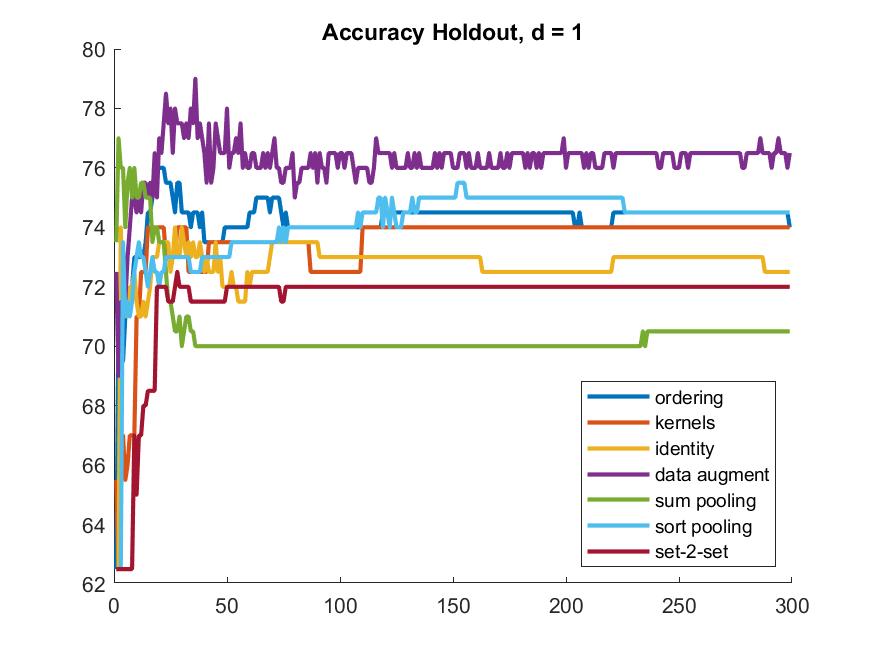}
	\includegraphics[width=0.49\linewidth]{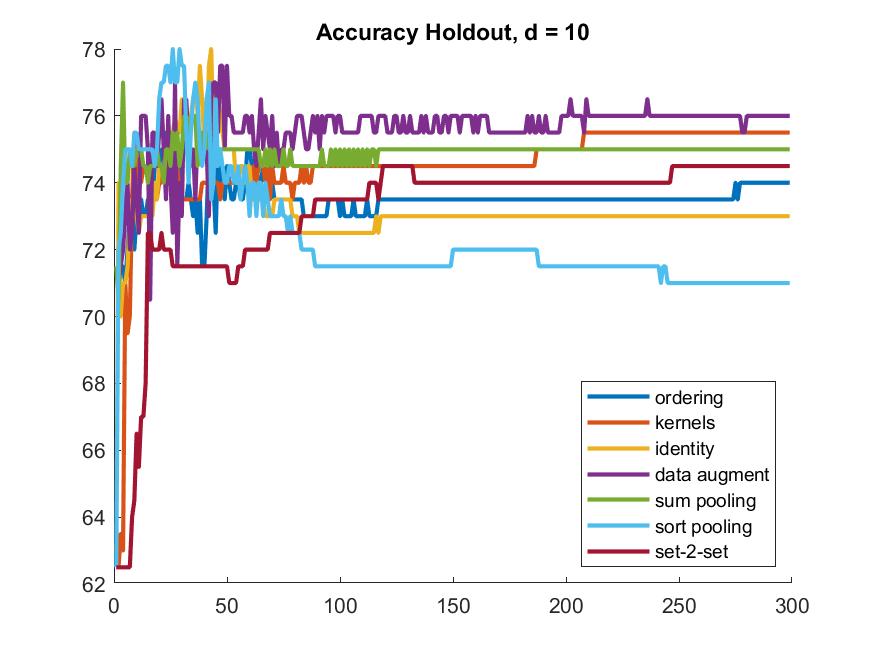}
    \\
	\includegraphics[width=0.49\linewidth]{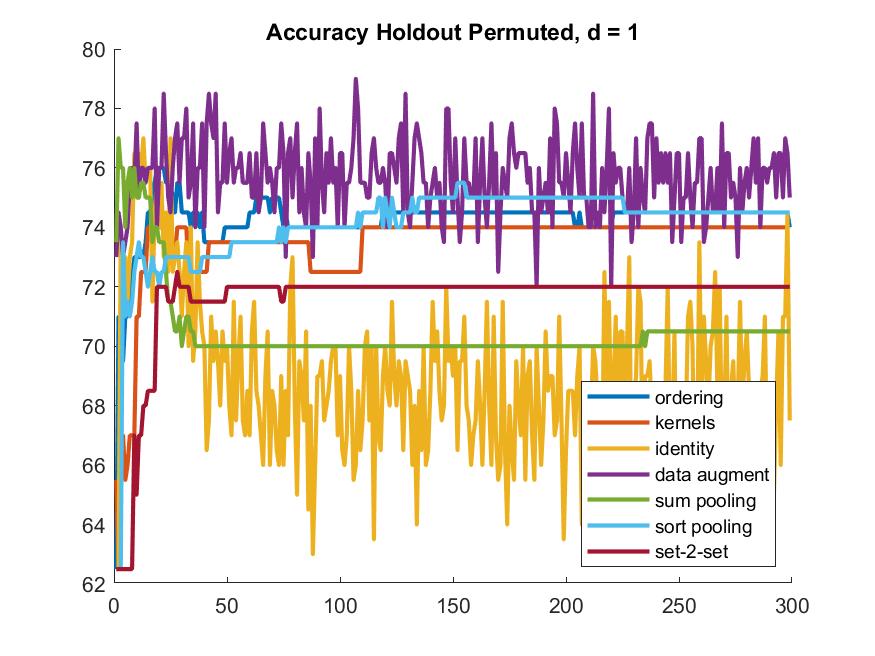}
	\includegraphics[width=0.49\linewidth]{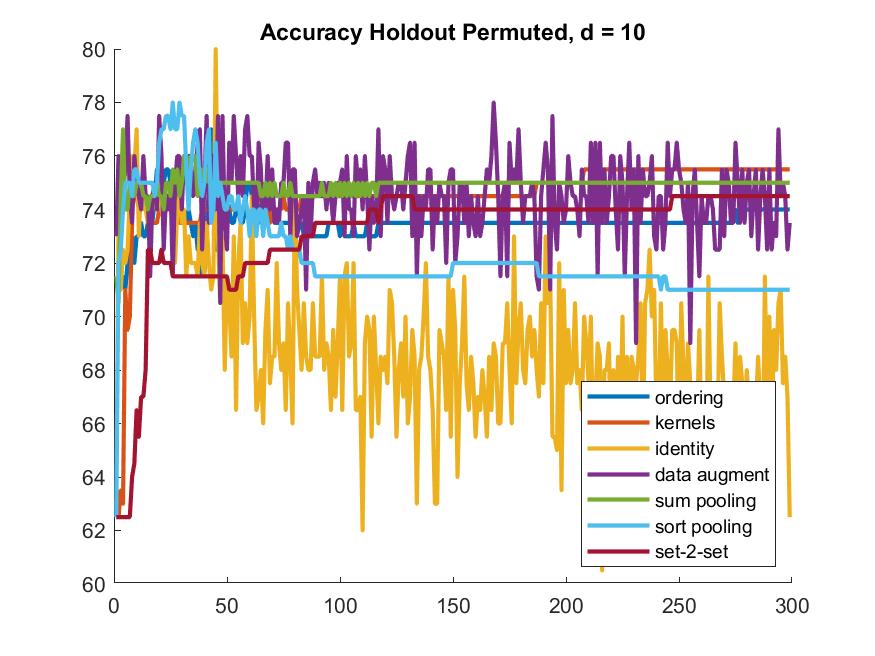}
	\caption[caption me]{Accuracy ACC(\%) for enzyme/non-enzyme classification on PROTEINS\_FULL dataset using; $d=1$ left column; $d=10$ right column.}
	\label{fig:prot3}
\end{figure}

\begin{figure}[p]
	\includegraphics[width=0.49\linewidth]{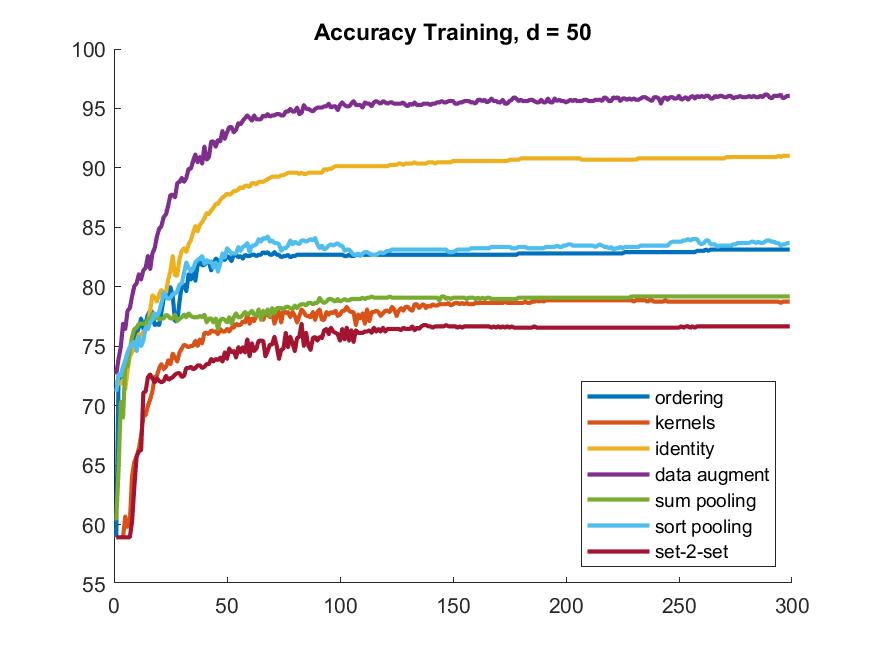}
	\includegraphics[width=0.49\linewidth]{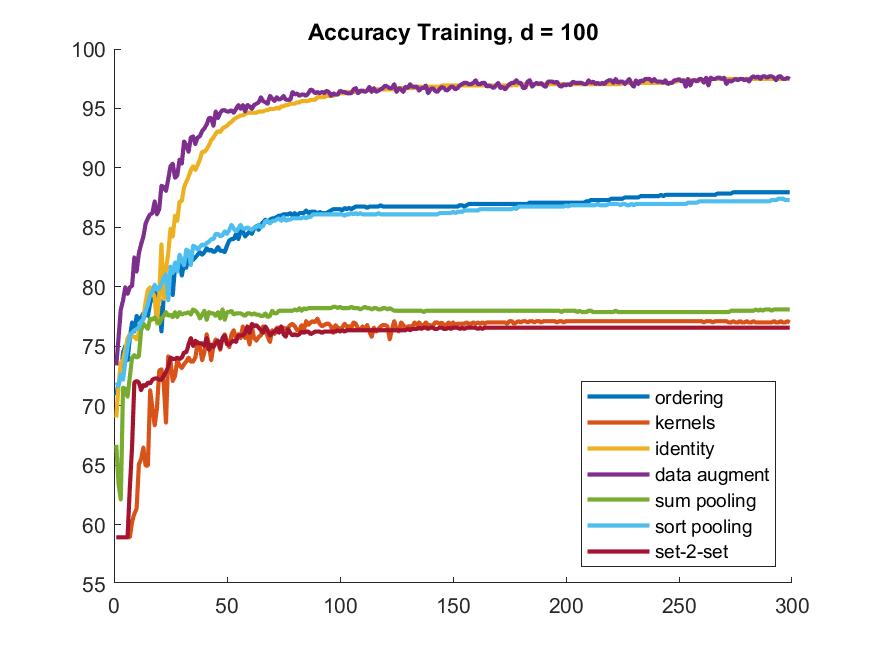}
	\\
	\includegraphics[width=0.49\linewidth]{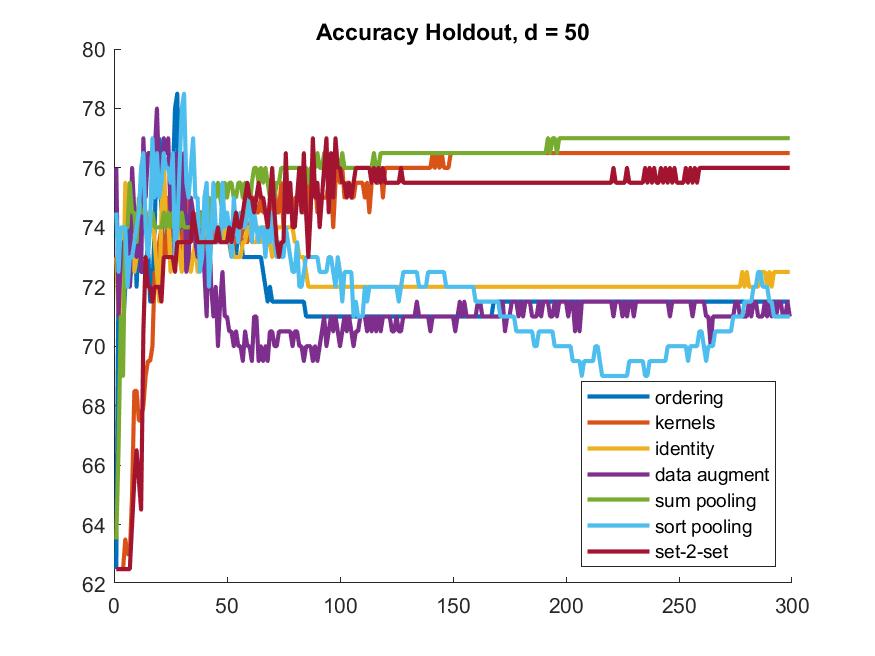}
	\includegraphics[width=0.49\linewidth]{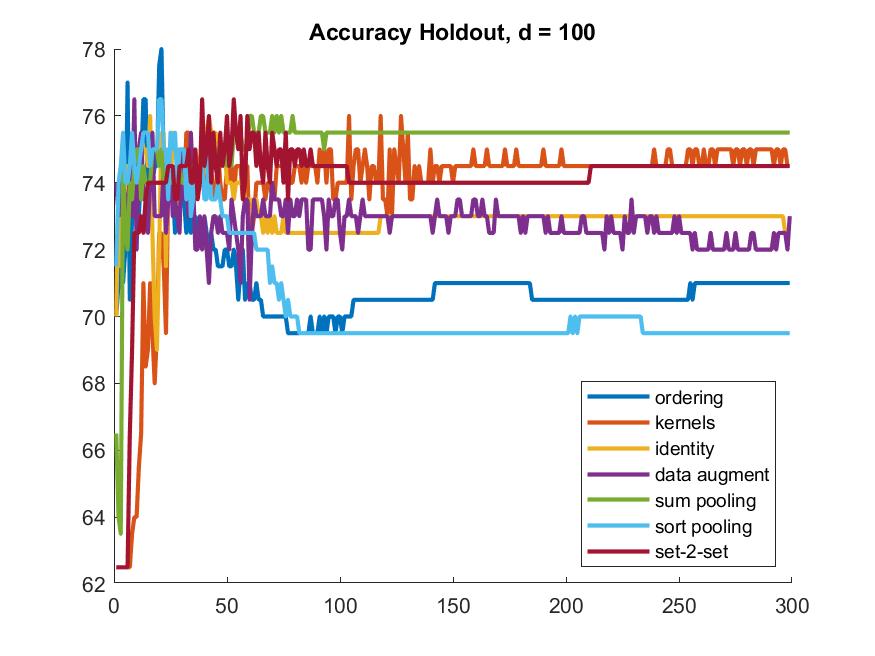}
    \\
	\includegraphics[width=0.49\linewidth]{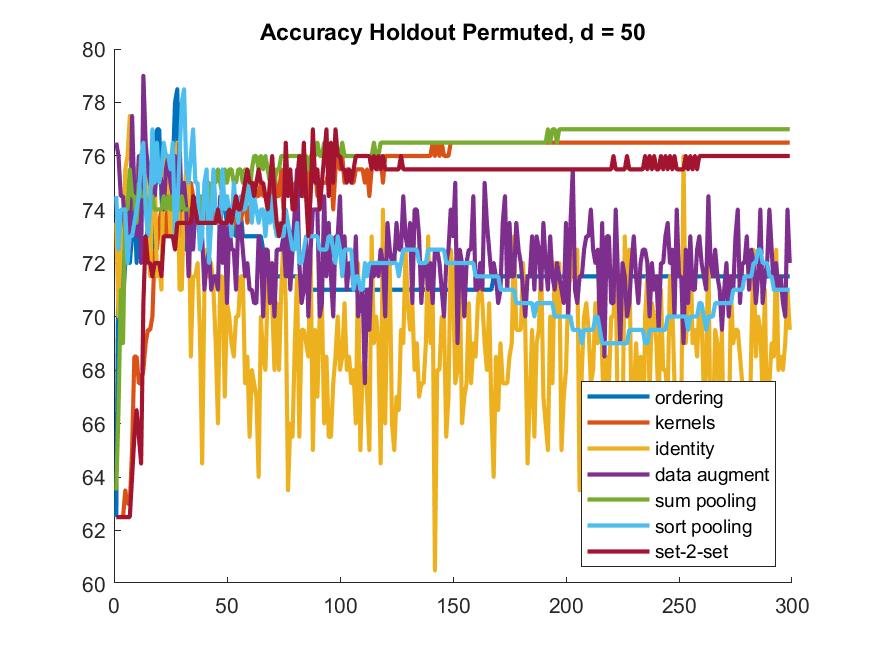}
	\includegraphics[width=0.49\linewidth]{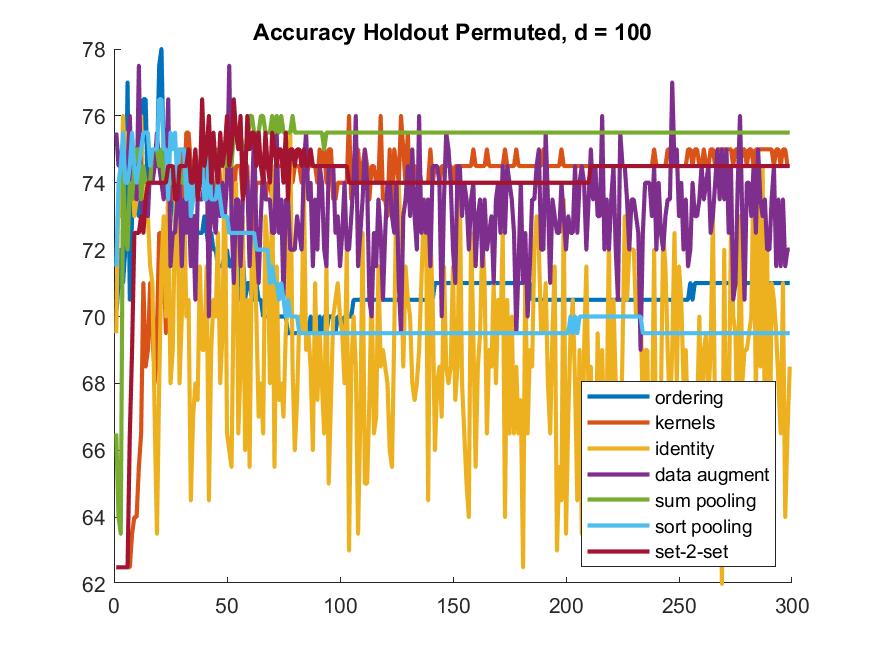}
	\caption[caption me]{Accuracy ACC(\%) for enzyme/non-enzyme classification on PROTEINS\_FULL dataset using; $d=50$ left column; $d=100$ right column.}
	\label{fig:prot4}
\end{figure}

\begin{figure}[p]
	\includegraphics[width=0.49\linewidth]{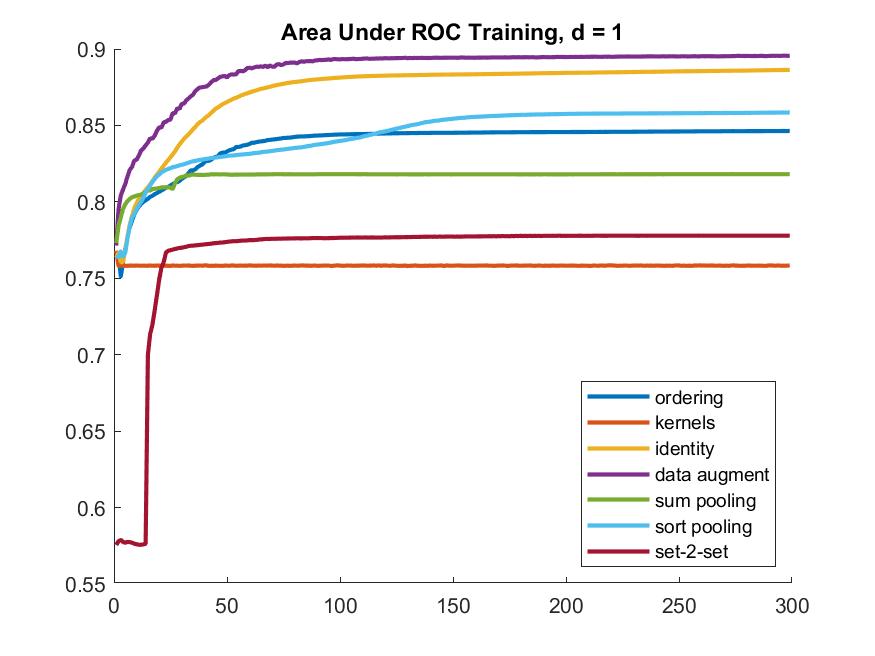}
	\includegraphics[width=0.49\linewidth]{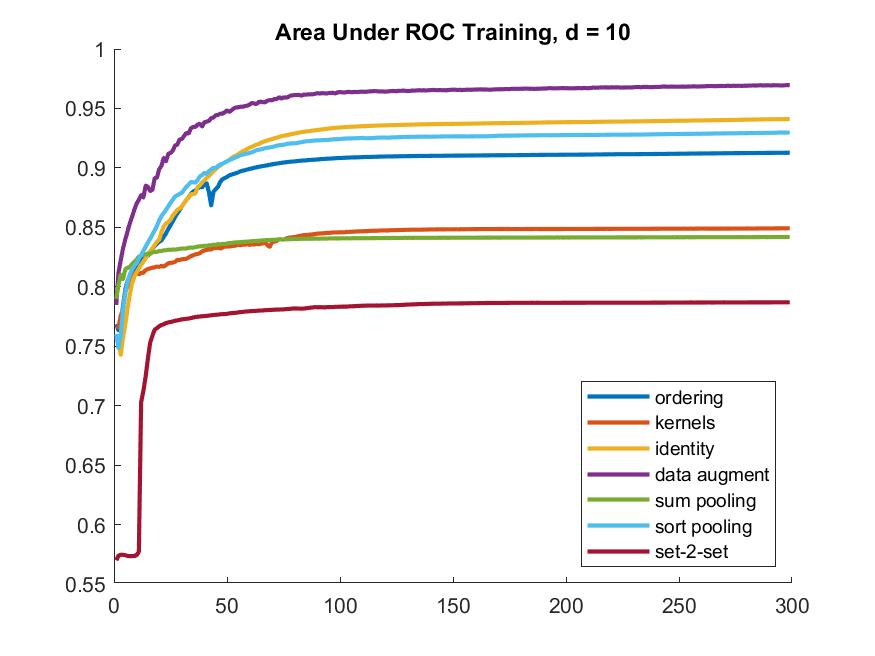}
	\\
	\includegraphics[width=0.49\linewidth]{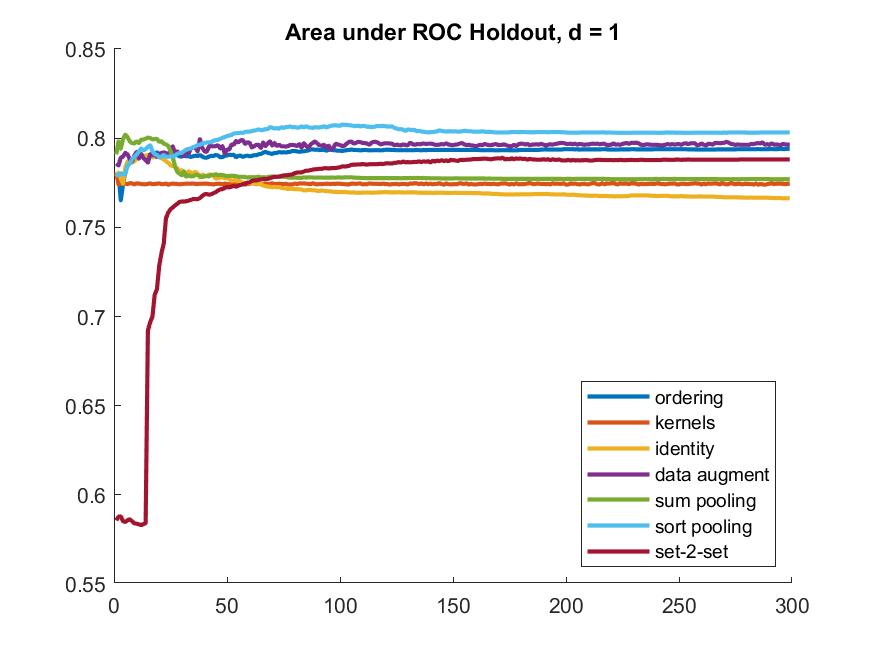}
	\includegraphics[width=0.49\linewidth]{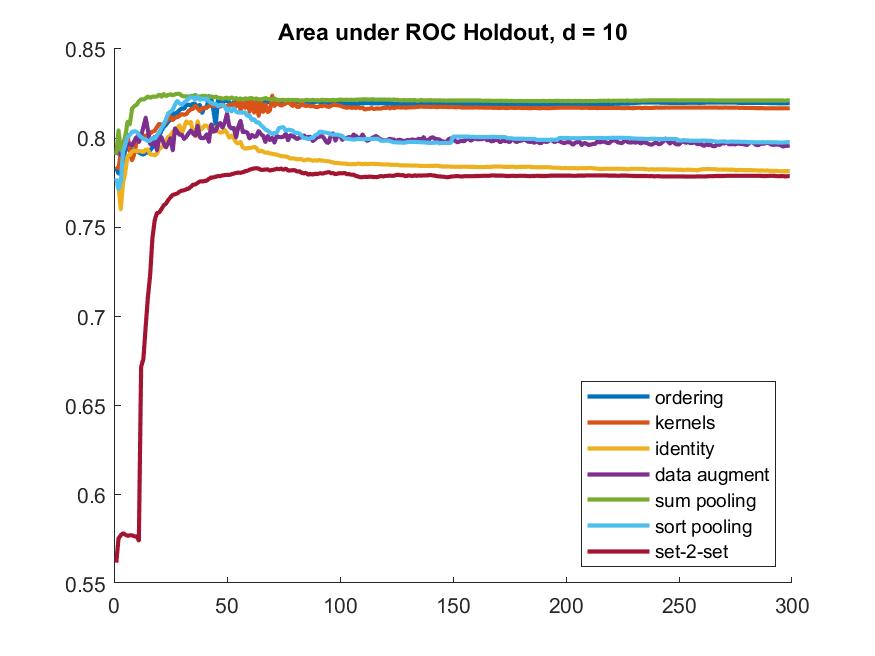}
    \\
	\includegraphics[width=0.49\linewidth]{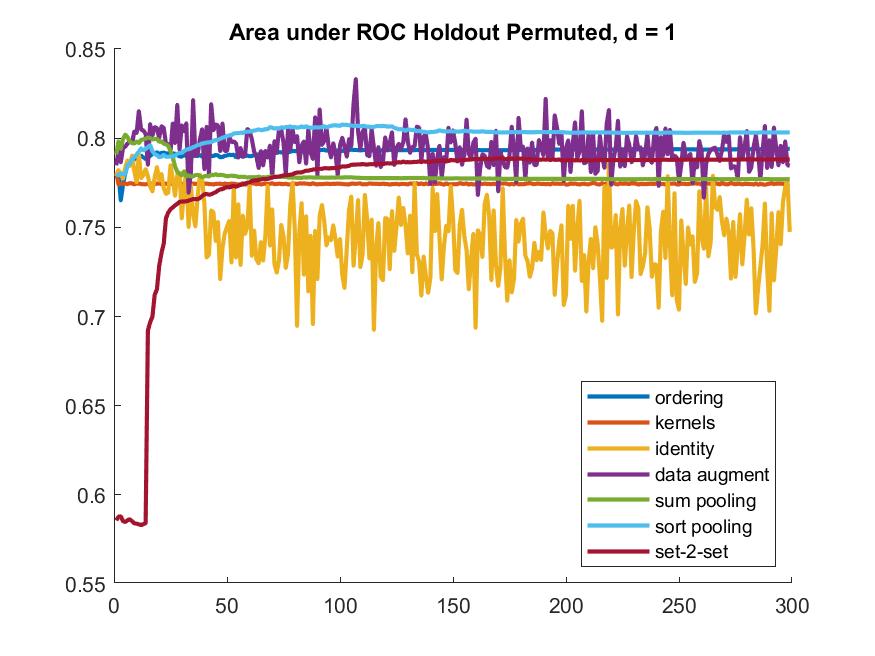}
	\includegraphics[width=0.49\linewidth]{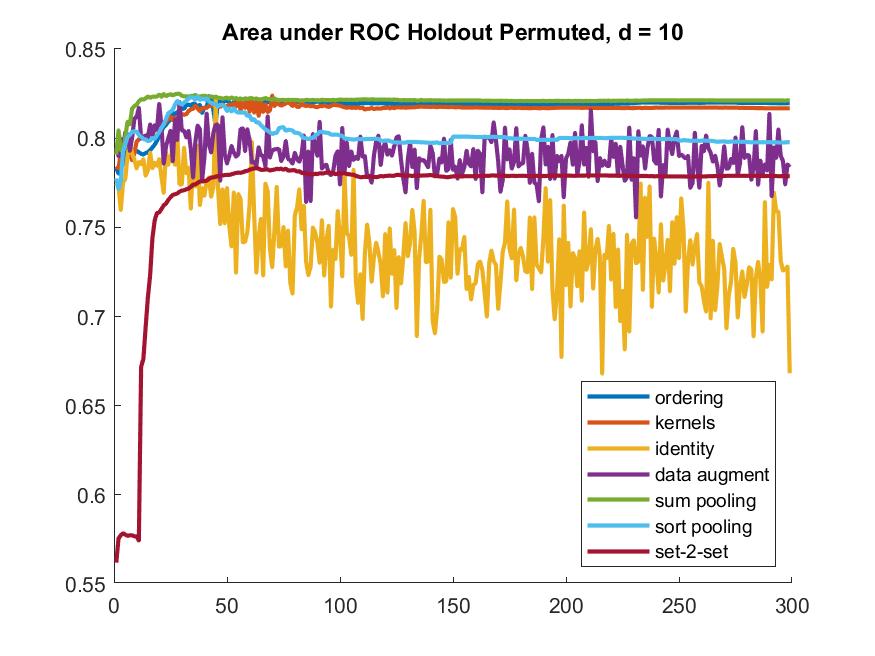}
	\caption[caption me]{Area under the ROC curve (AUC) for enzyme/non-enzyme classification on PROTEINS\_FULL dataset using; $d=1$ left column; $d=10$ right column.}
	\label{fig:prot5}
\end{figure}

\begin{figure}[p]
	\includegraphics[width=0.49\linewidth]{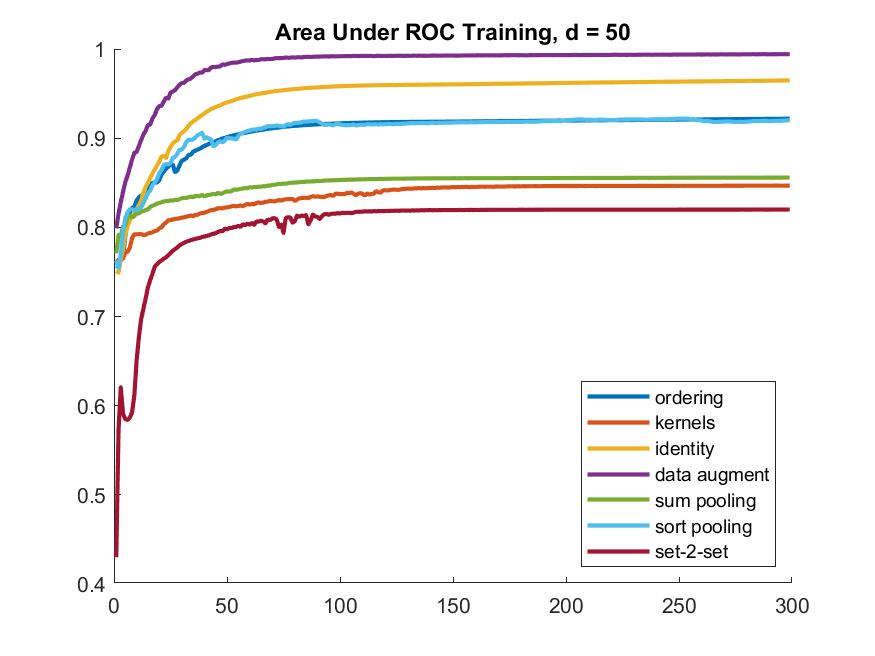}
	\includegraphics[width=0.49\linewidth]{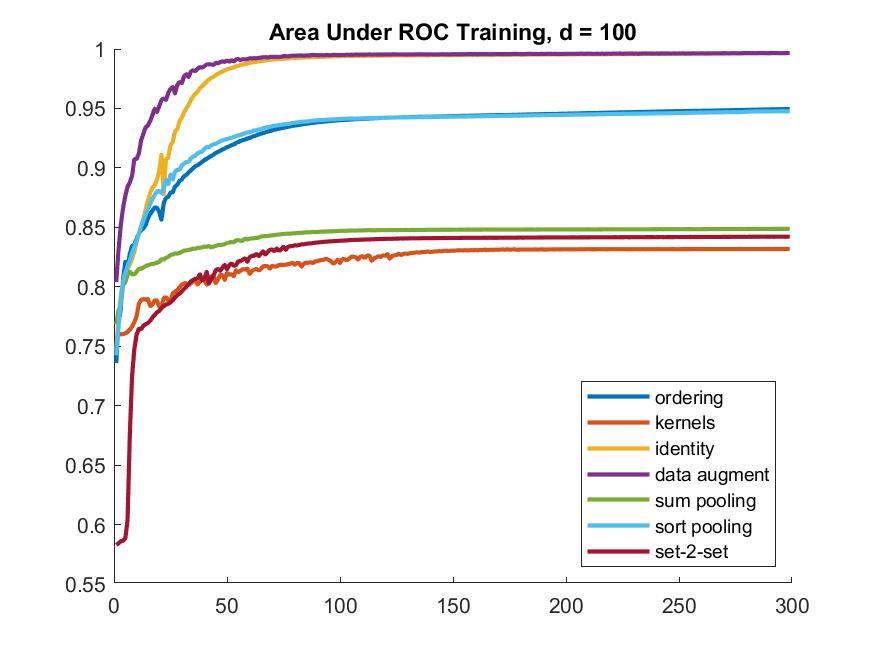}
	\\
	\includegraphics[width=0.49\linewidth]{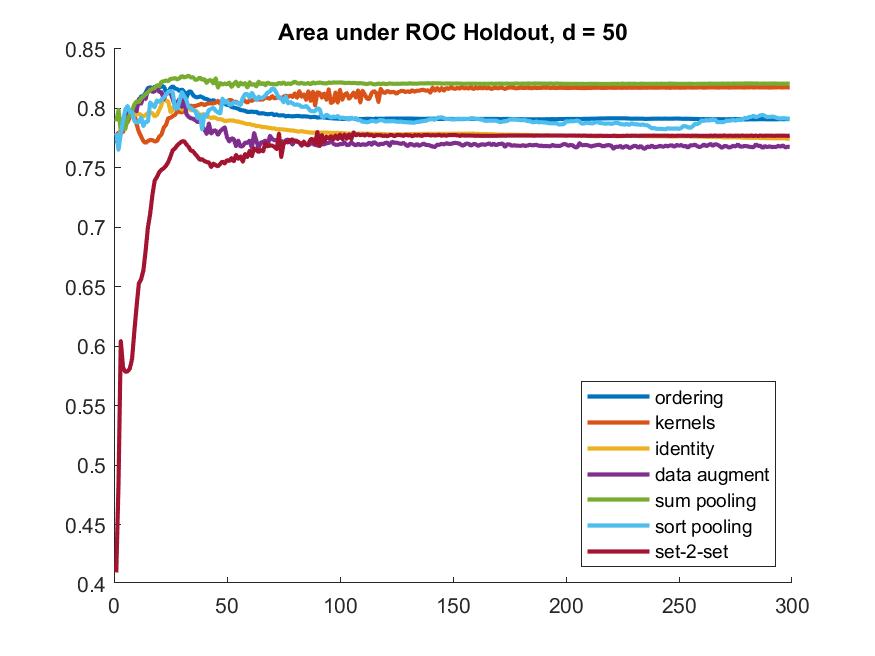}
	\includegraphics[width=0.49\linewidth]{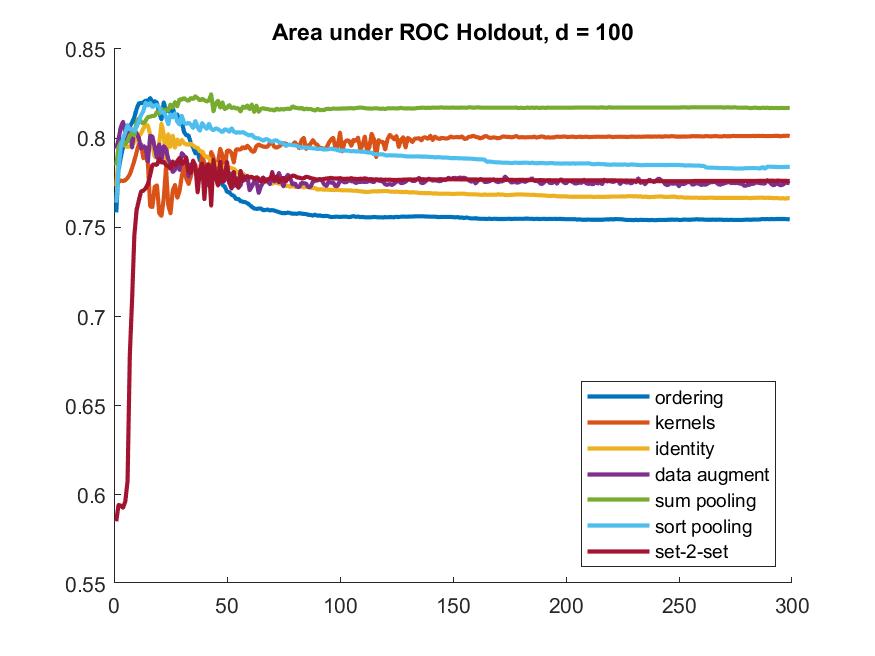}
    \\
	\includegraphics[width=0.49\linewidth]{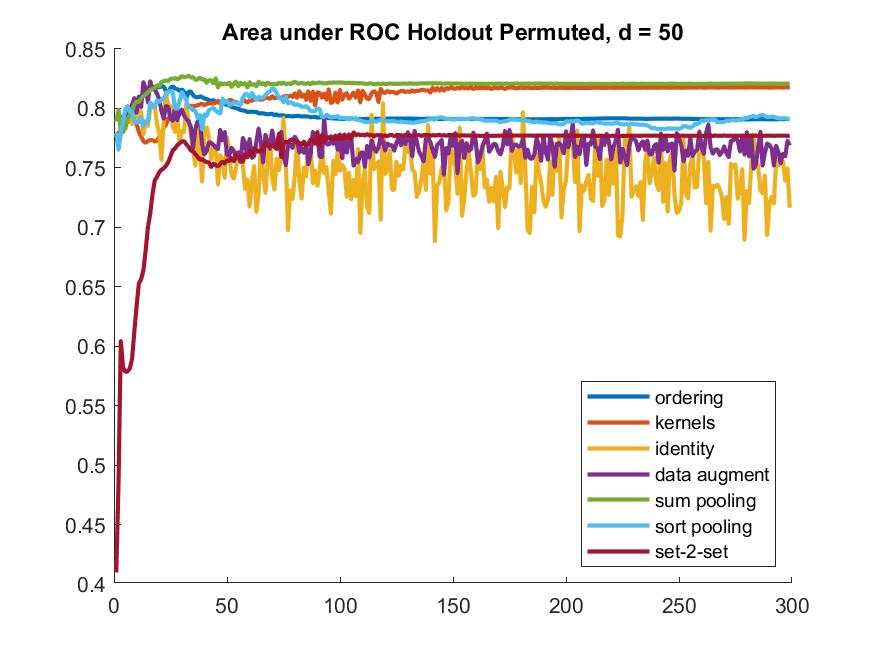}
	\includegraphics[width=0.49\linewidth]{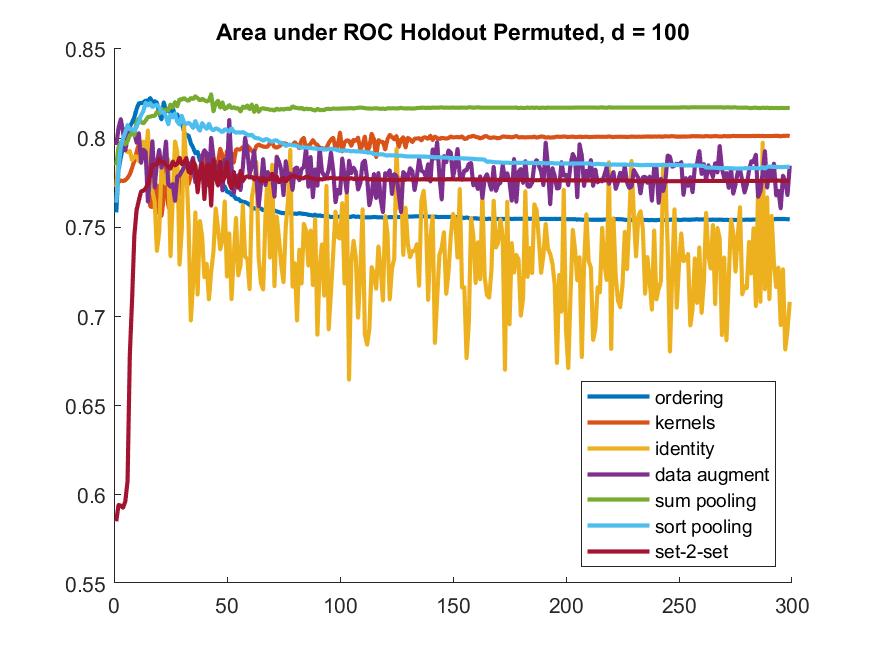}
	\caption[caption me]{Area under ROC curve (AUC) for enzyme/non-enzyme classification on PROTEINS\_FULL dataset using; $d=50$ left column; $d=100$ right column.}
	\label{fig:prot6}
\end{figure}

\begin{figure}[p]
	\includegraphics[width=0.49\linewidth]{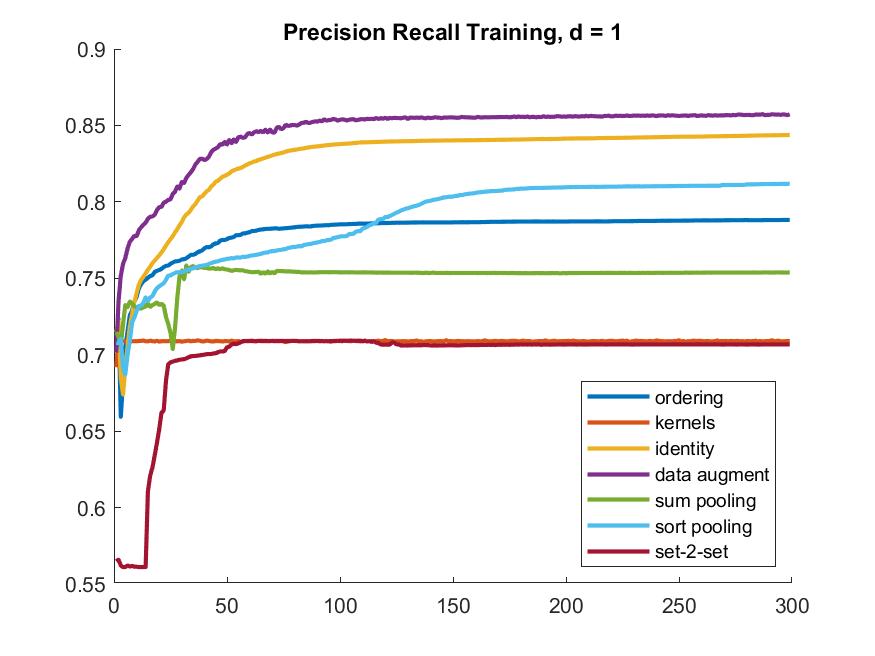}
	\includegraphics[width=0.49\linewidth]{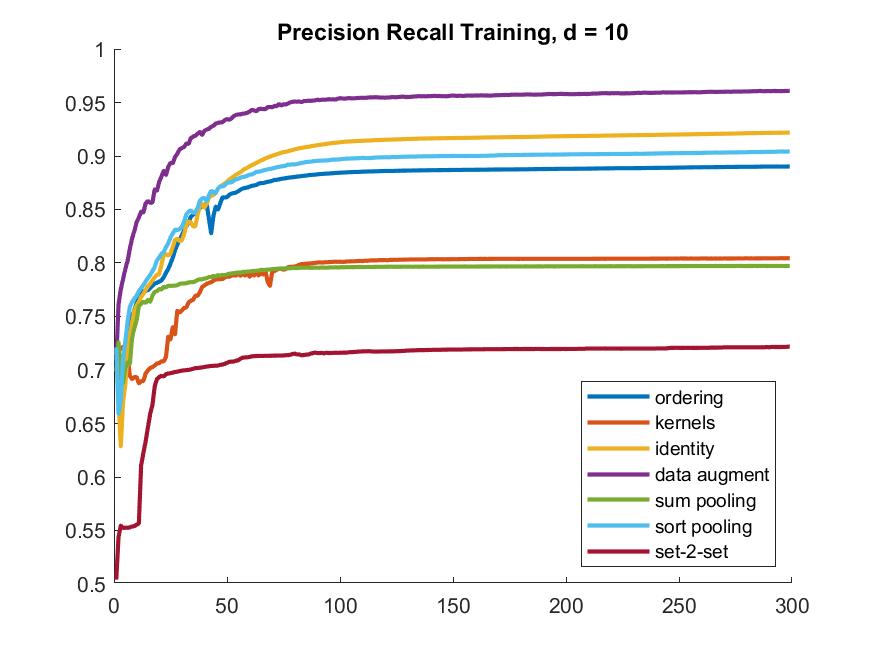}
	\\
	\includegraphics[width=0.49\linewidth]{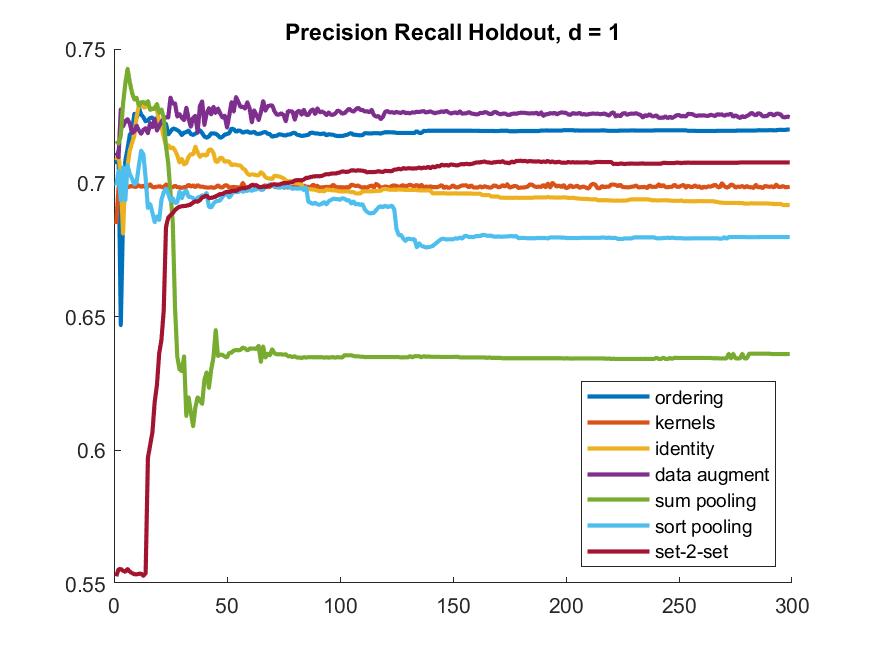}
	\includegraphics[width=0.49\linewidth]{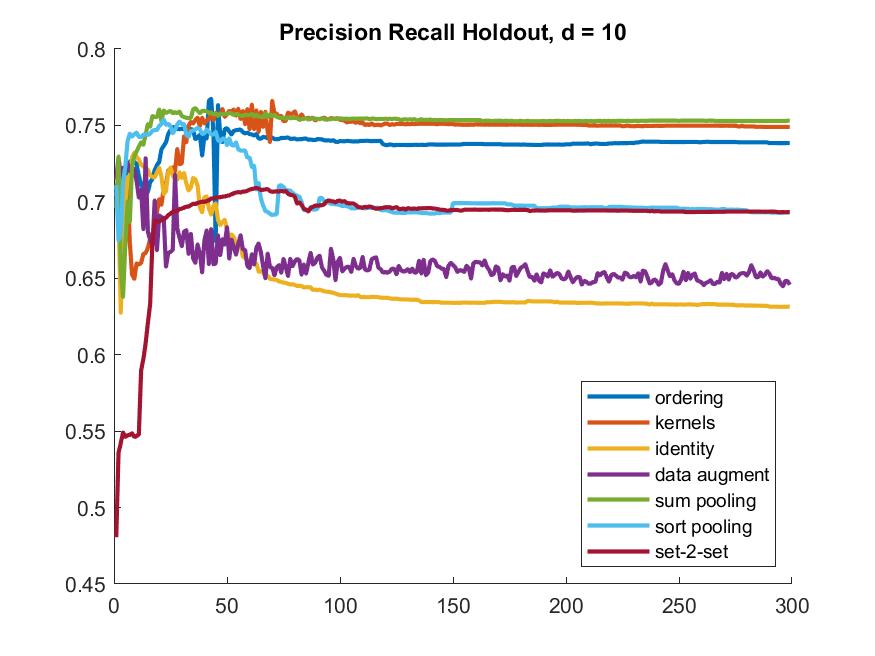}
    \\
	\includegraphics[width=0.49\linewidth]{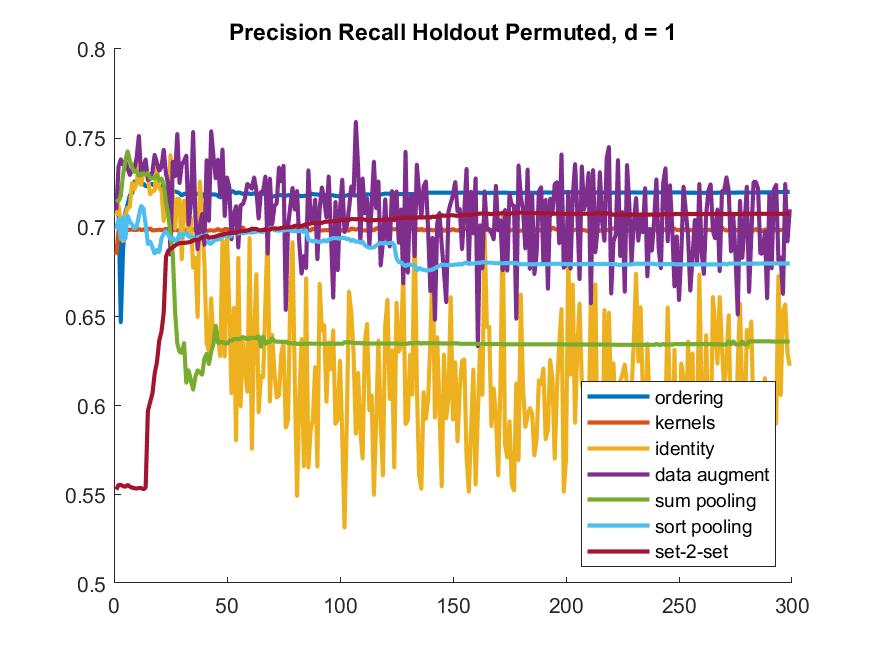}
	\includegraphics[width=0.49\linewidth]{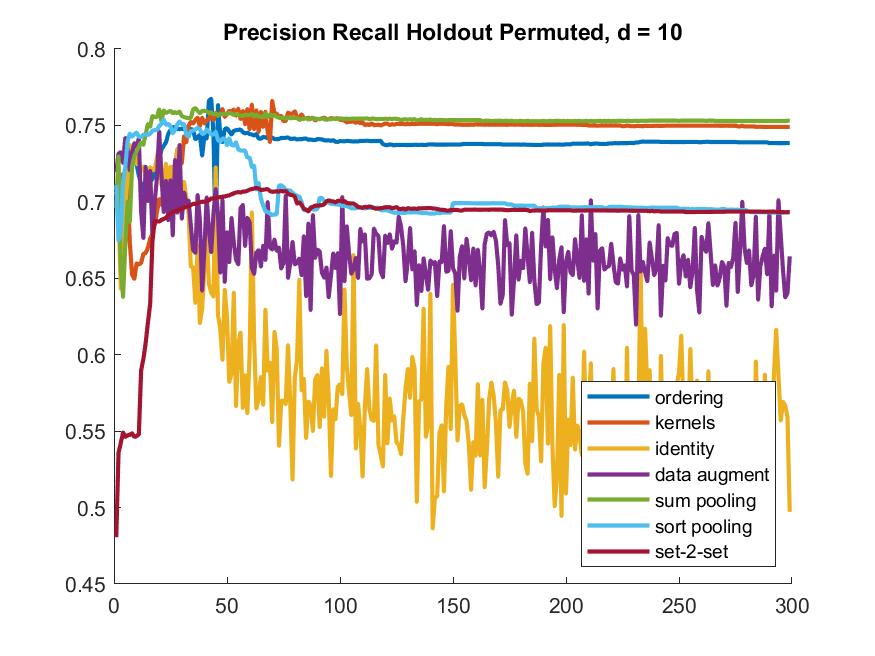}
	\caption[caption me]{Average precision  (AP) for enzyme/non-enzyme classification on PROTEINS\_FULL dataset using; $d=1$ left column; $d=10$ right column.}
	\label{fig:prot7}
\end{figure}

\begin{figure}[p]
	\includegraphics[width=0.49\linewidth]{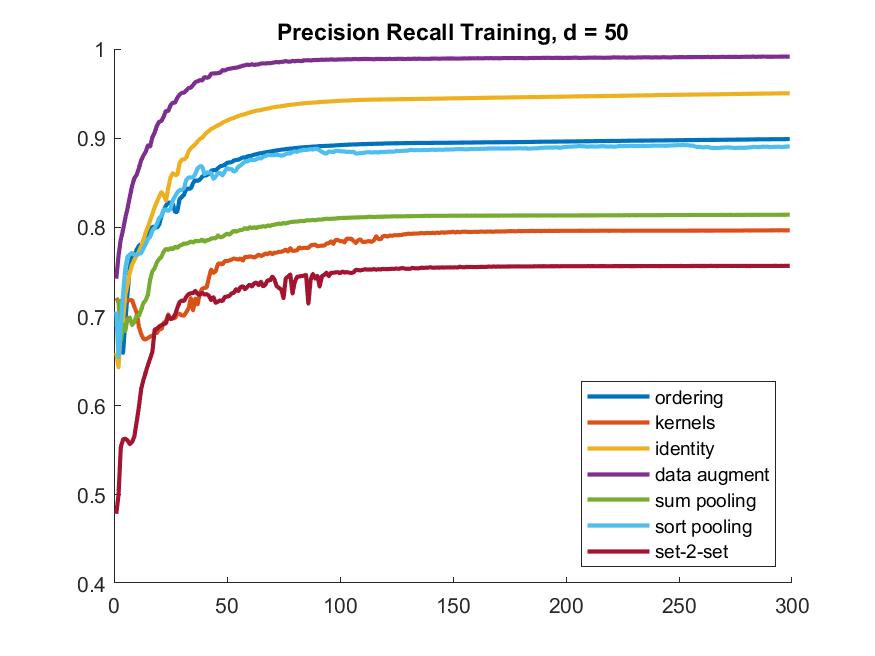}
	\includegraphics[width=0.49\linewidth]{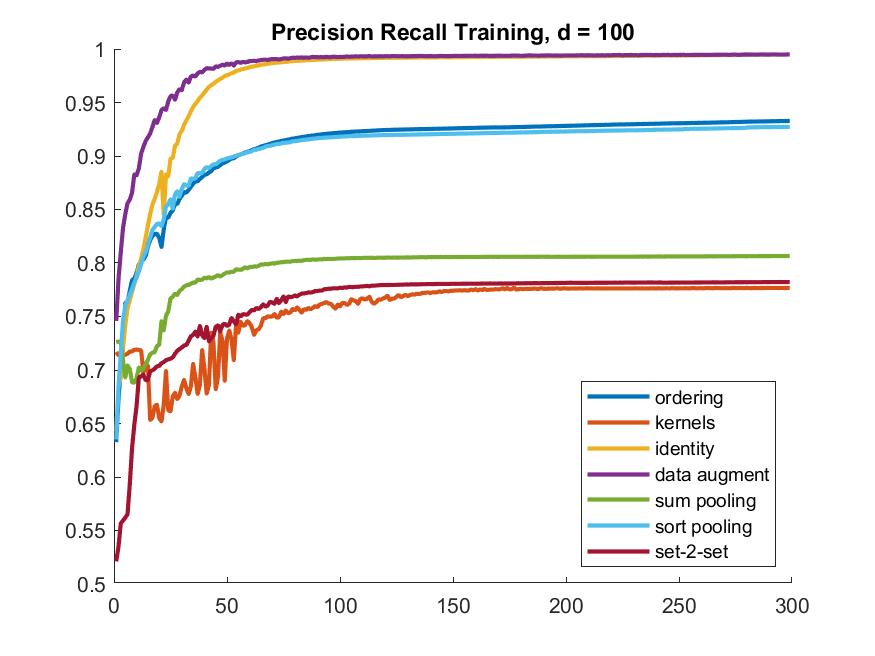}
	\\
	\includegraphics[width=0.49\linewidth]{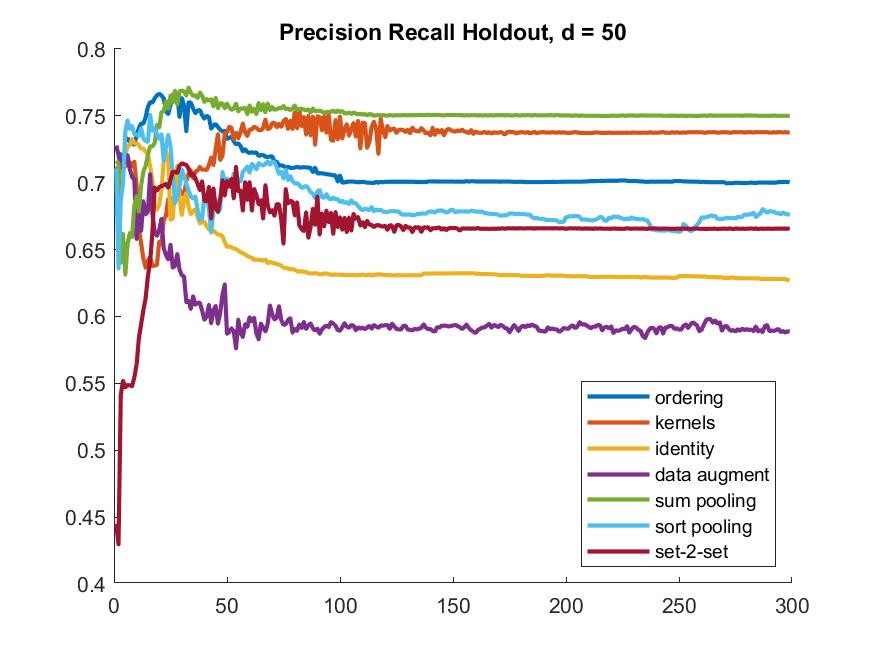}
	\includegraphics[width=0.49\linewidth]{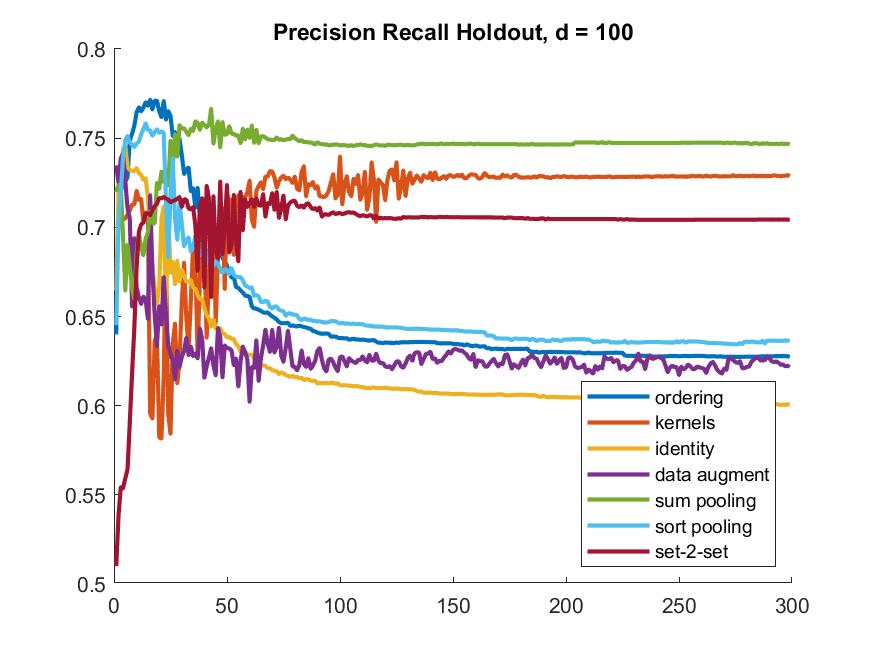}
    \\
	\includegraphics[width=0.49\linewidth]{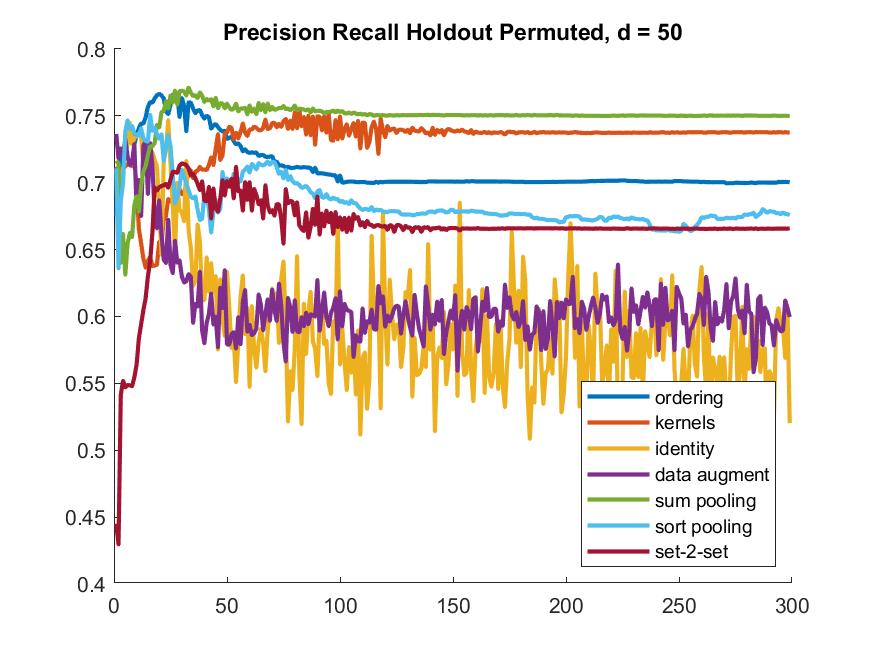}
	\includegraphics[width=0.49\linewidth]{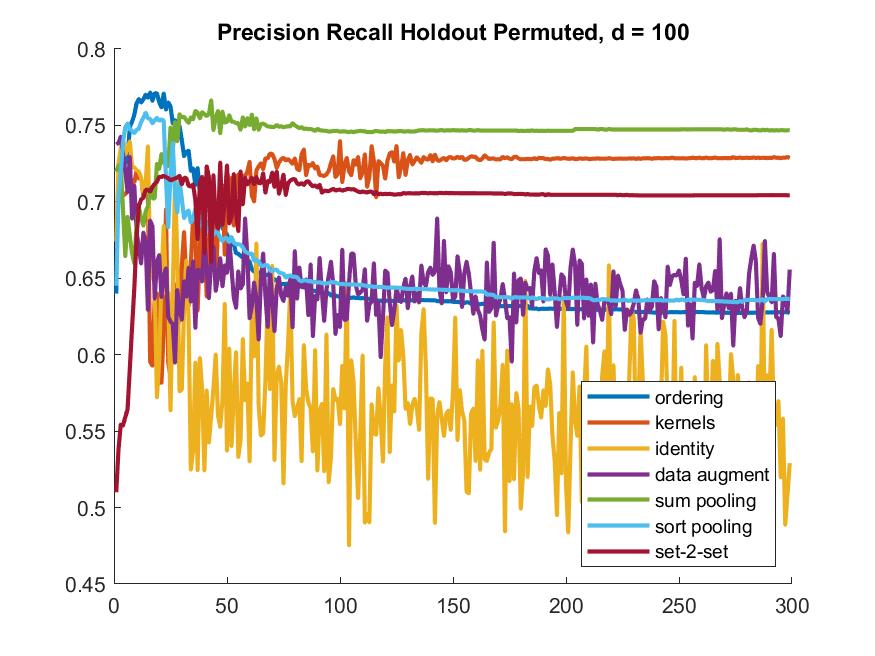}
	\caption[caption me]{Average precision (AP) for enzyme/non-enzyme classification on PROTEINS\_FULL dataset using; $d=50$ left column; $d=100$ right column.}
	\label{fig:prot8}
\end{figure}

\newpage
\section{Results for the QM9 dataset}

\begin{table}[pthb]
	\scalebox{1.0}{
	\begin{tabular}{|w{c}{2cm}|m{1.3cm} m{1.3cm} m{1.3cm} m{1.4cm} m{1.3cm} m{1.3cm} m{1.3cm}|}
		\hline
		d = 1 & ordering & kernels & identity & data augment & sum-pooling & sort-pooling & set-2-set  \\
		\hline
		Training & 0.302 & 0.867 & 0.320 & 0.281 & 0.349 & 0.309 & 0.389 \\
		\hline
		Holdout & 0.304 & 0.868 & 0.331 & 0.285 & 0.344 & 0.313 & 0.385 \\
		\hline
        Holdout Perm & 0.304 & 0.868 & 2.433 & 0.298 & 0.344 & 0.313 & 0.385 \\
		\hline
	\end{tabular}} 
	\caption{Mean Absolute Error (MAE) for regression of the electron energy gap $\Delta\varepsilon=LUMO-HOMO$ (eV) of the seven algorithms on QM9 dataset after 300 epochs for embedding dimension $d=1$}
	\label{table:a}
\end{table}

\begin{table}[!hbtp]
	\scalebox{1.0}{
	\begin{tabular}{|w{c}{2cm}|m{1.3cm} m{1.3cm} m{1.3cm} m{1.4cm} m{1.3cm} m{1.3cm} m{1.3cm}|}
		\hline
		d = 10 & ordering & kernels & identity & data augment & sum-pooling & sort-pooling & set-2-set  \\
		\hline
		Training & 0.220 & 0.219 & 0.182 & 0.175 & 0.214 & 0.226 & 0.282\\
		\hline
		Holdout & 0.232 & 0.222 & 0.244 & 0.208 & 0.223 & 0.278 & 0.287\\
		\hline
        Holdout Perm & 0.232 & 0.222 & 1.099 & 0.216 & 0.223 & 0.278 & 0.287\\
		\hline
	\end{tabular}} 
	\caption{Mean Absolute Error (MAE) for regression of the electron energy gap $\Delta\varepsilon=LUMO-HOMO$ (eV) of the seven algorithms on QM9 dataset after 300 epochs for embedding dimension $d=10$}
	\label{table:b}
\end{table}

\begin{table}[!hbtp]
	\scalebox{1.0}{
	\begin{tabular}{|w{c}{2cm}|m{1.3cm} m{1.3cm} m{1.3cm} m{1.4cm} m{1.3cm} m{1.3cm} m{1.3cm}|}
		\hline
		d = 50 & ordering & kernels & identity & data augment & sum-pooling & sort-pooling & set-2-set  \\
		\hline
		Training & 0.163 & 0.257 & 0.163 & 0.172 & 0.182 & 0.166 & 0.196\\
		\hline
		Holdout & 0.191 & 0.258 & 0.234 & 0.212 & 0.204 & 0.227 & 0.211\\
		\hline
        Holdout Perm & 0.191 & 0.258 & 1.607 & 0.219 & 0.204 & 0.277 & 0.211\\
		\hline
	\end{tabular}} 
	\caption{Mean Absolute Error (MAE) for regression of the electron energy gap $\Delta\varepsilon=LUMO-HOMO$ (eV) of the seven algorithms on QM9 dataset after 300 epochs for embedding dimension $d=50$}
	\label{table:c}
\end{table}

\begin{table}[!hbtp]
	\scalebox{1.0}{
	\begin{tabular}{|w{c}{2cm}|m{1.3cm} m{1.3cm} m{1.3cm} m{1.4cm} m{1.3cm} m{1.3cm} m{1.3cm}|}
		\hline
		d = 100 & ordering & kernels & identity & data augment & sum-pooling & sort-pooling & set-2-set  \\
		\hline
		Training & 0.155 & 0.269 & 0.139 & 0.164 & 0.178 & 0.199 & 0.173\\
		\hline
		Holdout & 0.187 & 0.267 & 0.227 & 0.206 & 0.201 & 0.239 & 0.201\\
		\hline
        Holdout Perm & 0.187 & 0.267 & 1.086 & 0.213 & 0.201 & 0.239 & 0.201\\
		\hline
	\end{tabular}} 
	\caption{Mean Absolute Error (MAE) for regression of the electron energy gap $\Delta\varepsilon=LUMO-HOMO$ (eV) of the seven algorithms on QM9 dataset after 300 epochs for embedding dimension $d=100$}
	\label{table:d}
\end{table}



\begin{figure}[hp]
	\includegraphics[width=0.49\linewidth]{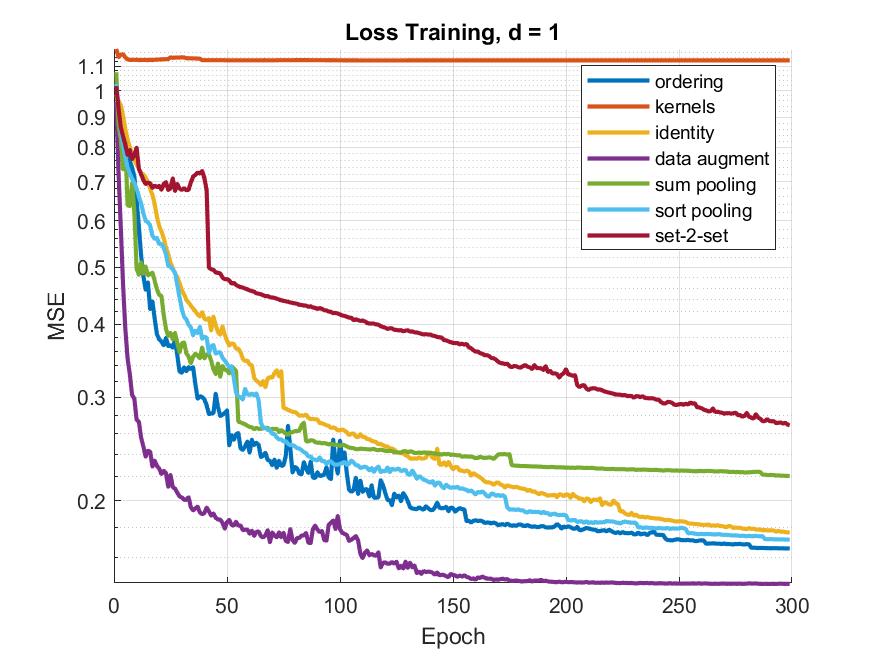}
	\includegraphics[width=0.49\linewidth]{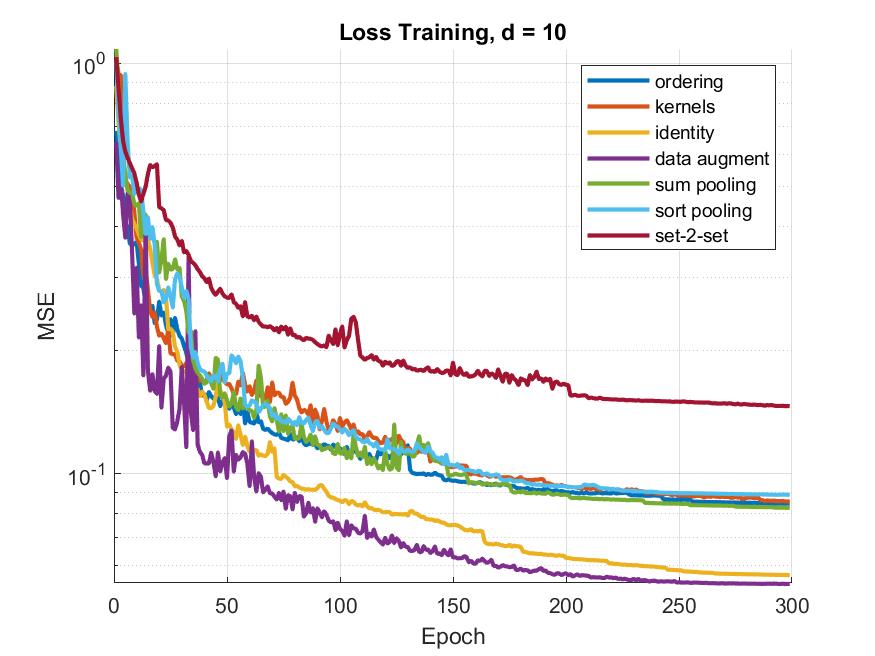}
 \\
	\includegraphics[width=0.49\linewidth]{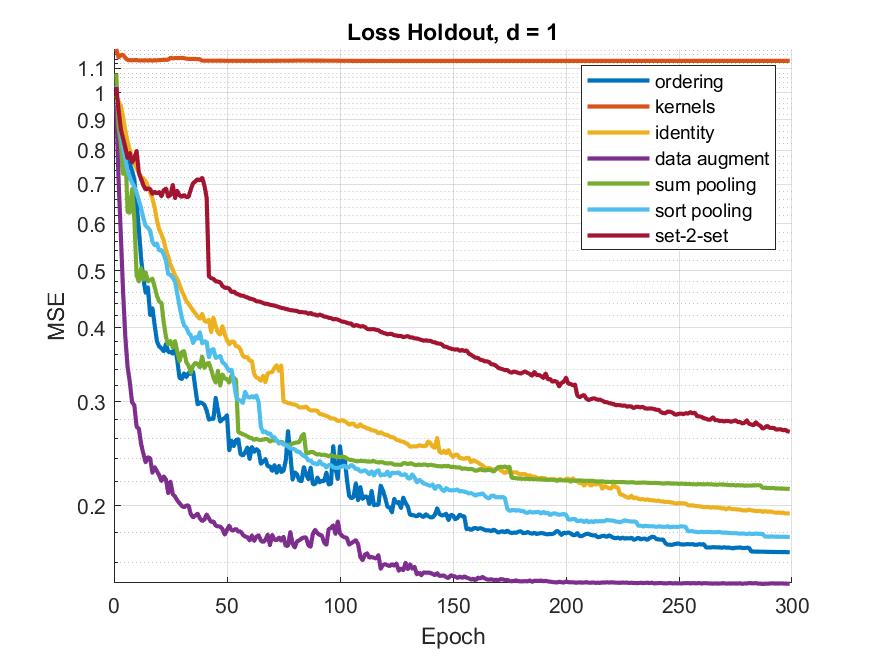}
	\includegraphics[width=0.49\linewidth]{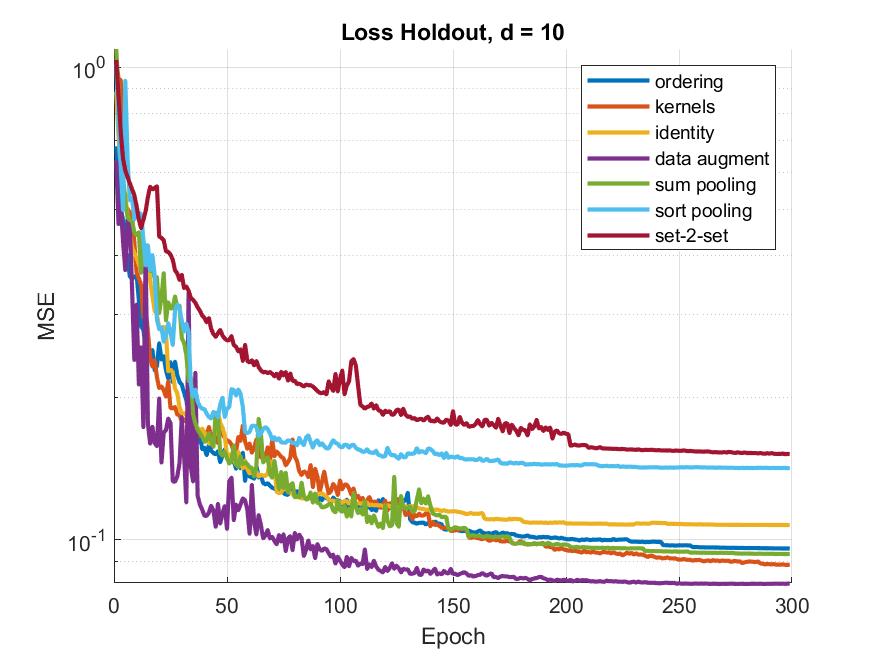}
 \\
	\includegraphics[width=0.49\linewidth]{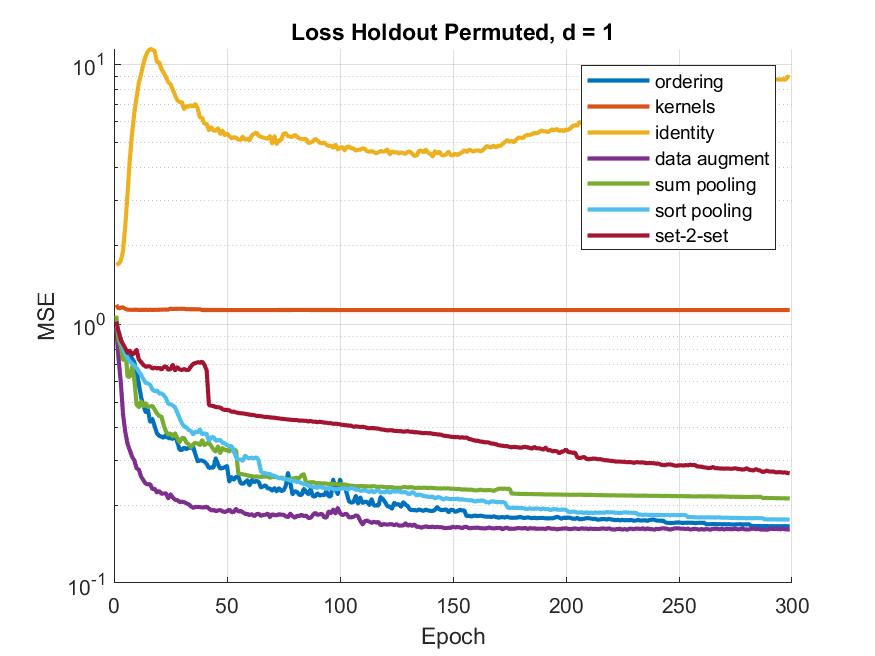}
	\includegraphics[width=0.49\linewidth]{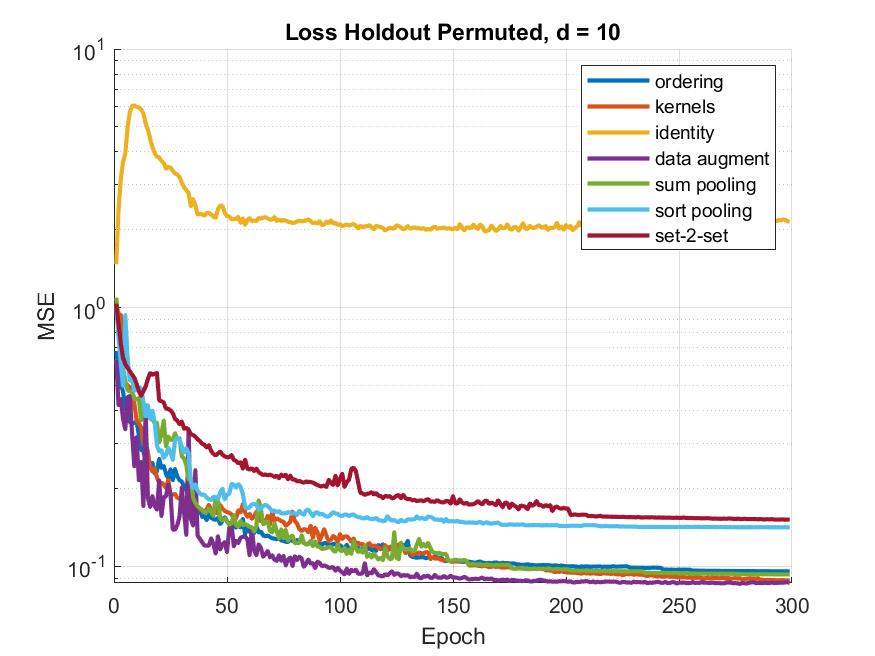}
	\caption[caption me]{Mean-Square Error (MSE) Results for regression of the electron gap energy $\Delta\varepsilon=LUMO-HOMO$ (eV) on QM9 dataset using the seven algorithms; $d=1$ left column; $d=10$ right column.}
	\label{fig:a1}
\end{figure}

\begin{figure}[hp]
\includegraphics[width=0.49\linewidth]{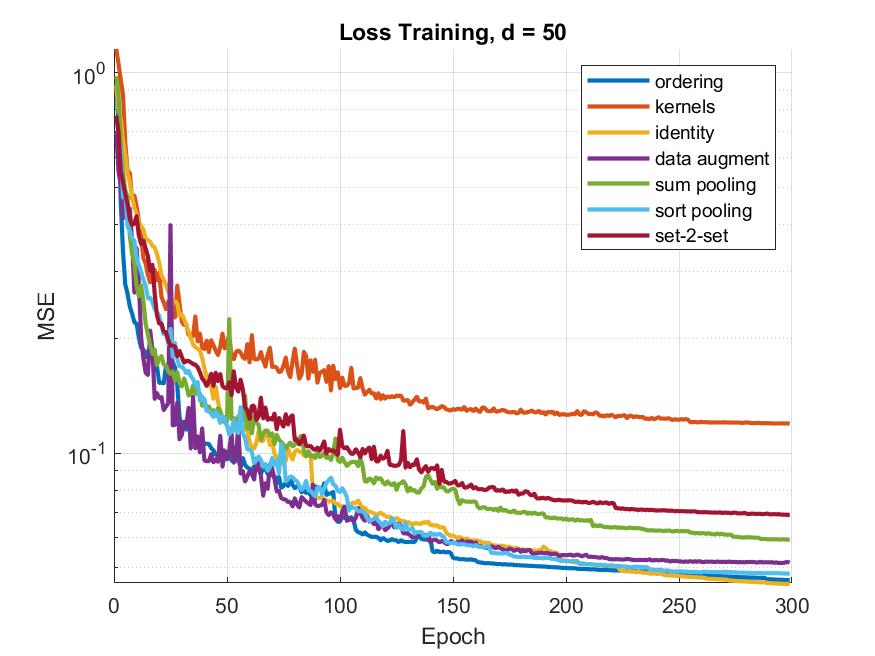}
\includegraphics[width=0.49\linewidth]{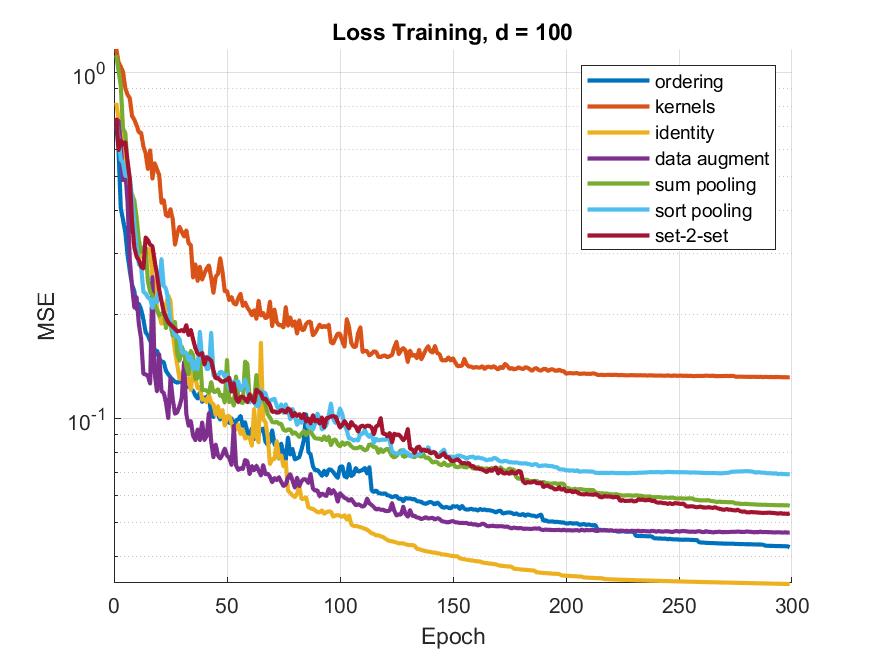} \\
\includegraphics[width=0.49\linewidth]{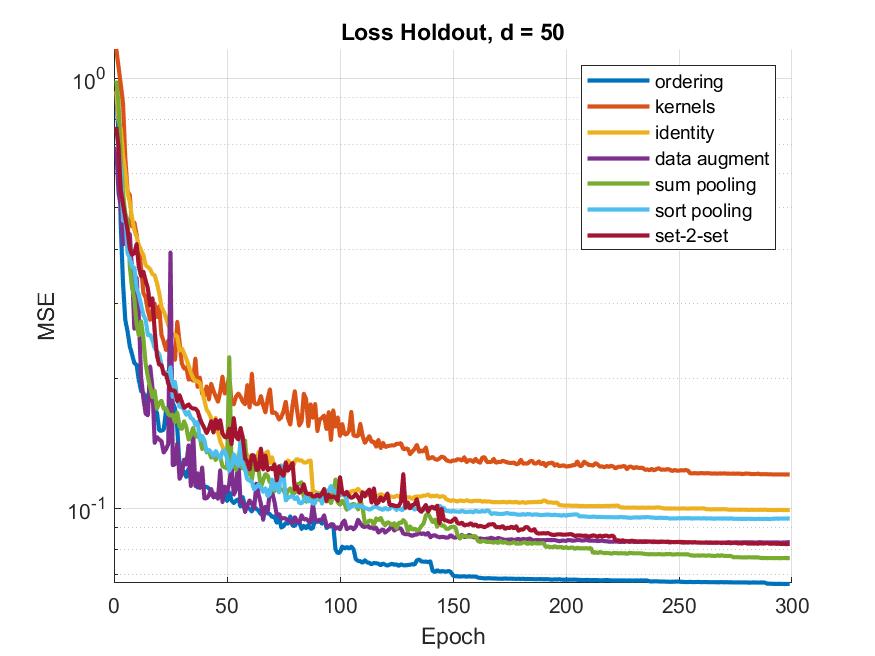} 
\includegraphics[width=0.49\linewidth]{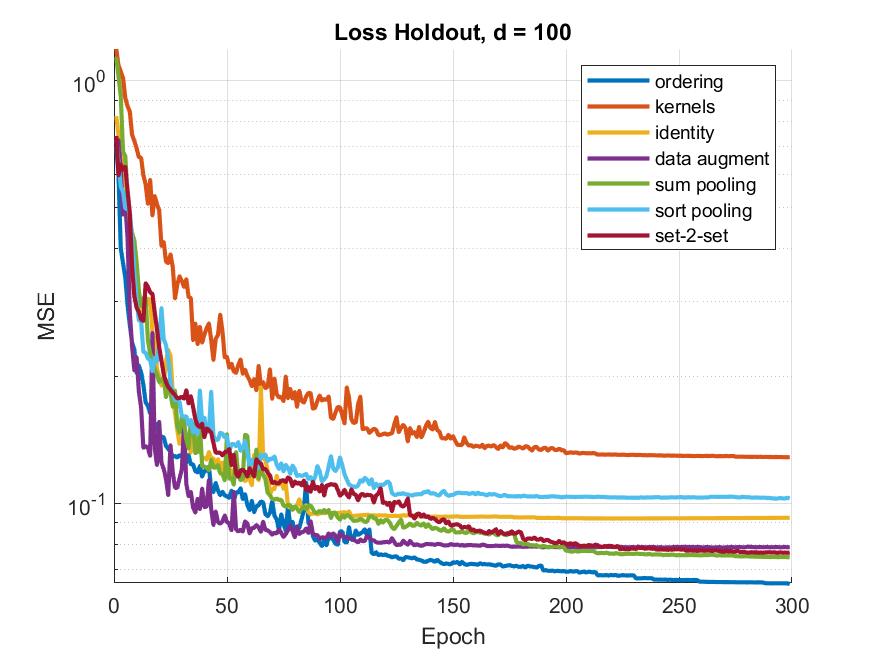} \\
\includegraphics[width=0.49\linewidth]{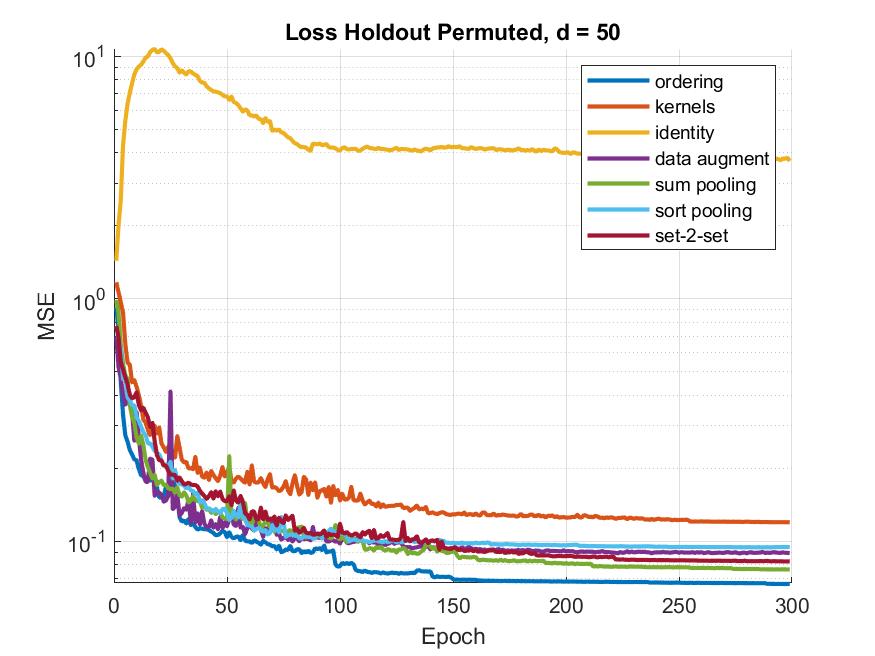}	
\includegraphics[width=0.49\linewidth]{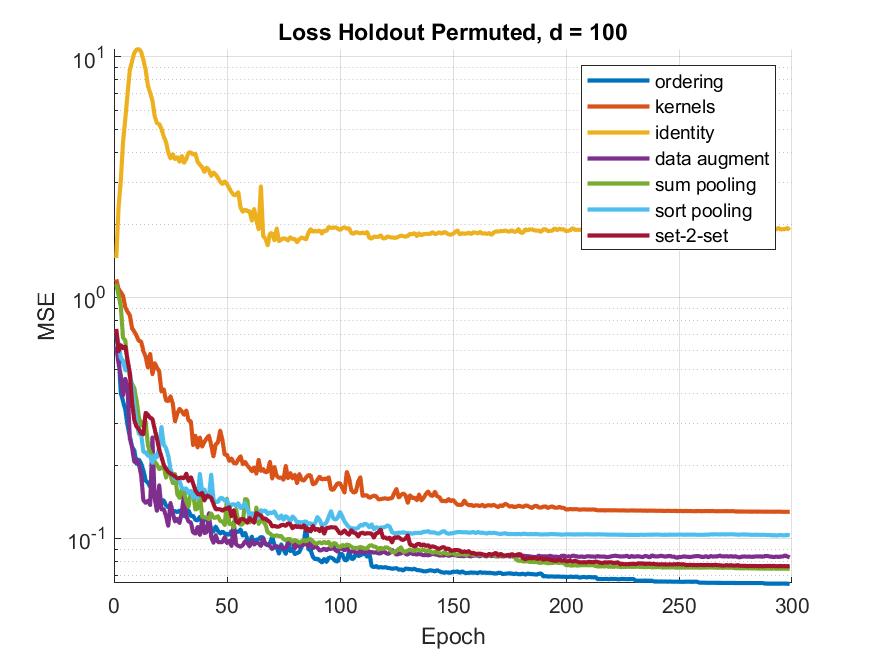}
\caption[caption me]{Mean-Square Error (MSE) Results for regression of the electron gap energy $\Delta\varepsilon=LUMO-HOMO$ (eV) on QM9 dataset using the seven algorithms; $d=50$ left column; $d=100$ right column.}
	\label{fig:a2}
\end{figure}

\begin{figure}[hp]
	\includegraphics[width=0.49\linewidth]{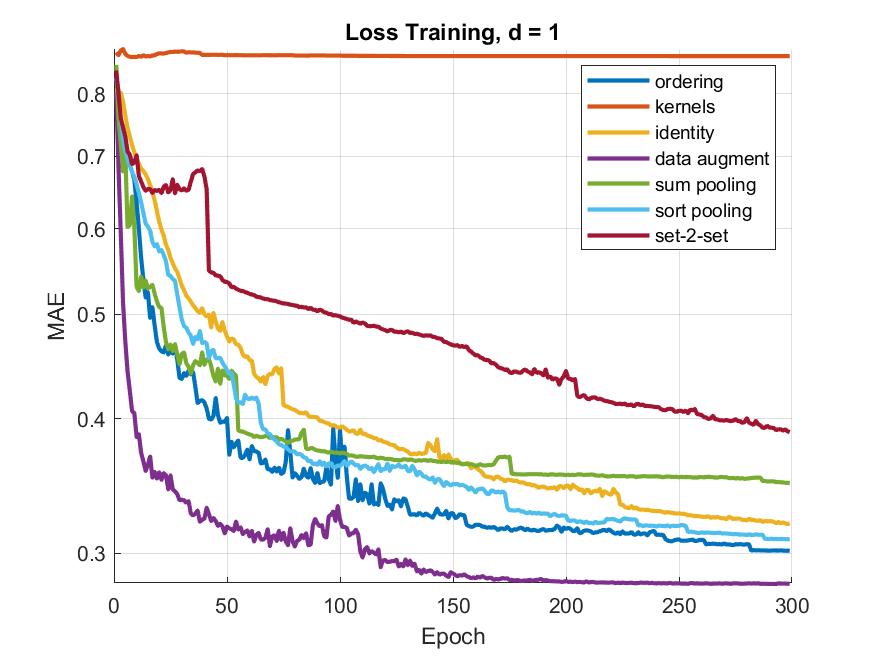}
	\includegraphics[width=0.49\linewidth]{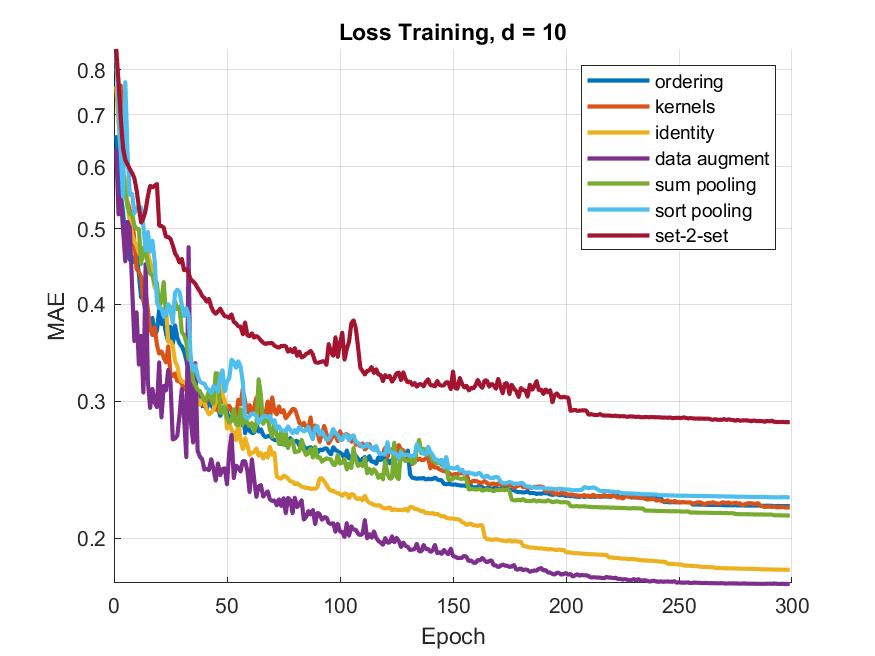} \\
	\includegraphics[width=0.49\linewidth]{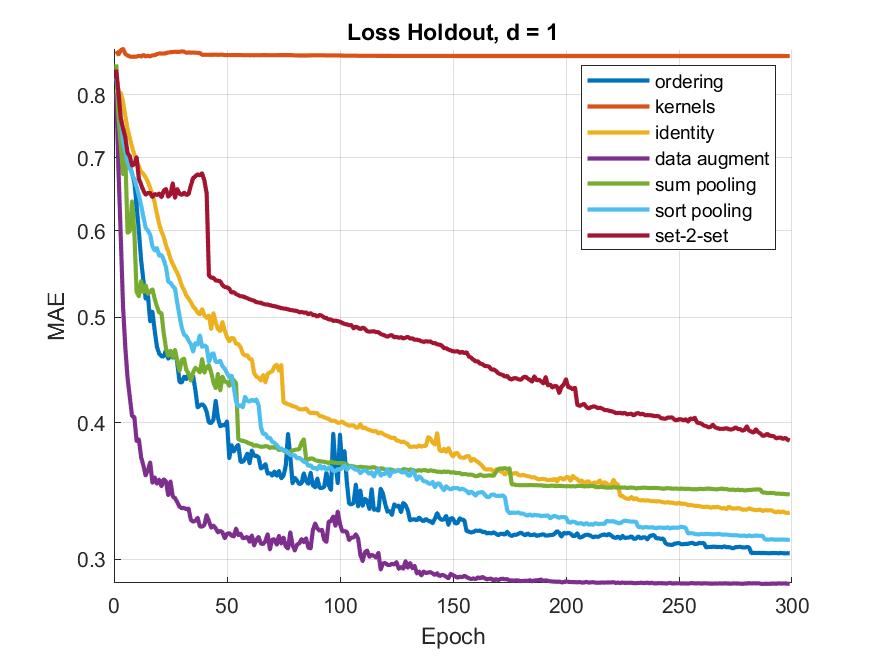}
	\includegraphics[width=0.49\linewidth]{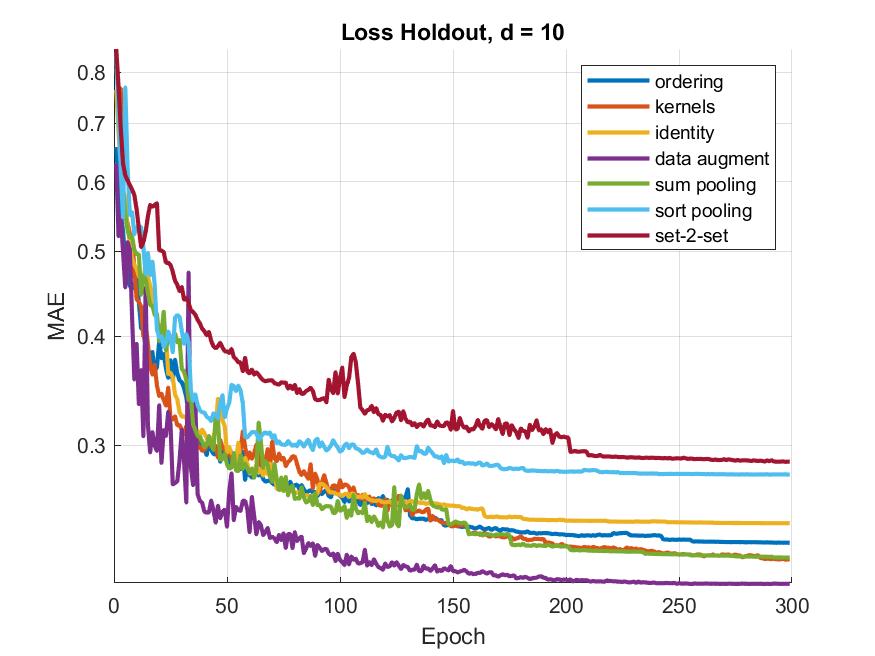} \\
	\includegraphics[width=0.49\linewidth]{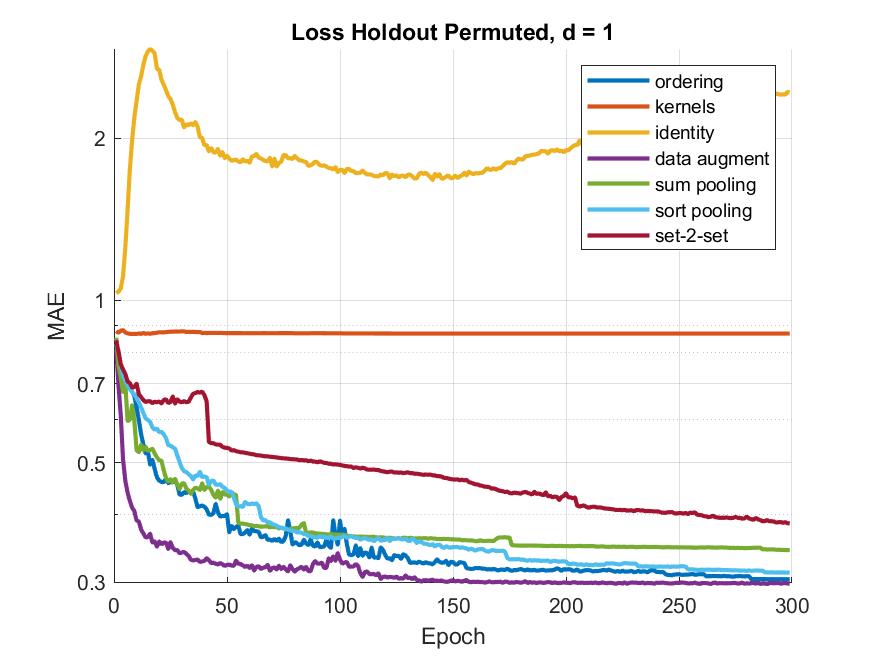}
	\includegraphics[width=0.49\linewidth]{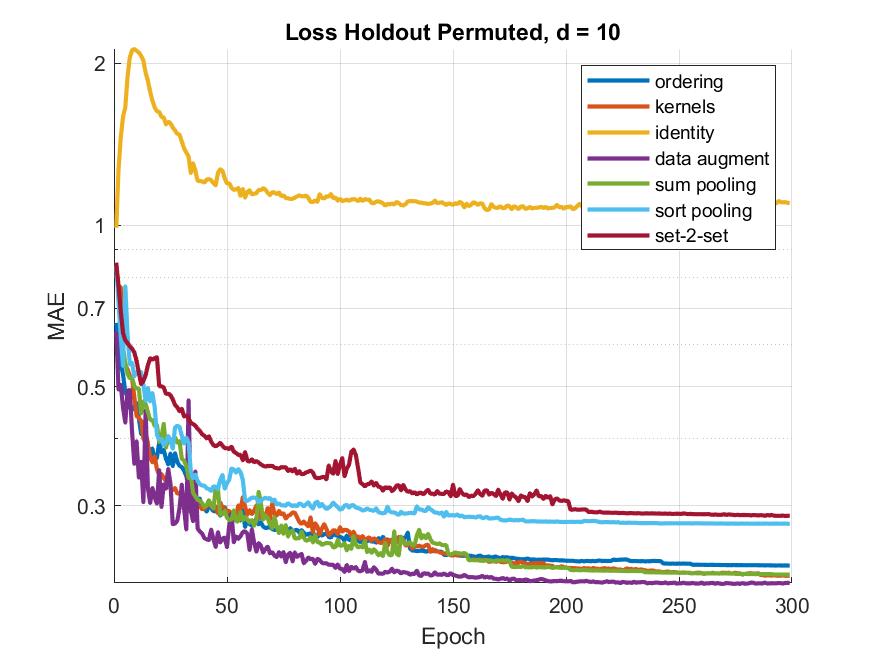}
	\caption[caption me]{Mean-Absolute Error (MAE) Results for regression of the electron gap energy $\Delta\varepsilon=LUMO-HOMO$ (eV) on QM9 dataset using the seven algorithms; $d=1$ left column; $d=10$ right column.}
	\label{fig:b1}
\end{figure}

\begin{figure}[hp]
	\includegraphics[width=0.49\linewidth]{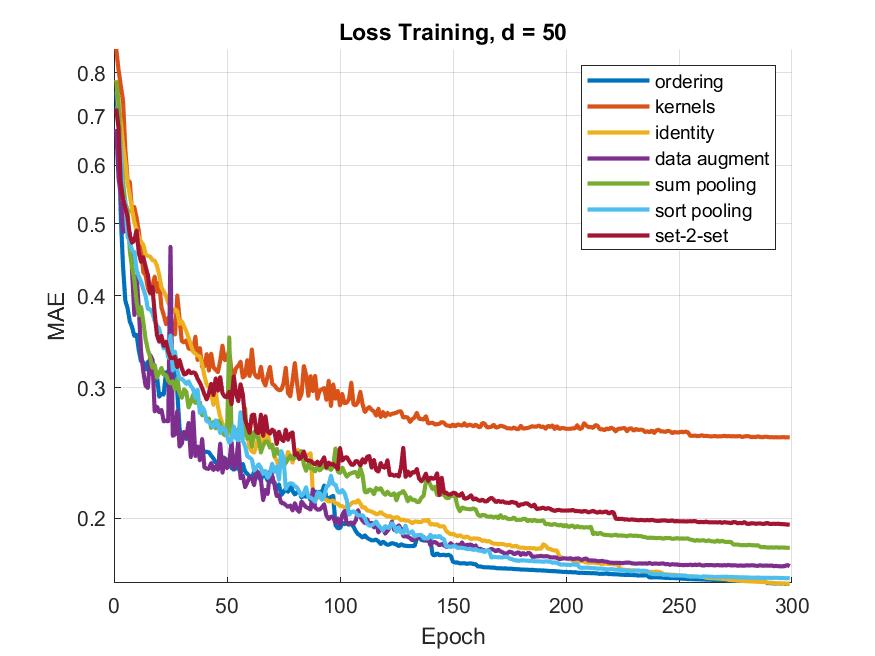}
	\includegraphics[width=0.49\linewidth]{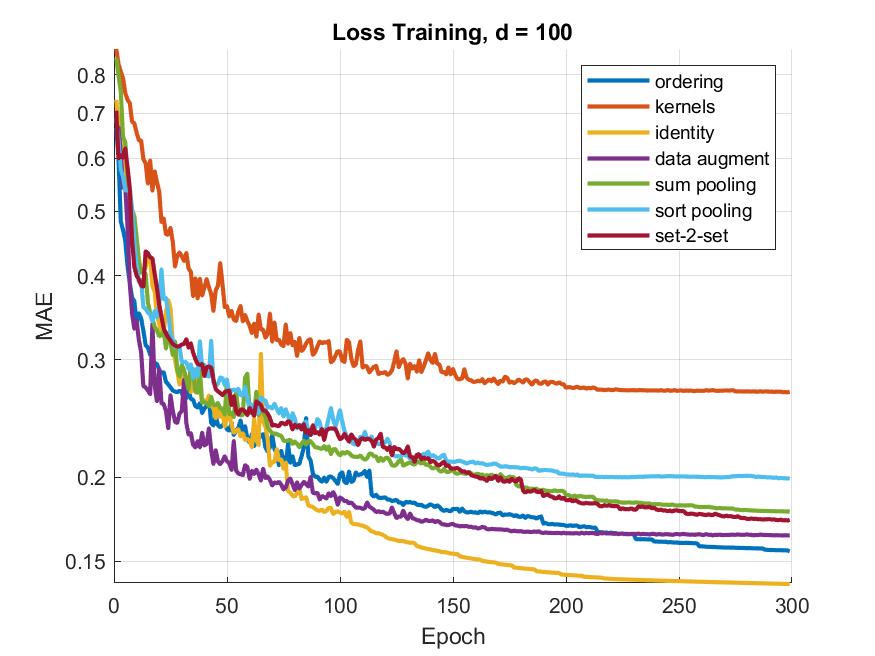} \\
	\includegraphics[width=0.49\linewidth]{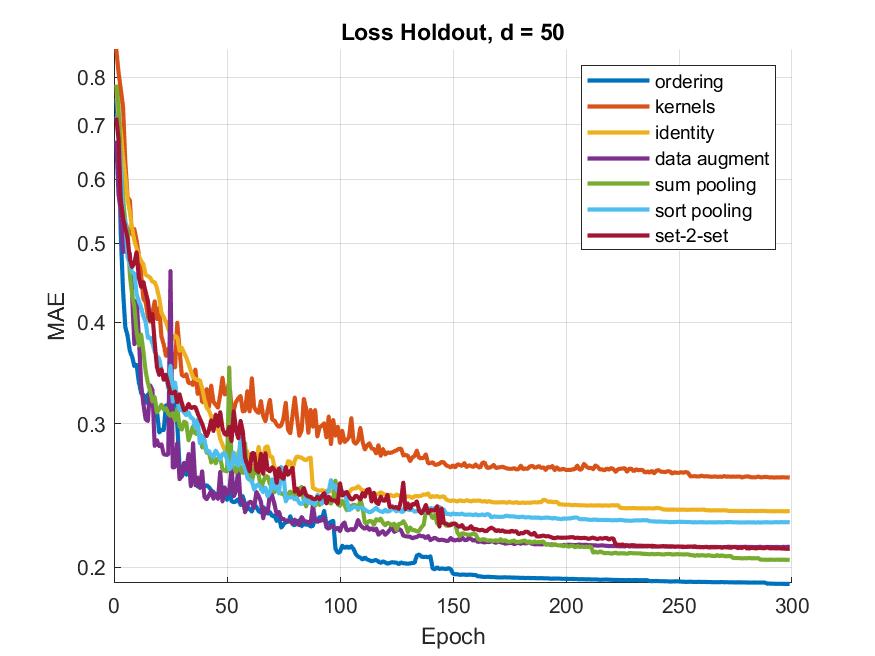}
	\includegraphics[width=0.49\linewidth]{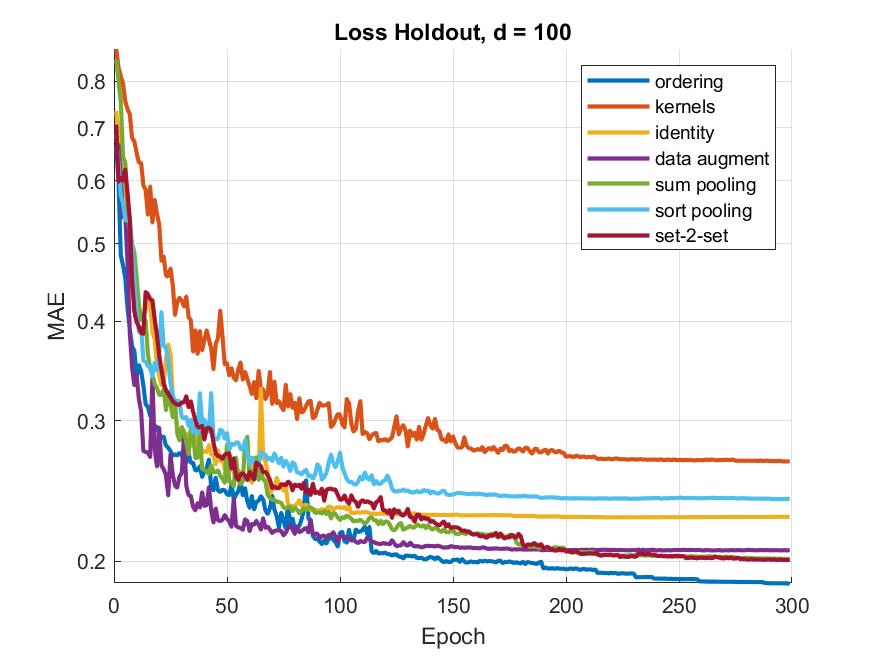} \\
	\includegraphics[width=0.49\linewidth]{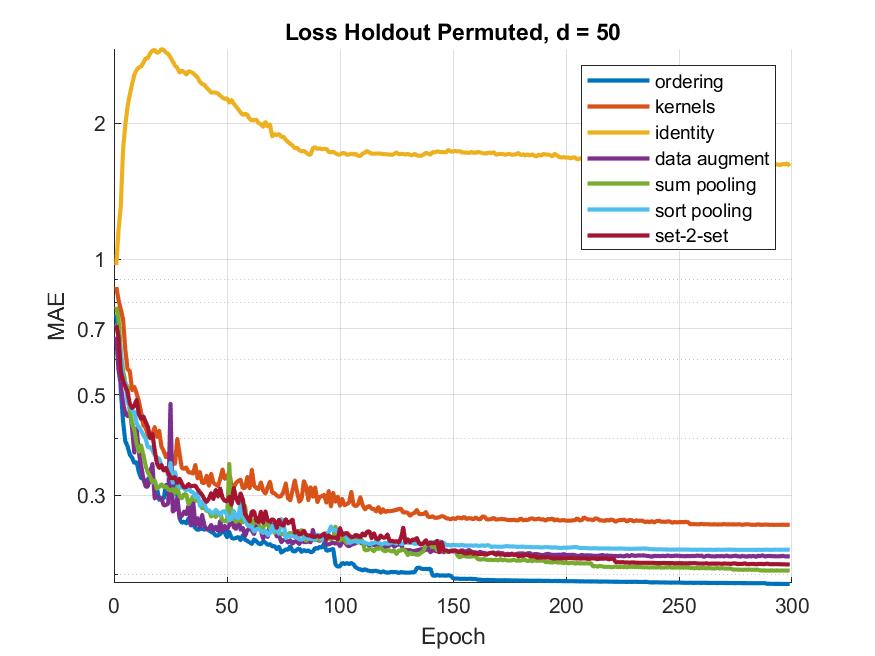}
	\includegraphics[width=0.49\linewidth]{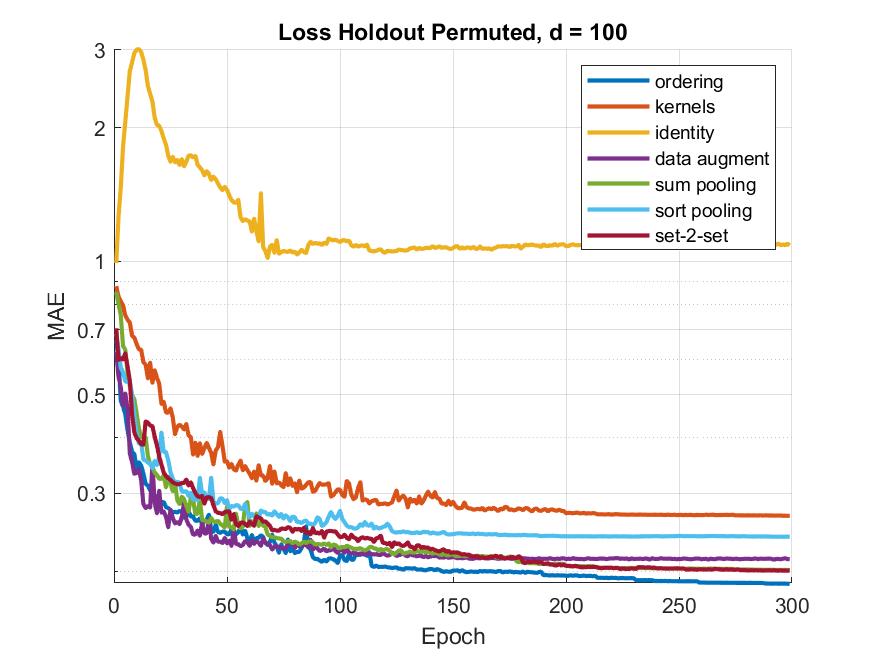}
	\caption[caption me]{Mean-Absolute Error (MAE) Results for regression of the electron gap energy $\Delta\varepsilon=LUMO-HOMO$ (eV) on QM9 dataset using the seven algorithms; $d=1$ left column; $d=10$ right column.}
	\label{fig:b2}
\end{figure}

\end{document}